\documentstyle[
11pt,twoside,amsfonts,amssymb,amsmath,amsthm,latexsym]{article}
{\tiny }
\newcommand {\rel} {{\mathbb R}}
\newcommand {\com} {{\mathbb C}}

\newcommand {\nat} {{\mathbb N}}

\newcommand {\sphere} {{\mathbb S}}
\newcommand {\Will} {{\mathcal{W} }}
\newcommand {\Wil} {{\mathcal{E} }}
\newcommand {\Area} {{\mathcal{A} }}
\newcommand {\Lift} {{\mathcal{L} }}
\newcommand {\NN} {{\mathcal{N} }}

\newcommand {\dom} {{\mathcal{D} }}
\newcommand {\OO} {{\mathcal{O} }}
\newcommand {\PP} {{\mathcal{P} }}
\newcommand {\Tan} {{\mathcal{T} }}
\newcommand {\F} {{\mathcal{F} }}
\newcommand {\R} {{\mathcal{R} }}
\newcommand {\CC} {{\mathcal{C} }}
\newcommand {\V} {{\mathcal{V} }}

\newcommand {\HH} {{\mathcal{H} }}
\newcommand {\Mill} {{\mathcal{M} }} 
\newcommand {\MR} {{\mathcal{MR} }}

\begin{document}
	\thanks{Institute: Mathematics Department, 
		Technion: Israel Institute of Technology,
		3200003 Haifa, Israel.\\
		Phone: 00972 549298486.
		E-mail: rubenj at technion.ac.il}
	
	\newtheorem{theorem}{Theorem}[section]
	\newtheorem{definition}{Definition}[section]
	\newtheorem{proposition}{Proposition}[section]
	\newtheorem{lemma}{Lemma}[section]
	\newtheorem{corollary}{Corollary}[section]
	\newtheorem{remark}{Remark}[section]
	\newtheorem{example}{Example}
	
	\author{Ruben Jakob}
	
	\title{Functional analytic properties and regularity of the 
	M\"obius-invariant Willmore flow in $\rel^n$}
	
	\maketitle
	
\begin{abstract}
In this article we continue the author's investigation of
the M\"obius-invariant Willmore flow
moving parametrizations of umbilic-free tori in $\rel^n$ 
and in the $n$-sphere $\sphere^n$.
In the main theorems of this article we prove basic 
properties of the evolution operator of the 
``DeTurck modification'' of the M\"obius-invariant Willmore flow and of its Fr\'echet derivative by means of a combination of the author's results about this topic with the theory of bounded $\HH_{\infty}$-calculus for linear elliptic operators 
due to Amann, Denk, Duong, Hieber, Pr\"uss and Simonett with 
Amann's and Lunardi's work on semigroups and interpolation theory. Precisely, we prove local real analyticity of the evolution operator $[F\mapsto \PP^*(\,\cdot\,,0,F)]$ of the ``DeTurck modification'' of the M\"obius-invariant Willmore flow in a small open ball in $W^{4-\frac{4}{p},p}(\Sigma,\rel^n)$, for any 
$p\in (3,\infty)$, about any fixed smooth parametrization 
$F_0:\Sigma \longrightarrow \rel^n$ of a compact and umbilic-free torus in $\rel^n$.
We prove moreover that the entire maximal flow line 
$\PP^*(\,\cdot\,,0,F_0)$, starting to move in a 
smooth and umbilic-free initial immersion $F_0$, is real analytic for positive times, and that therefore the Fr\'echet derivative $D_{F}\PP^*(\,\cdot\,,0,F_0)$ of the evolution operator in $F_0$ can be uniquely extended to a family of continuous linear operators $G^{F_0}(t_2,t_1)$ in $L^p(\Sigma,\rel^n)$, whose ranges are dense in $L^{p}(\Sigma,\rel^n)$, for every fixed pair of times 
$t_2\geq t_1$ within the interval of maximal existence $(0,T_{\textnormal{max}}(F_0))$.   
\end{abstract}

MSC-codes: 53C42,\,  35K46, \,35R01,\, 47B12,\, 58J35

\section{Introduction and state of the art}    \label{Intro}

In \cite{Jakob_Moebius_2016} the author has
proved short-time existence and uniqueness of smooth flow lines of the following non-linear, Willmore-type evolution equation:
\begin{equation}  \label{Moebius.flow}
	\partial_t f_t = -\frac{1}{2} \frac{1}{|A^0_{f_t}|^4} \,
	\Big{(} \triangle_{f_t}^{\perp} \vec H_{f_t} + Q(A^{0}_{f_t})(\vec H_{f_t}) \Big{)}
	\equiv -\frac{1}{|A^0_{f_t}|^4} \,\nabla_{L^2}\Will(f_t).
\end{equation}
Here, $\{f_t\}_{t \in \rel_+}$ denotes a differentiable family of $W^{4,2}$-immersions mapping some arbitrarily fixed compact smooth torus $\Sigma$ into $\rel^n$ or into $\sphere^n$, for $n\geq 3$, without any umbilic points. 
As already pointed out in the author's article \cite{Jakob_Moebius_2016}, the ``umbilic-free condition'' $|A^0_{f_t}|^2>0$ on $\Sigma$ implies that $\chi(\Sigma)=0$, where $\chi$ denotes the topological Euler-characteristic,
which forces the geometric flow (\ref{Moebius.flow}) to be only well-defined on differentiable families of sufficiently smooth umbilic-free tori, immersed into $\rel^n$ or $\sphere^n$, with $n \geq 3$. In flow equation (\ref{Moebius.flow}), $\Will$ denotes the Willmore-functional
\begin{equation} \label{Willmore.functional}
	\Will(f):= \int_{\Sigma} K^M_f + 
	\frac{1}{4} \, \mid \vec H_f \mid^2 \, d\mu_f,
\end{equation}
which can more generally be considered on $C^{2}$-immersions
$f:\Sigma \longrightarrow M$, mapping any closed smooth Riemannian orientable surface $\Sigma$ into an arbitrary 
smooth Riemannian manifold $M$, where $K^M_f(x)$ denotes the
sectional curvature of $M$ w.r.t. the ``immersed tangent plane'' $Df_x(T_x\Sigma)$ in $T_{f(x)}M$. The most prominent  cases are: $K_f\equiv 0$ if $M=\rel^n$ and $K_f\equiv 1$ if $M=\sphere^{n}$. Regarding the aims of this article, we will only have to consider the simplest case $M=\rel^n$, with 
$n\geq 3$. We endow the torus $\Sigma$ with the pullback $f^*g_{\textnormal{euc}}$ w.r.t. $f$ of the Euclidean metric 
of $\rel^n$, i.e. with coefficients 
$g_{ij}:=\langle \partial_i f, \partial_j f \rangle_{\rel^n}$, and we let $A_f$ denote the second fundamental form of the immersion $f:\Sigma \longrightarrow \rel^{n}$, defined on 
pairs of tangent vector fields $X,Y$ on $\Sigma$ by:
\begin{eqnarray}  \label{second.fundam.form}
	A_{f}(X,Y) := D_X(D_Y(f)) - P^{\textnormal{Tan}(f)}(D_X(D_Y(f)))
	\equiv (D_X(D_Y(f)))^{\perp_{f}},                                           
\end{eqnarray}
where $D_X(V)\lfloor_{x}$ denotes the classical derivative 
of a vector field $V:\Sigma \longrightarrow \rel^{n}$
in direction of the tangent vector field $X \in \Gamma(T\Sigma)$ at a point $x\in \Sigma$, and where
$$  
P^{\textnormal{Tan}(f)}:
\bigcup_{x \in \Sigma} \{x\} \times \rel^{n}
\longrightarrow \bigcup_{x \in \Sigma} \{x\} \times \textnormal{T}_{f(x)}(f(\Sigma))
$$
denotes the bundle morphism, which projects
the entire ambient space $\rel^n$ orthogonally into its subspaces $\textnormal{T}_{f(x)}(f(\Sigma))$ --
the tangent spaces of the immersion $f$ in the points $f(x)$ for every $x \in \Sigma$ -- and where $^{\perp_{f}}$ 
abbreviates the bundle morphism 
$\textnormal{Id}_{\rel^n} - P^{\textnormal{Tan}(f)}$. 
Furthermore, $A_{f}^0$ denotes the tracefree part of 
$A_{f}$, i.e.
$$
A_{f}^0(X,Y):= A_{f}(X,Y) - 
\frac{1}{2} \,g_f(X,Y)\, \vec H_{f}
$$
and 
$$
\vec H_{f}:=\textnormal{trace}(A_{f}) 
\equiv g^{ij}_f \,A_{f}(\partial_i,\partial_j) 
$$
(``Einstein's summation convention'') denotes the mean curvature vector of the immersion $f:\Sigma \longrightarrow \rel^n$. Finally, $Q(A_f)$ respectively $Q(A^0_f)$ operates 
on vector fields $\phi$ which are sections of the normal 
bundle of $f$, i.e. which ``are normal along $f$'', by assigning: 
$$ 
Q(A_f)(\phi):= g^{ij}_f \, g^{kl}_f\, A_f(\partial_i,\partial_k) \langle A_f(\partial_j,\partial_l),\phi \rangle_{\rel^n},
$$ 
which is by definition again a section of the normal bundle of $f$.  \\
Furthermore, Weiner computed in \cite{Weiner}, Section 2, that the first variation of the Willmore functional in a smooth immersion $f:\Sigma \longrightarrow \rel^n$, $n\geq 3$, in direction of a smooth section $\Psi$ of the normal bundle 
of $f$:
\begin{eqnarray} \label{first_variation}
	\delta \Will(f,\Psi):=
	\partial_t(\Will(f+t \,\Psi))\lfloor_{t=0}
	= \frac{1}{2} \int_{\Sigma}
	\langle \triangle^{\perp}_{f} \vec H_{f} + Q(A^{0}_{f})
	(\vec H_{f}), \Psi \rangle_{\rel^n}\, d\mu_{f} \quad \\
	=: \int_{\Sigma} \langle \nabla_{L^2}\Will(f),\Psi \rangle_{\rel^n} \, d\mu_{f}. \nonumber
\end{eqnarray}
For the reader's convenience we recall the following 
Proposition from \cite{Jakob_Moebius_2016}:
\begin{proposition} \label{trafo_heat_expression_2}
	\begin{itemize}
		Any differentiable family $\{f_t\}$ of $C^{4}$-immersions 
		$f_t:\Sigma \longrightarrow \rel^n$
		without umbilic points, i.e. with $\mid A^0_{f_t} \mid^2>0$ on $\Sigma$\, $\forall \, t\in [0,T)$, solves the flow equation  
		\begin{eqnarray} \label{invarianter_Waermefluss} 
		\partial_t f_t = -\frac{1}{2} \mid A^0_{f_t} \mid^{-4} \, \big{(} \triangle^{\perp}_{f_t} \vec H_{f_t} + 
		Q(A^{0}_{f_t}) (\vec H_{f_t}) \big{)} 
		\equiv - \mid A^0_{f_t} \mid^{-4} \,\nabla_{L^2}\Will(f_t)    \quad
		\end{eqnarray} 
		if and only if its composition $\Phi(f_t)$ with any applicable M\"obius-transformation $\Phi$ of $\rel^n$ solves the same flow equation, thus, if and only if  
		$$
		\partial_t (\Phi(f_t))= - \mid A^0_{\Phi(f_t)} \mid^{-4} \, \nabla_{L^2}\Will(\Phi(f_t))
		$$
		holds $\forall \, t\in [0,T)$ and for every $\Phi \in \textnormal{M\"ob}(\rel^n)$, for which $\Phi(f_t)$ is well-defined on $\Sigma \times [0,T)$. 
	\end{itemize}
\end{proposition} 
\noindent
Obviously, the above proposition suggests to term the geometric flow (\ref{invarianter_Waermefluss}) the ``M\"obius-invariant Willmore flow'' (MIWF). Actually, the author was able to show in the Appendix of \cite{Jakob_Moebius_2016},
that the expression on the right hand side of equation 
(\ref{invarianter_Waermefluss}) is the analytically simplest modification of the classical differential operator 
$[f \mapsto \nabla_{L^2}\Will(f)]$ from 
line (\ref{first_variation}), which yields a ``conformally invariant flow''. Looking at the right hand side of 
evolution equation (\ref{Moebius.flow}), 
one can easily see its ``degeneracy'' and thus 
imagine, that the corresponding flow produces many 
singular flow lines and thus several interesting technical challenges. Now, comparing the MIWF (\ref{Moebius.flow}) 
to the classical Willmore flow, given by the equation 
$$
\partial_t f_t = -\frac{1}{2} \,
\Big{(} \triangle_{f_t}^{\perp} \vec H_{f_t} 
+ Q(A^{0}_{f_t})(\vec H_{f_t}) \Big{)}
\equiv - \,\nabla_{L^2}\Will(f_t)
$$
the MIWF has the big advantage, that its flow lines can be 
conformally mapped - by stereographic projection -  
either from $\rel^n$ into $\sphere^n$ or from $\sphere^n$ into $\rel^n$, depending on the concrete objectives of the investigation. 
In particular the ``$3$-sphere'' is a simply-connected compact Lie-group - diffeomorphic to SU$(2)$ - can be interpreted as the set of quaternions of length $1$, is fibered by the Hopf-fibration $\pi:\sphere^3 \longrightarrow \sphere^2$ and contains the Clifford-torus - an embedded minimal surface in $\sphere^3$, which is the global minimizer of the 
Willmore functional (\ref{Willmore.functional}) 
among all compact surfaces of genus $\geq 1$ 
immersed into $\sphere^3$, on 
account of \cite{Marques.Neves.2014}, Theorem A.
This particular mathematical situation actually 
plays a key role in the proof of the first main 
theorem of the author's preprint \cite{Ruben.MIWF.IV},  
yielding a sufficient condition for global existence and full convergence of smooth flow lines of the MIWF in $\sphere^3$; 
see here point (2) of the list below. 
Moreover, even restricted to the rather narrow class 
of ``Hopf-tori'' in $\sphere^3$, the MIWF seems to 
develop ``singularities'' - see Definition \ref{Umbilic.free} below - since it reduces via the Hopf-fibration $\pi:\sphere^3 \longrightarrow \sphere^2$ to the ``degenerate'' parabolic flow: 
\begin{eqnarray}  \label{elastic.energy.flow}
	\partial_t \gamma_t =
	- \frac{1}{(\kappa_{\gamma_t}^2+1)^2} \,
	\Big{(} 2 \, \Big{(}\nabla^{\perp}_{\frac{\gamma_t'}
		{|\gamma_t'|}} \Big{)}^2(\vec{\kappa}_{\gamma_t})
	+ |\vec{\kappa}_{\gamma_t}|^2 \vec{\kappa}_{\gamma_t}
	+ \vec{\kappa}_{\gamma_t} \Big{)}                   
	\equiv -\frac{1}{(\kappa_{\gamma_t}^2+1)^2}  \,\nabla_{L^2}\Wil(\gamma_t)     
\end{eqnarray}
for smooth closed, regular curves 
$\gamma_t:\sphere^1\longrightarrow \sphere^2$, 
where ``$\vec{\kappa}_{\gamma_t}$'' denotes the curvature 
vectors of the curves $\gamma_t$, and 
``$\nabla_{L^2} \Wil$'' denotes the $L^2$-gradient of twice the classical elastic energy, i.e. of the functional
$$
\Wil(\gamma):=
\int_{\sphere^1} 1+|\vec{\kappa}_{\gamma}|^2 \, \,d\mu_{\gamma}.
$$
See Propositions 3.2 and 3.3 in \cite{Ruben.MIWF.II} for 
technical details of this particular ``dimension reduction''.
The author's most recent research in \cite{Ruben.MIWF.V} reveals, that this sort of ``degeneracy'' of the MIWF 
does actually not improve under the additional energy condition ``$\Will(F_0)<8\pi$'' on an initial parametrization $F_0$ of a smooth Hopf-torus $\pi^{-1}(\textnormal{trace}(\gamma_0))$ 
in $\sphere^3$. Hence, the MIWF presents us with 
many challenging questions and puts the mathematical 
machinery of modern Geometric Analysis to the test. 
\footnote{One should remark here, that the classical Willmore flow in $\sphere^3$ actually turns out to have only global smooth flow lines, which start moving in smooth parametrizations of arbitrary smooth Hopf-tori in $\sphere^3$, see Theorem 2.1 in \cite{Ruben.MIWF.II}.}  
On the other hand, we can also use ``$\rel^n$'' as 
ambient space for flow lines of the MIWF, and this allows us 
to apply heavy machinery from both linear and non-linear ``Functional Analysis'', 
in order to investigate the evolution operator of the MIWF 
- as an operator between appropriately chosen Banach spaces - 
via its linearization, and this is actually the main 
technical focus of this article and also of 
the author's paper \cite{Jakob_Moebius_2016}.    
So far, the author has proved four basic results 
about the MIWF:
\begin{itemize} 
	\item[1)] Short-time existence and uniqueness of smooth flow lines $\{\PP(t,0,F_0)\}_{t \geq 0}$ of the M\"obius-invariant Willmore flow (\ref{Moebius.flow}), which start moving in $C^{\infty}$-smooth and umbilic-free initial immersions $F_0:\Sigma \longrightarrow \rel^n$, see Theorem 1 in \cite{Jakob_Moebius_2016} for the precise statement.   
	\item[2)] Global existence and full $C^k(\Sigma)$-convergence 
	of any flow line $\{\PP(t,0,F_0)\}_{t \geq 0}$ of the MIWF (\ref{Moebius.flow}) to a smooth parametrization of a Clifford-torus in $\sphere^3$, provided $\{\PP(t,0,F_0)\}_{t \geq 0}$ starts moving in a $C^{\infty}$-smooth immersion 
	$F_0:\Sigma \longrightarrow \sphere^3$ which is sufficiently close to a fixed smooth parametrization of 
	a Clifford-torus in a $h^{2+\beta}(\Sigma,\rel^4)$-norm. Here, ``$h^{2+\beta}(\Sigma,\rel)$'' 
	denotes the ``little H\"older space'' of differentiation 
	order $2+\beta$. See Theorem 1.1 in \cite{Ruben.MIWF.IV} for the precise statement.   
	\item[3)] Stability of the above ``full convergence property''
	of flow lines of the MIWF (\ref{Moebius.flow}) w.r.t. 
	perturbations of their umbilic-free initial immersions 
	$F:\Sigma \longrightarrow \sphere^3$   
	in any $C^{4,\gamma}(\Sigma,\rel^4)$-norm. 
	See here Theorem 1.2 in \cite{Ruben.MIWF.IV}.
	\item[4)] Global existence and full $C^{k-1,\alpha}(\Sigma)$-convergence 
	of any flow line $\{\PP(t,0,F_0)\}_{t \geq 0}$ of the MIWF (\ref{Moebius.flow}) into some $C^k$-local minimizer of the Willmore functional, which starts moving 
	in a $C^{\infty}$-smooth immersion 
	$f_0:\Sigma \longrightarrow \rel^3$ that is 
	sufficiently close to a fixed $C^k$-local minimizer 
	of the Willmore functional in a  
	$C^{k,\alpha}(\Sigma,\rel^3)$-norm.
	See Theorem 1.3 in \cite{Ruben.MIWF.IV} for the precise statement. 
\end{itemize}
In this article, we will combine the theory 
of bounded $\HH_{\infty}$-calculus for linear elliptic differential operators of higher order due to Amann, Denk, Duong, Hieber, Pr\"uss and Simonett in \cite{Amann.1995}, \cite{Amann.Hieber.Simonett}, \cite{Denk.Hieber.Pruess.1} and 
\cite{Duong.Simonett.1997} with Amann's and Lunardi's 
work on non-autonomously generated
``semigroups'' and ``interpolation theory'' in 
\cite{Amann.2}, \cite{Lunardi} and \cite{Lunardi.interpolation}, 
and with modern regularity criteria due to Shao and Simonett \cite{Shao.2013}, \cite{Shao.2015}, \cite{Shao.Simonett}, 
in order to prove below in Theorems \ref{real.analytic.flow},
\ref{Frechensbergo} and \ref{Functional.prop.Frechet.derivative} basic properties of the evolution operator $\{\PP^*(t,0,F_0)\}_{t \geq 0}$ of the ``DeTurck modification'' (\ref{de_Turck_equation}) of the MIWF (\ref{Moebius.flow}) in $\rel^n$ and of its ``linearization''. 
Just as in the author's first paper \cite{Jakob_Moebius_2016} 
about the MIWF, a certain adaption of the ``DeTurck trick'' will also be used in this paper, in order to be able to apply 
well developed theories dealing with linear respectively 
quasilinear parabolic differential equations.
Moreover we should also mention here, that certain 
variants of Proposition \ref{A.priori.estimate} 
and of Theorem \ref{Psi.of.class.C_1} below in Section \ref{Preparation} play key r\^oles in the proof of the above mentioned Theorem 1.3 of the author's preprint \cite{Ruben.MIWF.IV} about the MIWF; 
see point (4) in the list above.   \\
Returning to the above mentioned open question 
about ``singularities'' of the MIWF, we should also note here, 
that the statements of the five main theorems of this 
article hold mutatis mutandis also for the classical 
Willmore flow in $\rel^n$. Hence modulo obvious changes,
Theorems \ref{real.analytic.flow}, \ref{Frechensbergo} and \ref{Functional.prop.Frechet.derivative} of this article 
can be formulated for the classical Willmore flow in 
$\rel^n$ and can then be combined with the main theorem of the author's article \cite{Ruben.MIWF.II} about flow lines of the classical Willmore flow in $\rel^4$, 
whose initial immersions parametrize smooth 
Hopf-tori in $\sphere^3$, particularly with the aim 
to follow the lines of Andrews' article \cite{Andrews.2002} 
and to prove, that generically every such flow line 
of the classical Willmore flow converges fully and smoothly 
to some smooth parametrization of the Clifford-torus in $\sphere^3$ - modulo some M\"obius-transformation of 
$\sphere^3$ - without imposing any initial Willmore 
energy condition. Actually, motivated by Andrews' article \cite{Andrews.2002}, Theorems 
\ref{real.analytic.flow}, \ref{Frechensbergo} 
and \ref{Functional.prop.Frechet.derivative} - 
particularly the final part of Theorem \ref{Functional.prop.Frechet.derivative} - have been  
discovered and proved by the author exactly with the aim, 
to control - with ``relatively high precision'' - 
the deviation of flow lines 
$\{\PP^*(t,0,F)\}_{t\geq 0}$ in 
(\ref{solution.operator.1}) of the 
``DeTurck-modified'' flow equation (\ref{de_Turck_equation}) 
w.r.t. perturbations of their initial immersions 
$F$ in small open balls $B_{\rho}^{4,p}(F_0)$ in $W^{4,p}(\Sigma,\rel^n)$ at arbitrarily large times $t>>1$, provided every flow line 
$\{\PP^*(t,0,F)\}_{t\geq 0}$ starting in 
$F \in B_{\rho}^{4,p}(F_0)\cap C^{\infty}(\Sigma,\rel^n)$ 
with $\rho>0$ sufficiently small, 
certainly exists globally, i.e. for all $t \geq 0$.       
It is quite a remarkable consequence of this article, 
that the functional analytic statements of Theorems \ref{Frechensbergo} and \ref{Functional.prop.Frechet.derivative} turn out to hold along flow lines of the ``rather degenerate'' MIWF, and not only along flow lines of the classical Willmore flow in $\rel^n$, provided we choose $p>3$ in those two theorems.  \\
This paper is organized as follows:
In Section \ref{Preparation} we firstly introduce 
some basic notation and prove fundamental properties of the
non-linear differential operator (\ref{D_F_0}) corresponding 
to the ``DeTurck modification'' (\ref{de_Turck_equation}) of the MIWF (\ref{Moebius.flow}) and of its linearization in the 
setting of parabolic $L^p$-spaces. Moreover, in Section \ref{Preparation} we also prove a new short-time existence 
result for evolution equation (\ref{de_Turck_equation}) and
real analyticity - both in space and in time - of maximal flow lines to evolution equation (\ref{de_Turck_equation}) starting 
in smooth - but not necessarily real analytic - and 
umbilic-free initial immersions of $\Sigma$ into $\rel^n$. 
To this end we follow the strategies of the papers \cite{Escher.Pruess.Simonett}, \cite{Pruess.Simonett.2004},
\cite{Koehne.Pruess.Wilke.2010}, \cite{Shao.2013}, \cite{Shao.2015} and \cite{Shao.Simonett}, i.e. we use 
the technique by Escher, Pr\"uss, Shao and 
Simonett of parameter-dependent diffeomorphisms on 
closed manifolds, maximal regularity in both 
$L^p$- and $h^{\beta}$-spaces and the 
``Implicit Function Theorem'' for real analytic, 
non-linear operators. 
In Section \ref{Linearization} we prove the main results 
of this paper, whose proofs heavily rely on the 
local real analyticity of the evolution operator of 
evolution equation (\ref{de_Turck_equation}) - see 
the first part of Theorem \ref{Frechensbergo} - and the  
above mentioned real analyticity of smooth flow lines 
of evolution equation (\ref{de_Turck_equation}). 
In particular, this regularity result is a very effective tool, 
in order to prove the two final statements of this paper, Theorem \ref{Functional.prop.Frechet.derivative} (ii) and (iii), yielding the fact that the Fr\'echet derivative of the evolution operator of evolution equation (\ref{de_Turck_equation}) in any smooth and umbilic-free immersion $F_0$ can be interpreted as a 
$2$-parameter family ``$G^{F_0}(t,s)$'' 
of linear continuous operators in $L^p(\Sigma,\rel^n)$ 
mapping $W^{4,p}(\Sigma,\rel^n)$ into 
itself and having \underline{dense range in 
$L^{p}(\Sigma,\rel^n)$}, for any choice of 
the two parameters $s\leq t$ within the open interval of 
maximal existence $(0,T_{\textnormal{max}}(F_0))$ of 
the considered flow line $\PP^*(\,\cdot\,,0,F_0)$ 
starting in $F_0$ at time $s=0$, and for any fixed $p>3$.

\section{Preparation for the proofs of Theorems \ref{Frechensbergo} and \ref{Functional.prop.Frechet.derivative}}  \label{Preparation}

We follow the definition of Besov- and Sobolev-Slobodeckij-spaces $W^{s,p}(\rel^m_+)$ on half spaces $\rel^m_+$ in Sections 2.9.1 and 2.9.3 in \cite{Triebel} and define here similarly $L^p$- and $W^{s,p}$-functions 
on smooth domains respectively on compact, closed smooth 
surfaces, having values in $\rel^n$.
\begin{definition}
	Let $p \in (1,\infty)$ and $s\geq 2$ be two real numbers, 
	$n\in \nat$, and let $\Omega$ be a bounded domain in $\rel^2$ with smooth boundary and $M$ a compact, closed smooth surface, equipped with an atlas $\Area:=\{(\Omega_j,\psi_j)\}_{j=1,\ldots,N}$ 
	of coordinate patches $\Omega_j \subset M$ and smooth charts 
	$\psi_j:\Omega_j \stackrel{\cong}\longrightarrow B^2_1(0)$.   
	\begin{itemize} 
		\item [1)] We define the space of $L^p$-functions 
		on the compact surface $M$ by 
		$$ 
		L^p(M,\rel^n):=\{\,u:M \longrightarrow \rel^n\, |\, u \circ \psi_j^{-1} \in L^p(B^2_1(0),\rel^n) \,\, \textnormal{for} \,\,j=1,\ldots,N \, \},
		$$
		and we set 
		$$
		\parallel u \parallel_{L^p(M,\rel^n)}
		:= \max_{j=1,\ldots,N} \parallel u \circ \psi_j^{-1} \parallel_{L^p(B^2_1(0),\rel^n)}.
		$$   
		\item [2)] We define the Sobolev-Slobodeckij-space $W^{s,p}(\Omega,\rel^n)$ on the bounded domain $\Omega \subset \rel^2$ by  
		\begin{eqnarray*}
			W^{s,p}(\Omega,\rel^n):=\{ u \in L^p(\Omega,\rel^n) |\, 
			\exists \,g =(g_1,\ldots,g_n)  
			\in W^{s,p}(\rel^2,\rel^n) \\ \,\,\textnormal{s.t.} \,\,\, u(x)=g(x) \,\,
			\textnormal{for} \,\, \HH^2-\textnormal{almost every} \, x \in \Omega \, \}.\nonumber
		\end{eqnarray*}  
		Here, the Sobolev-Slobodeckij-space ``$W^{s,p}(\rel^2,\rel)$'' is defined as in Section 2.3.1 in \cite{Triebel}. 
		\item [3)] We define the Sobolev-Slobodeckij-space $W^{s,p}(M,\rel^n)$ on the compact surface $M$ by 
		$$ 
		W^{s,p}(M,\rel^n) := \{u \in L^p(M,\rel^n) \,| \,
		u \circ \psi_j^{-1} \in W^{s,p}(B^2_1(0),\rel^n) \,\, \textnormal{for}\,\,j=1,\ldots, N \}
		$$ 
		and its norm by
		$$
		\parallel u \parallel_{W^{s,p}(M,\rel^n)}:= 
		\max_{j=1,\ldots,N} \parallel  u \circ \psi_j^{-1} \parallel_{W^{s,p}(B^2_1(0),\rel^n)}
		$$   
		for any $u\in W^{s,p}(M,\rel^n)$.
	\end{itemize}
\qed 
\end{definition}  
\noindent
Moreover, we will need:
\begin{definition}  \label{Umbilic.free}
Let $\Sigma$ be a smooth compact torus and 
$n \geq 3$ an integer. 
\begin{itemize} 
\item[a)] We denote by $\textnormal{Imm}_{\textnormal{uf}}(\Sigma,\rel^n)$
the subset of $C^2(\Sigma,\rel^n)$ consisting of 
umbilic-free immersions, i.e.:
		$$ 	
		\textnormal{Imm}_{\textnormal{uf}}(\Sigma,\rel^n):=\{ f \in C^{2}(\Sigma,\rel^n)\,|\, f \,\,\textnormal{is an immersion with} \, \mid A^0_{f} \mid^{2}>0 \,\, \textnormal{on} \,\,\Sigma \,\}.
		$$
		\item[b)] A ``flow line'' of the MIWF (\ref{Moebius.flow}) in $\rel^n$ is a smooth family $\{f_t\}_{t\in [0,T)}$ of smooth immersions of $\Sigma$ into $\rel^n$, such that the resulting smooth function 
		$f:\Sigma \times [0,T)\longrightarrow \rel^n$ satisfies equation (\ref{Moebius.flow}) classically on $\Sigma \times [0,T)$, i.e. such that 
		$$
		\partial_t f_t(x) = -\frac{1}{2} \frac{1}{|A^0_{f_t}(x)|^4} \,
		\Big{(} \triangle_{f_t}^{\perp} \vec H_{f_t}(x) 
		+ Q(A^{0}_{f_t})(\vec H_{f_t})(x) \Big{)}
		$$
		holds pointwise in every $(x,t)\in \Sigma \times [0,T)$. 
		\item[c)] We use exactly the same terminology as in part (b) for smooth solutions $\{f_t\}_{t\in [0,T)}$ of the ``relaxed MIWF-equation'' in $\rel^n$, which is: 
		$$
		\big{(}\partial_t f_t(x)\big{)}^{\perp_{f_t}} 
		= -\frac{1}{2} \frac{1}{|A^0_{f_t}(x)|^4} \,
		\Big{(} \triangle_{f_t}^{\perp} \vec H_{f_t}(x) + Q(A^{0}_{f_t})(\vec H_{f_t})(x) \Big{)}
		$$ 
		for $(x,t)\in \Sigma \times [0,T)$.
		Here, $^{\perp_{f_t}}$ abbreviates the projection of 
		the velocity vector $\partial_t f_t(x)$ into the normal space of the immersion $f_t$ within $\rel^n$, for every fixed $x\in \Sigma$, as in (\ref{second.fundam.form}). 	
		\item[d)] Let $F_0:\Sigma \longrightarrow \rel^n$ be a smooth and umbilic-free immersion and $\{F_t\}_{t\in [0,T)}$ a flow line of the MIWF starting in $F_0$. We call $[0,T)$ the ``interval of maximal existence'' of the MIWF starting in $F_0$, if either $T=\infty$, 
		or if there holds $T<\infty$ and there is not: an  
		$\varepsilon>0$ and a smooth solution 
		$\{\tilde F_t\}_{t\in [0,T+\varepsilon)}$ 
		of the MIWF with $\tilde F_t= F_t$ on $\Sigma$ for 
		$t\in [0,T)$. In both cases the element $T \in \rel \cup \{\infty\}$ is uniquely determined by the initial immersion $F_0$, and we call it the ``maximal time of existence'' of the MIWF starting in $F_0$, in symbols ``$T_{\textnormal{max}}(F_0)$'',  respectively we call $\{F_t\}_{t\in [0,T_{\textnormal{max}}(F_0))}$ the 
		``maximal solution'' of evolution equation (\ref{Moebius.flow}) starting in $F_0$.
		\item[e)] If $T_{\textnormal{max}}(F_0)$ is finite, 
		then we also call $T_{\textnormal{max}}(F_0)$   
		``the singular time'' of the flow line $\{F_t\}$ 
		of the MIWF starting in $F_0$. In this case we also say, that the respective flow line of the MIWF ``forms a singularity'' as $t\nearrow T_{\textnormal{max}}(F_0)$.
	\end{itemize}
\qed 
\end{definition}
\noindent
We should note here, that Definition \ref{Umbilic.free}, 
(b)--(e), makes sense on account of ``short-time 
existence and uniqueness'' of flow lines of the MIWF, 
guaranteed by Theorem 1 in \cite{Jakob_Moebius_2016}.
Now, for any fixed immersion 
$G:\Sigma \longrightarrow \rel^n$ of class $C^{2}$ and a smooth chart $\psi$ of an arbitrary coordinate patch $\Sigma'$ of a fixed smooth compact torus $\Sigma$, we will denote throughout this article the resulting partial derivatives on $\Sigma'$ by $\partial_i$, $i=1,2$, the coefficients 
$g_{ij}:=\langle \partial_{i}G, \partial_{j}G \rangle_{\rel^n}$ of the first fundamental form of $G$ w.r.t. $\psi$ and the associated Christoffel-symbols 
$(\Gamma_G)^m_{kl}:=g^{mj}_G \,\langle \partial_{kl} G,\partial_jG \rangle_{\rel^n}$ of $(\Sigma',G^*(g_{\textnormal{eu}}))$. Moreover, we define the first (covariant) derivatives by $\nabla^G_i(v):=\nabla^G_{\partial_i}(v):=\partial_i(v)$, $i=1,2$, and the second covariant derivatives by
\begin{equation} \label{double_nabla} 
	\nabla_{kl}^G(v)\equiv \nabla_k^{G} \nabla_l^{G}(v)
	:=\partial_{kl}v - (\Gamma_G)^m_{kl} \,\partial_mv
\end{equation} 
of any function $v \in C^2(\Sigma,\rel)$. 
Moreover, for any vector field $V \in C^2(\Sigma,\rel^n)$ we define the projections of its first derivatives onto the normal bundle of the immersed torus $G(\Sigma)$ by
$$
\nabla_{i}^{\perp_G}(V) \equiv (\nabla^G_i(V))^{\perp_G}:=\nabla^G_i(V)- P^{\textnormal{Tan}(G)}(\nabla^G_i(V))
$$ 
and the ``normal second covariant derivatives'' of 
$V$ w.r.t. the immersion $G$ by
$$
\nabla_{k}^{\perp_G} \nabla_{l}^{\perp_G}(V):= 
\nabla_k^{\perp_G} (\nabla_l^{\perp_G}(V))-(\Gamma_G)^m_{kl} \,\nabla_{m}^{\perp_G}(V). 
$$
Using these terms, we define the Beltrami-Laplacian 
w.r.t. $G$ by $\triangle_G(V):= g^{kl}_G \nabla_{kl}^G(V)$, its projection $(\triangle_G V)^{\perp_G}:=\big{(} g^{kl}_G \nabla_{k}^{G} \nabla^G_{l}(V)\big{)}^{\perp_G}$ onto the normal bundle of the surface 
$G(\Sigma)$ and the ``normal Beltrami-Laplacian'' by 
$\triangle_G^{\perp_G}(V):=g^{kl}_G \nabla_{k}^{\perp_G} \nabla_{l}^{\perp_G}(V)$. 
We shall note here that equations (\ref{second.fundam.form})
and (\ref{double_nabla}) together imply: 
\begin{equation}  \label{second.fund.i.j}
	(A_G)_{ij}=A_G(\partial_i,\partial_j)=\partial_{ij} G - (\Gamma_G)^m_{ij} \,\partial_mG
	= \nabla_i^{G} \nabla_j^{G}(G),
\end{equation} 
which shows that the second fundamental form $A_G$ is a covariant tensor field of degree $2$ and that there holds:
\begin{equation}  \label{mean_curvat_laplacian}  
	\vec H_G = g^{ij}_G\, (A_G)_{ij} = g^{ij}_G \,\nabla_i^{G} \nabla_j^{G}(G) = \triangle_G(G)
\end{equation} 
for the mean curvature of the immersion $G$.
Now we recall from \cite{Jakob_Moebius_2016}, that the main problem about equation (\ref{invarianter_Waermefluss}) is its non-parabolicity. We have 
\begin{equation}    \label{triangle.perp} 
	\triangle^{\perp}_{f} \vec H_{f} + Q(A^{0}_{f})(\vec H_{f})
	=(\triangle_{f} \vec H_{f})^{\perp_f} + 
	2\,Q(A_f)(\vec H_f) 
	- \frac{1}{2} \mid \vec H_f \mid^2\,\vec H_f  
\end{equation} 
and by (\ref{mean_curvat_laplacian}): 
\begin{eqnarray}   \label{leading_term_of_invariant}
	(\triangle_{f} \vec H_{f})^{\perp_f}
	=  g^{ij}_{f} \, g^{kl}_{f} 
	\, \nabla_i^{f} \nabla_j^{f} \nabla_k^{f} \nabla_l^{f}(f)- 
	g^{ij}_{f} \, g^{kl}_{f} \, 
	\langle \nabla_i^{f} \nabla_j^{f} \nabla_k^{f} \nabla_l^{f}(f), 
	\partial_m f \rangle\, g^{mr}_{f} \, \partial_r(f) 
\end{eqnarray}
for any $W^{4,2}$-immersion $f:\Sigma \longrightarrow \rel^n$, which shows that the leading term of the right-hand side of (\ref{invarianter_Waermefluss}),
i.e. of $[\{f_t\} \mapsto \mid A^0_{f_t} \mid^{-4} \,\nabla_{L^2} \Will(f_t)]$, is not uniformly elliptic (of fourth order), even if $\mid A^0_{f_t} \mid^2$ should stay positive on the torus $\Sigma$ for all times $t\in [0,T]$. In order to overcome this unpleasant obstruction we have applied in \cite{Jakob_Moebius_2016} the ``DeTurck trick'' in the following way:
We fix some $C^{\infty}$-smooth immersion $F_0$ of $\Sigma$ 
into $\rel^n$ arbitrarily and compute:
\begin{eqnarray*}
	\nabla_k^{f} \nabla_l^{f}(f)
	= \big{(}\partial_{kl}f - (\Gamma_{F_0})^m_{kl} \partial_m(f)\big{)}
	+ C^m_{kl}(f,F_0) \, \partial_m(f)                \\ 
	= \nabla_k^{F_0} \nabla_l^{F_0}(f) + C^m_{kl}(f,F_0)\, \partial_m(f)
\end{eqnarray*}
for 
$$
C^m_{kl}(f,F_0) := \big{(}(\Gamma_{F_0})^m_{kl} - (\Gamma_f)^m_{kl}\big{)} 
\qquad \textnormal{on}  \quad \Sigma'.
$$
It is important to note here, that the difference 
$(\Gamma_{F_0})^m_{kl} - (\Gamma_f)^m_{kl}$ is a tensor field of $3$rd degree and that therefore the difference 
\begin{eqnarray}    \label{difference.Christoffel} 
	g^{ij}_{f} \, g^{kl}_{f} 
	\, \Big{(} \nabla_i^{f} \nabla_j^{f} \nabla_k^{f} \nabla_l^{f}(f) 
	- \nabla_i^{f} \nabla_j^{f} \big{(}C^m_{kl}(f,F_0) \, \partial_m(f)\big{)} \Big{)}      
	= g^{ij}_{f} \, g^{kl}_{f}\,
	\, \nabla_i^{f} \nabla_j^{f} \nabla_k^{F_0} \nabla_l^{F_0}(f) 
\end{eqnarray}
is a ``scalar'', i.e. does not depend on the choice of the local chart $\psi$,
and thus yields a globally well-defined differential operator of fourth order, for $W^{4,2}$-immersions $f:\Sigma \longrightarrow \rel^n$ again.
Moreover, one can easily verify by (\ref{double_nabla}) and by the derivation formulae in (\ref{derivation_formulae}) below, that $f \mapsto  g^{ij}_{f} \, g^{kl}_{f}\,
\, \nabla_i^{f} \nabla_j^{f} \nabla_k^{F_0} \nabla_l^{F_0}(f)$
is a non-linear operator of fourth order whose leading term is 
$g^{ij}_{f} \, g^{kl}_{f} \nabla_i^{F_0} \nabla_j^{F_0} \nabla_k^{F_0} \nabla_l^{F_0}(f)$,
which is locally elliptic in the sense of Section 1 
in \cite{Mantegazza}. Since the highest order term 
of the above expression 
$\nabla_i^{f} \nabla_j^{f}(C^m_{kl}(f,F_0) \, \partial_m(f))$  
is given by $\nabla_i^{f} \nabla_j^{f} \,
\big{(}(\Gamma_{F_0})^m_{kl} - (\Gamma_{f}\big{)}^m_{kl}) \, \partial_m(f)$ and since the latter term is a section of the tangent bundle of the immersion $f$, we are thus naturally led to firstly consider flow lines of the evolution equation  
\begin{eqnarray}  \label{de_Turck_equation}
	\partial_t(f_t) = - \frac{1}{2} \mid A^0_{f_t} \mid^{-4} \,
	\big{(} 2\, \delta \Will(f_t) + \Tan_{F_0}(f_t) \big{)}  
	=: \Mill_{F_0}(f_t),    \qquad
\end{eqnarray} 
for some arbitrarily fixed $C^{\infty}-$smooth  
immersion $F_0:\Sigma \longrightarrow \rel^n$,
where the symbol $[f \mapsto \Tan_{F_0}(f)]$ denotes a 
globally well-defined differential operator of fourth order, which via locally defined coordinates is concretely given 
by the formula  
\begin{eqnarray}    \label{tangential.correction} 
	\Tan_{F_0}(f):=  g^{ij}_{f} \,g^{kl}_{f}\, g^{mr}_{f} \, 
	\langle \nabla_i^{f} \nabla_j^{f} \nabla_k^{f} \nabla_l^{f}(f), 
	\partial_m f \rangle\, \partial_r f                 \\
	- g^{ij}_{f} g^{kl}_{f}\, \nabla_i^{f} \nabla_j^{f} \,
	\big{(}(\Gamma_{F_0})^m_{kl} - (\Gamma_{f})^m_{kl}\big{)} \, \partial_m(f),                   \nonumber
\end{eqnarray}
for any $W^{4,2}$-immersion $f:\Sigma \longrightarrow \rel^n$. On account of formula (\ref{difference.Christoffel}) the right-hand side ``$\Mill_{F_0}(f_t)$'' of equation (\ref{de_Turck_equation}) can thus be expressed by: 
\begin{eqnarray}  \label{D_F_0}
	\Mill_{F_0}(f_t)(x) = 
	- \frac{1}{2} \mid A^0_{f_t} \mid^{-4} \,
	g^{ij}_{f_t} \, g^{kl}_{f_t}\, \nabla_i^{f_t} \nabla_j^{f_t} 
	\nabla_k^{F_0} \nabla_l^{F_0}(f_t)(x)   
	+ B(x,D_xf_t,D^2_xf_t,D_x^3 f_t),       
\end{eqnarray} 
for $(x,t) \in \Sigma' \times [0,T]$, and this differential operator is locally elliptic in the sense of 
Section 1 in \cite{Mantegazza}, as well. Here, the symbols 
$D_xf_t, D_x^2f_t, D_x^3 f_t$ abbreviate the matrix-valued functions $(\partial_{1}f_t,\partial_{2}f_t)$, 
$(\nabla^{F_0}_{ij}f_t)_{i,j \in \{1,2\}}$ and 
$(\nabla^{F_0}_{ijk}f_t)_{i,j,k \in \{1,2\}}$, and the 
$n$ components of the lower order term 
$B(\,\cdot \,,D_xf_t,D^2_x f_t, D_x^3 f_t)$ are ``scalars'', i.e. well-defined 
functions on $\Sigma' \times [0,T]$ whose values do not depend on the choice of the chart $\psi$ of $\Sigma'$. Hence, there has to exist some well-defined function $B: \Sigma \times \rel^{2n} \times \rel^{4n} \times \rel^{8n} \to \rel^n$ whose $n$ components are rational functions in their $14n$ real variables, such that equation (\ref{D_F_0}) holds ``globally'' for any pair $(x,t) \in \Sigma \times [0,T]$, and we arrive at the non-linear equation 
\begin{eqnarray}  \label{de_Turck_equation_2}
	\partial_t(f_t) = - \frac{1}{2} \mid A^0_{f_t} \mid^{-4} \,
	g^{ij}_{f_t} \, g^{kl}_{f_t}
	\, \nabla_i^{f_t} \nabla_j^{f_t} \nabla_k^{F_0} \nabla_l^{F_0}(f_t)
	+ B(\,\cdot\,,D_xf_t,D_x^2f_t,D_x^3f_t)  \qquad
\end{eqnarray} 
of fourth order with quasilinear leading term, 
called the ``DeTurck modification of the MIWF'',
whose linearization in any fixed family of 
$W^{4,p}$-immersions $f_t:\Sigma \longrightarrow \rel^n$ 
belonging to the parabolic $L^p$-space 
\begin{equation}   \label{X.T}
	X_{T,p}:= W^{1,p}([0,T];L^p(\Sigma,\rel^n)) \cap L^p([0,T];W^{4,p}(\Sigma,\rel^n)), \quad 
\end{equation} 
being well-defined for any $p \in (1,\infty)$, will turn out to yield a uniformly parabolic operator from $X_{T,p}$ to $L^p([0,T];L^{p}(\Sigma,\rel^n))$ - 
see Theorem \ref{Psi.of.class.C_1} below -  
provided we choose here $p>3$ and provided the immersions $f_t$ satisfy $\mid A^0_{f_t} \mid^2 >0$ on $\Sigma$, for every 
$t \in [0,T]$. The mathematical significance of the choice 
``$p>3$'' lies in embeddings (\ref{embedding.Sobol.Hoelder}), (\ref{C0.Hoelder.embedding.XT}) and (\ref{Hoelder.embedding.XT}) below, which might fail to be continuous for $p\in (1,3]$; see here also the footnote addressing embedding (\ref{Hoelder.embedding.XT}). 
It should be emphasized here, that embeddings (\ref{embedding.Sobol.Hoelder}), (\ref{trace.time.embedding}), (\ref{C0.Hoelder.embedding.XT}) and (\ref{Hoelder.embedding.XT}) play key roles in the proofs of Theorems 
\ref{Psi.of.class.C_1}--\ref{Functional.prop.Frechet.derivative} below.    \\
Now we fix some additional smooth immersion 
$U_0 \in \textnormal{Imm}_{\textnormal{uf}}(\Sigma,\rel^n) 
\cap C^{\infty}(\Sigma,\rel^n)$, 
consider the maximal flow line 
$\{\PP(t,0,U_0)\}_{t \in [0,T_{\textnormal{max}}(U_0))}$ 
of the MIWF-equation (\ref{Moebius.flow}) and choose some $0<T<T_{\textnormal{max}}(U_0)$,
where $T_{\textnormal{max}}(U_0)$ has been defined in 
Definition \ref{Umbilic.free}.
We recall from the proof of Theorem 1 in \cite{Jakob_Moebius_2016}, 
that there is a unique smooth family of smooth diffeomorphisms 
$\phi_t \equiv \phi_t^{U_0,F_0}:\Sigma \longrightarrow \Sigma$ with $\phi_0=\textnormal{id}_{\Sigma}$, $t\in [0,T]$, such that the reparametrization $\{\PP(t,0,U_0)\circ \phi_t^{U_0,F_0}\}_{t \in [0,T]}$ satisfies equation (\ref{de_Turck_equation}) respectively (\ref{de_Turck_equation_2}).
Now we fix some $p\in (3,\infty)$ and some open neighborhood 
$W_{U_0,T,p}$ of the above flow line 
$\{\PP(t,0,U_0)\circ \phi_t^{U_0,F_0}\}_{t\in [0,T]}$ 
in the space $X_{T,p}$ from line (\ref{X.T}), and we shall follow Spener's strategy in \cite{Spener}, to use the 
fact that there holds for the Banach space of traces of 
elements in $X_{T,p}$ at the initial time $t=0$: 
\begin{eqnarray}  \label{Trace.X.T}
	\gamma_0(X_{T,p}) \equiv \textnormal{trace}(X_{T,p})
	= (L^p(\Sigma,\rel^n),W^{4,p}(\Sigma,
	\rel^n))_{1-\frac{1}{p},p}              \\
	= B^{4-\frac{4}{p}}_{p,p}(\Sigma,\rel^n)     
	=W^{4-\frac{4}{p},p}(\Sigma,\rel^n),    \nonumber
\end{eqnarray} 
by Corollary 1.14 in \cite{Lunardi.interpolation} respectively Theorem 4.10.2 in Chapter III of \cite{Amann.1995}, combined with Theorem 4.3.3 in \cite{Triebel} and Theorem 1.2.4 in \cite{Triebel.2} respectively formula (3.5) in \cite{Amann.2}, 
which motivates us to consider the non-linear operator
$$
\Psi^{F_0,T}: W_{U_0,T,p} \subset X_{T,p}
\longrightarrow W^{4-\frac{4}{p},p}(\Sigma,\rel^n) \times 
L^p([0,T];L^p(\Sigma,\rel^n)) =: Y_{T,p} 
$$
defined by: 
\begin{eqnarray}  \label{Psi}
	\Psi^{F_0,T}(\{f_t\}_{t\in [0,T]})
	:=(\gamma_0(\{f_t\}), \{\partial_t(f_t) - 
	\Mill_{F_0}(f_t)\}_{t\in [0,T]}), 
\end{eqnarray}
compare here to Section 1.5 in Chapter III of \cite{Amann.1995}. 
We shall at first investigate the properties of exactly 
this operator in the following theorem, Theorem \ref{Psi.of.class.C_1}, and then in Theorem \ref{short.time.solution} we shall precisely address 
unique short-time existence of solutions 
to evolution equation (\ref{de_Turck_equation_2})
in the space $X_{T,p}$. 
The most important techniques in these two theorems are the 
explicit use of the quasi-linearity of the operator 
$[f\mapsto \Mill_{F_0}(f)]$ in lines 
(\ref{de_Turck_equation}) and (\ref{D_F_0}) 
in combination with Theorem 2.1 in \cite{Koehne.Pruess.Wilke.2010}, and the ``linearization'' of the quasi-linear operator (\ref{Psi}) in appropriate points of its open domain $W_{U_0,T,p} \subset X_{T,p}$, yielding linear parabolic differential operators and 
thus the applicability of classical ``Linear Functional Analysis'' - in accordance with the classical approach to the Ricci-flow; see e.g. Chapter 5 in \cite{Andrews.Hopper.2011}.      
\begin{theorem}           \label{Psi.of.class.C_1} 
	Let $\Sigma$ be a smooth compact torus, 
	$F_0 \in C^{\infty}_{\textnormal{Imm}}(\Sigma,\rel^n)$ 
	a smooth immersion, $U_0 \in \textnormal{Imm}_{\textnormal{uf}}(\Sigma,\rel^n) 
	\cap C^{\infty}(\Sigma,\rel^n)$ a smooth and umbilic-free initial immersion of $\Sigma$ into $\rel^n$, 
	and let $T\in (0,T_{\textnormal{max}}(U_0))$ and 
	$p\in (3,\infty)$ be arbitrarily chosen, where $T_{\textnormal{max}}(U_0)$ has been defined in 
	Definition \ref{Umbilic.free}. 
	Then there is a sufficiently small open neighbourhood $W_{U_0,T,p}$ in the Banach space $X_{T,p}$ about the smooth solution $\{\PP(t,0,U_0)\circ \phi_t^{U_0,F_0}\}_{t\in [0,T]}$ of the modified MIWF-equation (\ref{de_Turck_equation_2}), such that the following statements hold: 
	\begin{itemize}
	\item[1)] For any family of immersions 
	$\{f_t\} \in W_{U_0,T,p}$ there holds the identity
	\begin{equation}  \label{K.F_0.f_t.f_t}
	\Mill_{F_0}(\{f_t\}) = K^{F_0}(\{f_t\}).(\{f_t\}) \quad  
	\textnormal{in} \,\,\,L^p([0,T];L^p(\Sigma,\rel^n)),  
	\end{equation} 
	where $K^{F_0}$ is the non-linear operator from 
	line (\ref{A.F.0}), here applied to 
    $f_t \in \textnormal{Imm}_{\textnormal{uf}}(\Sigma,\rel^n)$ 
	for each $t \in [0,T]$.  
	    \item[2)] The map $\Psi^{F_0,T}:W_{U_0,T,p} \longrightarrow Y_{T,p}$, defined in line (\ref{Psi}), is of class $C^{\omega}$ on the open subset $W_{U_0,T,p}$ of the Banach space $X_{T,p}$ from line (\ref{X.T}).
		\item[3)] In any fixed element $\{f_t\} \in W_{U_0,T,p}$ the Fr\'echet derivative of the second component of $\Psi^{F_0,T}$ is a linear, uniformly parabolic operator of order $4$, whose leading operator acts on each component of $f=\{f_t\}$ separately:
		\begin{eqnarray} \label{Frechet.of.Psi}
		(D(\Psi^{F_0,T})_2(f)).(\eta) \equiv 
		\partial_t(\eta) - D(\Mill_{F_0})(f).(\eta) \quad \\
		= \partial_t(\eta) + \frac{1}{2} \mid A^0_{f_t} \mid^{-4} \,g^{ij}_{f_t} \, g^{kl}_{f_t} \, \nabla^{F_0}_{ijkl}(\eta)      
		+ B_3^{ijk} \cdot \nabla^{F_0}_{ijk}(\eta) + B_2^{ij} \cdot \nabla^{F_0}_{ij}(\eta) 
		+ B_1^{i} \cdot \nabla^{F_0}_{i}(\eta)          \nonumber
		\end{eqnarray} 
		on $\Sigma \times [0,T]$, for any element $\eta=\{\eta_t\}$ of the tangent space
		$T_{f}W_{U_0,T,p}=X_{T,p}$. Here, the coefficients 
		$\mid A^0_{f_t} \mid^{-4} \,g^{ij}_{f_t} \, g^{kl}_{f_t}$ of the leading order term are of class 
		$C^{0,\alpha}([0,T];C^{0,\alpha}(\Sigma,\rel))$, $B^{ij}_2$ and $B^i_1$ are the coefficients of $\textnormal{Mat}_{n,n}(\rel)$-valued, contravariant tensor fields of degrees $2$ and $1$, which depend on $x$, $D_xf_t,D^2_xf_t,D^3_xf_t$ and on $D^4_xf_t$ and are of class $L^{p}([0,T];L^p(\Sigma,
		\textnormal{Mat}_{n,n}(\rel)))$. 
		They depend on the spatial derivatives of third order 
		$\{\nabla^{F_0}_{ijk}f_t\}$ and of fourth order $\{\nabla^{F_0}_{ijkl}f_t\}$ of $\{f_t\}$ affine linearly. Finally $B_{3}^{ijk}$ are the coefficients 
		of a $\textnormal{Mat}_{n,n}(\rel)$-valued, 
		contravariant tensor field of degree $3$, which depends on $x$, $D_xf_t$ and $D^2_xf_t$ only and is of class $C^{0,\alpha}([0,T];C^{0,\alpha}(\Sigma,
		\textnormal{Mat}_{n,n}(\rel)))$, for $\alpha(p)\in (0,1)$ as in (\ref{Hoelder.embedding.XT}) below.
		\item[4)] The Fr\'echet derivative of $\Psi^{F_0,T}$ yields a topological isomorphism 
		$$
		D\Psi^{F_0,T}(f): T_{f}W_{U_0,T,p}=X_{T,p} \stackrel{\cong} \longrightarrow Y_{T,p}
		$$
		in any fixed family of immersions 
		$f \equiv \{f_t\} \in C^{4}(\Sigma \times [0,T],\rel^n) \cap W_{U_0,T,p}$.
	\end{itemize}	
\end{theorem} 
\qed
\noindent 
\begin{theorem}  \label{short.time.solution}
Let $\Sigma$ be a smooth compact torus, $p\in (3,\infty)$, and let $F_0 \in C^{\infty}_{\textnormal{Imm}}(\Sigma,\rel^n)$ 
and $U_0 \in \textnormal{Imm}_{\textnormal{uf}}(\Sigma,\rel^n)
\cap W^{4-\frac{4}{p},p}(\Sigma,\rel^n)$ 
be fixed immersions.  
\begin{itemize} 
\item[1)]  
There are sufficiently small numbers $T=T(U_0)>0$ and  
$\varepsilon=\varepsilon(U_0)>0$, such that  
for every initial immersion $u_0 \in W^{4-\frac{4}{p},p}(\Sigma,\rel^n)$ satisfying
\begin{equation}  \label{small}
\parallel U_0 - u_0 \parallel_{W^{4-\frac{4}{p},p}(\Sigma,\rel^n)}
<\varepsilon 
\end{equation}  
there is a unique solution $\{u(t,u_0)\}_{t\in [0,T]}$ 
of the non-linear Cauchy problem 
\begin{equation}  \label{initial.value.problem}
\partial_t f_t = \Mill_{F_0}(f_t)    \quad \textnormal{on}  \quad \Sigma \times [0,T] \quad \textnormal{with} \quad 
f_0 = u_0  \quad \textnormal{on}  \quad \Sigma,	
\end{equation}  
from lines (\ref{de_Turck_equation}) and (\ref{de_Turck_equation_2}), i.e. of the problem $\Psi^{F_0,T}(\{f_t\})=(u_0,0)$, within the Banach space
$X_{T,p}$. 	
\item[2)] Let $T=T(U_0)>0$ and  
$\varepsilon=\varepsilon(U_0)>0$ be as in the first 
part of the theorem. Then for every initial immersion 
$u_0 \in W^{4-\frac{4}{p},p}(\Sigma,\rel^n)$ 
satisfying condition (\ref{small}) the unique short-time solution $\{u(t,u_0)\}_{t\in [0,T]}$ of Cauchy problem 
(\ref{initial.value.problem}) possesses a unique extension 
$\{\PP^*(t,0,u_0)\}_{t\in [0,t^{+}(u_0))}$
to a maximal, half-sided open interval of existence 
$[0,t^{+}(u_0))$, which means that $t^{+}(u_0)>T$, 
that for every $T^*\in [T,t^{+}(u_0))$ the restriction 
$\{\PP^*(t,0,u_0)\}_{t\in [0,T^*]}$ of $\{\PP^*(\,\cdot\,,0,u_0)\}$ 
to $\Sigma \times [0,T^*]$ is the unique solution 
of Cauchy problem (\ref{initial.value.problem}) in $X_{T^*,p}$, 
and that $t^{+}(u_0)$ is maximal in $\rel_+ \cup \infty$ 
with this property.
\qed
\end{itemize}     
\end{theorem} 
\noindent
Moreover, for only smooth and umbilic-free initial 
immersions $U_0$ we will prove, that the maximal 
smooth solution of Cauchy problem 
(\ref{initial.value.problem}) becomes even real analytic 
for positive times, i.e. immediately after its initial time. 
This is the decisive regularity result 
of this paper. In its proof we will combine slight modifications of several results of Theorem \ref{Psi.of.class.C_1} respectively techniques of its 
proof with the modern method in \cite{Escher.Pruess.Simonett}, 
similarly to the procedure in Sections 3--6 of 
\cite{Shao.2015} or Sections 3--4 in \cite{Shao.2013}.  
Earlier respectively simpler methods, as e.g. 
in \cite{Angenent.1990} or in \cite{Lamm.Koch.2012}, 
apply in different mathematical situations 
and cannot yield the statement of Theorem 
\ref{real.analytic.flow} (ii) below. 
\begin{theorem}   \label{real.analytic.flow} 
Suppose that $\Sigma$ is a smooth compact torus, 
that $F_0$ is a fixed $C^{\infty}$-smooth immersion of 
$\Sigma$ into $\rel^n$, that
$U_0 \in \textnormal{Imm}_{\textnormal{uf}}(\Sigma,\rel^n)$
$\cap C^{\infty}(\Sigma,\rel^n)$
is an arbitrary smooth and umbilic-free initial 
immersion, and that $p\in (3,\infty)$.
\begin{itemize} 
	\item[1)] The maximal solution  
	$\{\PP^*(t,0,U_0)\}_{t\in [0,t^{+}(U_0))}$
	of the quasilinear parabolic Cauchy problem (\ref{initial.value.problem}) - here 
	with $u_0=U_0$ - from the second part of Theorem 
	\ref{short.time.solution} is of class  
    $C^{\infty}(\Sigma \times [0,t^{+}(U_0)),\rel^n)$,
    and it is therefore a maximal flow line of 
    evolution equation (\ref{de_Turck_equation})   
	in the sense of Definition \ref{Umbilic.free}, (d), 
	with $T_{\textnormal{max}}(U_0)=t^{+}(U_0)$. 
	\item[2)] Suppose that the smooth torus $\Sigma$ is  additionally endowed with a complex structure, and 
	that only the immersion $F_0$ is additionally of class  
    $C^{\omega}(\Sigma,\rel^n)$. 
    Then the maximal flow line 
    $\{\PP^*(t,0,U_0)\}_{t\in [0,T_{\textnormal{max}}(U_0))}$
    of evolution equation (\ref{de_Turck_equation}) 
    from the first part of this theorem is additionally 
    of class $C^{\omega}(\Sigma \times (0,T_{\textnormal{max}}(U_0)),\rel^n)$.
\qed
\end{itemize} 
\end{theorem}  
\noindent  		
We firstly prepare the proofs of 
Theorems \ref{Psi.of.class.C_1}--\ref{real.analytic.flow} 
by means of the following elementary lemma.
\begin{lemma}  \label{preparation}
\begin{itemize} 	 
\item[1)] For every $p \in (1,\infty)$, the differential 
operator $[f \mapsto  \Mill_{F_0}(f)]$ from line (\ref{de_Turck_equation}) is a well-defined map from 
$\textnormal{Imm}_{\textnormal{uf}}
(\Sigma,\rel^n) \cap W^{4,p}(\Sigma,\rel^n)$
to $L^{p}(\Sigma,\rel^n)$, 
and it has a quasilinear structure, which means here 
precisely that there holds:
\begin{equation}  \label{A.F_0.f.f}
	\Mill_{F_0}(f) = K^{F_0}(f).(f),  \,\,\,
	\textnormal{for every} \,\, f\in \textnormal{Imm}_{\textnormal{uf}}(\Sigma,\rel^n) 
	\cap W^{4,p}(\Sigma,\rel^n),  
\end{equation} 
where $K^{F_0}$ is a non-linear operator 
\begin{equation} \label{A.F.0} 
K^{F_0}: \textnormal{Imm}_{\textnormal{uf}}(\Sigma,\rel^n) 
\longrightarrow \Lift(W^{4,p}(\Sigma,\rel^n),
L^{p}(\Sigma,\rel^n)),
\end{equation} 
whose restriction to $\textnormal{Imm}_{\textnormal{uf}}(\Sigma,\rel^n)  
\cap W^{4-\frac{4}{p},p}(\Sigma,\rel^n)$, endowed 
with the $W^{4-\frac{4}{p},p}$-norm, is locally 
Lipschitz continuous, if here additionally $p > 3$.  
\item[2)] For every fixed 
$h \in \textnormal{Imm}_{\textnormal{uf}}
(\Sigma,\rel^n)$ and every $p\in (1,\infty)$, 
the value $-K^{F_0}(h):W^{4,p}(\Sigma,\rel^n) 
\longrightarrow L^{p}(\Sigma,\rel^n)$ of operator 
(\ref{A.F.0}) is of maximal $L^p$-regularity on 
$[0,T]$ w.r.t. the pair $(W^{4,p}(\Sigma,\rel^n)$,  
$L^{p}(\Sigma,\rel^n))$, i.e. $-K^{F_0}(h)$ is of 
class $\MR_p([0,T];W^{4,p}(\Sigma,\rel^n),
L^{p}(\Sigma,\rel^n))$ in the notation of 
\cite{Amann.2004} or \cite{Koehne.Pruess.Wilke.2010}.  
\end{itemize}  
\end{lemma} 
\proof 
\begin{itemize} 
    \item[1)] First of all we would like to refine 
	the rough, non-linear formulation (\ref{D_F_0}) 
	of the differential operator $\Mill_{F_0}$ 
	appearing in line (\ref{de_Turck_equation}). 
	To this end, we recall from lines 
	(\ref{triangle.perp})--(\ref{tangential.correction}) that   
	\begin{eqnarray}    \label{M.F_0} 
	\Mill_{F_0}(f) = 
	-\frac{1}{2}\,\frac{1}{|A^0_{f}|^4} 
	\Big{(} g^{ij}_{f} \, g^{kl}_{f} 
	\, \nabla_i^{f} \nabla_j^{f} \nabla_k^{f} \nabla_l^{f}(f)                 
	- g^{ij}_{f} g^{kl}_{f}\, \nabla_i^{f} \nabla_j^{f} \,
	\big{(}(\Gamma_{F_0})^m_{kl} - (\Gamma_{f})^m_{kl}\big{)} \, \partial_m(f)                \nonumber           \\
	+ 2\,Q(A_f)(\vec H_f) 
	- \frac{1}{2} \mid \vec H_f \mid^2\,\vec H_f \Big{)} \nonumber\\  
	= -\frac{1}{2}\,\frac{1}{|A^0_{f}|^4} \,
	\Big{(} g^{ij}_{f} \, g^{kl}_{f} 
	\, \nabla_i^{f} \nabla_j^{f} \nabla_k^{F_0} \nabla_l^{F_0}(f)     	          
	+ g^{ij}_{f} g^{kl}_{f}\,\,
	\big{(}(\Gamma_{F_0})^m_{kl} - (\Gamma_{f})^m_{kl}\big{)} \,  \nabla_i^{f} \nabla_j^{f}(\partial_m(f))  \nonumber  \\
	+\, g^{ij}_{f} g^{kl}_{f}\,\,
	\nabla_i^{f} \big{(}(\Gamma_{F_0})^m_{kl} - (\Gamma_{f})^m_{kl}\big{)} \, 
	\nabla_j^{f}(\partial_m(f))	           
	+ 2\,Q(A_f)(\vec H_f) 
	- \frac{1}{2} \mid \vec H_f \mid^2\,\vec H_f \Big{)}, \quad 
	\end{eqnarray} 
    for any fixed umbilic-free immersion $f \in W^{4,p}(\Sigma,\rel^n)$ with $p \in (1,\infty)$.  
    Using now formulae (\ref{double_nabla})--(\ref{mean_curvat_laplacian})
    and the first two lines of the general derivation formulae (\ref{derivation_formulae}) below for covariant 
    tensor fields, one can deduce - similarly to the 
    technical reasoning in the proof of the second part 
    of Theorem \ref{Psi.of.class.C_1} - the following more systematic formulation of the right hand side in (\ref{M.F_0}):   
    \begin{eqnarray}    \label{M.F_0.2} 
    \Mill_{F_0}(f)(x) =
    -\frac{1}{2}\,\frac{1}{|A^0_{f}|^4(x)} \,
    g^{ij}_{f} \, g^{kl}_{f} 
    \, \nabla_i^{f} \nabla_j^{f} \nabla_k^{F_0} \nabla_l^{F_0}(f)(x)     	          \nonumber     \\
    + \NN^{F_0}(x,D_xf(x),D_x^2f(x)) 
    \cdot D^3_xf(x) \\
    + \CC^{F_0}(x,D_xf(x),D_x^2f(x)) \cdot D_x^2f(x)      \nonumber     \\
    + \dom^{F_0}(x,D_xf(x),D_x^2f(x)) \cdot D_xf(x) 
    \nonumber
    \end{eqnarray}  
for $x \in \Sigma$ and for every 
fixed umbilic-free immersion $f \in W^{4,p}(\Sigma,\rel^n)$, where the functions 
$\NN^{F_0}:\Sigma \times \rel^{2n} \times \rel^{4n}  \longrightarrow \textnormal{Mat}_{n,8n}(\rel)$, 
$\CC^{F_0}: \Sigma \times \rel^{2n} \times \rel^{4n} \longrightarrow \textnormal{Mat}_{n,4n}(\rel)$ and 
$\dom^{F_0}: \Sigma \times \rel^{2n} \times \rel^{4n} \longrightarrow \textnormal{Mat}_{n,2n}(\rel)$
have the same algebraic structures as the function
$M^{F_0}$ in formula (\ref{F.third.derivatives}) below. 
In particular, the right-hand side of formula 
(\ref{M.F_0.2}) depends only affine 
linearly on the covariant derivatives of $f$
of fourth and third order. Hence, formula (\ref{M.F_0.2}) 
shows especially, that the non-linear differential operator  
$\Mill_{F_0}$ maps the intersection  
$\textnormal{Imm}_{\textnormal{uf}}
(\Sigma,\rel^n) \cap W^{4,p}(\Sigma,\rel^n)$
into $L^{p}(\Sigma,\rel^n)$.    
Now we are able to guess the desired non-linear map 
$K^{F_0}$ appearing in assertions (\ref{A.F_0.f.f}) and 
(\ref{A.F.0}) to be: 
\begin{eqnarray} \label{A.F_0} 
(K^{F_0}(h).\eta)(x):=
-\frac{1}{2}\,\frac{1}{|A^0_{h}|^4} \,
g^{ij}_{h} \, g^{kl}_{h} 
\, \nabla_i^{h} \nabla_j^{h} \nabla_k^{F_0} \nabla_l^{F_0}(\eta)(x)     	           \nonumber     \\
+ \NN^{F_0}(x,D_xh(x),D_x^2h(x)) 
\cdot D^3_x\eta(x)          \nonumber \\
+ \CC^{F_0}(x,D_xh(x),D_x^2h(x)) \cdot D_x^2\eta(x)     \nonumber     \\
+ \dom^{F_0}(x,D_xh(x),D_x^2h(x)) \cdot D_x\eta(x),   
\end{eqnarray}
for $x\in \Sigma$, and for any fixed umbilic-free immersion 
$h$ and $\eta \in W^{4,p}(\Sigma,\rel^{n})$, 
and formula (\ref{M.F_0.2}) shows us indeed, that the choice of 
the non-linear map $K^{F_0}$ in (\ref{A.F_0}) yields 
assertion (\ref{A.F_0.f.f}) for every umbilic-free immersion 
$f \in W^{4,p}(\Sigma,\rel^{n})$.  
Now, by formula (\ref{A.F_0})
the linear operator $K^{F_0}(h)$ has continuous coefficients 
and thus belongs to 
$\Lift(W^{4,p}(\Sigma,\rel^n),L^{p}(\Sigma,\rel^n))$, 
for every $h \in \textnormal{Imm}_{\textnormal{uf}}(\Sigma,\rel^n)$ 
and for every $p>1$. If we assume additionally that $p>3$, 
then we can use the fractional Sobolev embedding  
\begin{equation}  \label{embedding.Sobol.Hoelder}
	W^{4-\frac{4}{p},p}(\Sigma,\rel^n)
	\hookrightarrow C^{2,\alpha}(\Sigma,\rel^n), 
\end{equation}
for $0<\alpha <1$ with $\alpha \leq 2-\frac{6}{p}$,
and the classical mean value theorem, in order to 
infer again from formula (\ref{A.F_0}), that the 
map $K^{F_0}$ in (\ref{A.F.0}) is locally Lipschitz 
continuous on $\textnormal{Imm}_{\textnormal{uf}}(\Sigma,\rel^n)
\cap W^{4-\frac{4}{p},p}(\Sigma,\rel^n)$. 
\item[2)] From the first part of the lemma we know, 
that for any fixed 
$h\in \textnormal{Imm}_{\textnormal{uf}}(\Sigma,\rel^n)$
we obtain a bounded, linear differential operator $-K^{F_0}(h):W^{4,p}(\Sigma,\rel^n)
\longrightarrow L^{p}(\Sigma,\rel^n)$ with  
continuous coefficients. Moreover, again formula 
(\ref{A.F_0}) shows, that the operator 
$-K^{F_0}(h)$ is uniformly elliptic in the sense of 
Proposition \ref{A.priori.estimate} below. Hence, 
Proposition \ref{A.priori.estimate} guarantees 
us, that $\partial_t - K^{F_0}(h)$ maps the Banach space 
$X^0_T$ - defined below in (\ref{global.initial.value.problem}) - isomorphically onto $L^{p}([0,T],L^p(\Sigma,\rel^n))$, 
for any fixed $p \in (1,\infty)$, thus that $-K^{F_0}(h)$ 
is of class $\MR_p([0,T];W^{4,p}(\Sigma,\rel^n),
L^{p}(\Sigma,\rel^n))$ for every $p\in (1,\infty)$, 
using the notation of \cite{Amann.2004} or \cite{Koehne.Pruess.Wilke.2010}.  
\end{itemize}
\qed
	
\underline{Proof of Theorem \ref{Psi.of.class.C_1}}
\begin{itemize}
\item[1)] On account of the embedding 
\begin{equation} \label{trace.time.embedding} 
X_{T,p} \hookrightarrow C^0([0,T],W^{4-\frac{4}{p},p}(\Sigma,\rel^n))
\end{equation}
with $T$-independent embedding constant 
- see here Proposition 1.4.2 and Theorem 4.10.2 of Chapter III 
in \cite{Amann.1995}, Corollary 1.14 in \cite{Lunardi.interpolation} and formula (\ref{Trace.X.T}) above - and on account of embedding (\ref{embedding.Sobol.Hoelder}), 
we can derive identity (\ref{K.F_0.f_t.f_t}) 
immediately from identity (\ref{A.F_0.f.f}), for any 
family of immersions $\{f_t\} \in W_{U_0,T,p}$, provided 
the neighbourhood $W_{U_0,T,p}$ about the flow line 
$\{\PP(t,0,U_0)\circ \phi_t^{U_0,F_0}\}$ in $X_{T,p}$
has been chosen sufficiently small.        
\item[2)] The first component of the operator $\Psi^{F_0,T}$ from line (\ref{Psi}) is linear and continuous on account of Theorem 4.10.2 in Chapter III of \cite{Amann.1995}, thus it is of class $C^{\omega}$. In order to prove that also the second component of the operator $\Psi^{F_0,T}$ is of class $C^{\omega}$, we shall firstly prove its $C^1$-regularity, following the strategy of the proof of Theorem 2 in \cite{Jakob_Moebius_2016}, because this preparation 
considerably simplifies the proof of its real analyticity 
and also the explicit computation of the Fr\'echet derivative of the operator $\Psi^{F_0,T}$ in any fixed element $\{f_t\}$ of a sufficiently small neighbourhood $W_{U_0,T,p}$ about the flow line 
$\{\PP(t,0,U_0)\circ \phi_t^{U_0,F_0}\}$ in $X_{T,p}$. \\
$\alpha$) At first we aim at a slight improvement of 
formula (\ref{M.F_0.2}) above.  
To this end, we fix an arbitrary smooth chart 
$\psi:O \stackrel{\cong} \longrightarrow \Sigma'$ of an arbitrary coordinate patch $\Sigma'$ of $\Sigma$, which yields partial derivatives $\partial_m$, $m=1,2$, on $\Sigma'$. 
For any $C^2$-immersion $G:\Sigma \longrightarrow \rel^n$ the choice of $\psi$ yields the coefficients 
$g_{ij}:=\langle \partial_{i}G, \partial_{j}G \rangle_{\rel^n}$ of the first fundamental form of $G$ 
w.r.t. $\psi$ and the associated Christoffel-symbols 
$(\Gamma_G)^m_{kl} := g^{mj} \,
\langle \partial_{kl} G,\partial_jG \rangle_{\rel^n}$ of $(\Sigma',G^*(g_{\textnormal{eu}}))$. On account of (\ref{double_nabla}) and by the general derivation formulae 
\begin{eqnarray}  \label{derivation_formulae}
		\nabla_i^{G} (\omega_k) = \partial_i(\omega_k) - (\Gamma_G)^m_{ik}\, \omega_m \nonumber \\
		\nabla_i^G(\lambda_{jk})=\partial_i(\lambda_{jk}) - (\Gamma_G)^m_{ij} \, \lambda_{mk} 
		-(\Gamma_G)^m_{ik} \, \lambda_{jm} \\
		\nabla_i^G(\zeta_{jkl}) = \partial_i(\zeta_{jkl}) - (\Gamma_G)^m_{ij} \, \zeta_{mkl} 
		-(\Gamma_G)^m_{ik} \, \zeta_{jml} - (\Gamma_G)^m_{il} \, \zeta_{jkm} \nonumber
	\end{eqnarray}
on $\Sigma'$ for the coefficients $\omega_k$, $\lambda_{jk}$ and $\zeta_{jkl}$ (w.r.t. to the chart $\psi$) of covariant $C^1$-tensor fields $\omega$, $\lambda$ and $\zeta$ of degrees $1$, $2$ and $3$, one can verify that for any $W^{4,p}$-immersion $f:\Sigma \longrightarrow \rel^n$, for any fixed $C^{\infty}$-immersion $G:\Sigma \longrightarrow \rel^n$ and for fixed $i,j,k,l \in \{1,2\}$ there is a unique rational function $P^G_{(ijkl)} \in C^{\infty}(\Sigma')[v_1,\ldots,v_{2n},w_1,
\ldots,w_{4n},y_1,\ldots,y_{8n}]$ (with $n$ real components) 
in $14n$ real variables, whose coefficients are rational expressions involving the partial derivatives  
$\partial_i(G^1),\ldots, \partial_i(G^n), \partial_{ij}(G^1),\ldots, \partial_{ij}(G^n),  
\partial_{ijk}(G^1)$, $\ldots$, $\partial_{ijk}(G^n)$ of the components of $G$ up to third order, such that 
	\begin{eqnarray}  \label{4nabla_hurra}
	\nabla_i^{f} \nabla_j^{f} \nabla_k^{G} \nabla_l^{G}(f) 
	= \nabla_i^{G} \nabla_j^{G} \nabla_k^{G} \nabla_l^{G}(f)                  \qquad \\
	+P^G_{(ijkl)}(\nabla^{G}_{1}(f^1),\ldots,\nabla^{G}_{2}(f^n),\nabla^{G}_{11}(f^1),\ldots, \nabla^{G}_{22}(f^n),
	\nabla^{G}_{111}(f^1),\ldots, \nabla^{G}_{222}(f^n)) \nonumber 
	\end{eqnarray}
	holds on $\Sigma'$. We note here that the terms
	$$
	P^G_{(ijkl)}(\nabla^{G}_{1}(f^1),\ldots,\nabla^{G}_{2}(f^n),\nabla^{G}_{11}(f^1), \ldots, \nabla^{G}_{22}(f^n),
	\nabla^{G}_{111}(f^1), \ldots, \nabla^{G}_{222}(f^n))
	$$ 
	must be the coefficients of a covariant tensor field of fourth degree on $\Sigma'$, and moreover we note, that this expression is linear w.r.t. the derivatives 
	$\nabla^{G}_{rsm}(f^1), \ldots ,\nabla^{G}_{rsm}(f^n)$ of third order of the components of $f$, i.e. we have  
	\begin{eqnarray}    \label{P.third.derivatives}
		P^G_{(ijkl)}(\nabla^{G}_{1}(f^1),\ldots,\nabla^{G}_{2}(f^n),\nabla^{G}_{11}(f^1), \ldots, \nabla^{G}_{22}(f^n),
		\nabla^{G}_{111}(f^1), \ldots, \nabla^{G}_{222}(f^n))   
		             \nonumber       \\
		=M^{G}_{(ijkl)}(\nabla^{G}_{1}(f^1),\ldots,
		\nabla^{G}_{2}(f^n), 
		\nabla^{G}_{11}(f^1), \ldots, \nabla^{G}_{22}(f^n)) \cdot                 \nonumber      \\            	\cdot                                            
		\Big{(} \nabla^{G}_{111}(f^1),\nabla^{G}_{112}(f^1),\ldots,
		\nabla^{G}_{221}(f^n),\nabla^{G}_{222}(f^n)  \Big{)}^T                           \nonumber \\
		+ \tilde P^G_{(ijkl)}(\nabla^{G}_{1}(f^1),\ldots,
		\nabla^{G}_{2}(f^n),\nabla^{G}_{11}(f^1), \ldots, \nabla^{G}_{22}(f^n))     \qquad               
	\end{eqnarray}
	where the expression 
	$$
	M^{G}_{(ijkl)}(\nabla^{G}_{1}(f^1),\ldots,
	\nabla^{G}_{2}(f^n),\nabla^{G}_{11}(f^1), \ldots, \nabla^{G}_{22}(f^n))
	$$ 
	is a $\textnormal{Mat}_{n,8n}(\rel)$-valued function, whose entries are rational expressions involving the partial derivatives 
	$\partial_i(G^1),\ldots, \partial_i(G^n), \partial_{ij}(G^1),\ldots, \partial_{ij}(G^n)$,
	$\partial_{ijk}(G^1)$,$\ldots$, $\partial_{ijk}(G^n)$ of 
	$G$ up to third order and the derivatives $\nabla^{G}_{1}(f^1)$,$\ldots$, $\nabla^{G}_{2}(f^n)$, 
	$\nabla^{G}_{11}(f^1),\ldots, \nabla^{G}_{22}(f^n)$ of $f$ up to only second order. Compare here also with the statement of Lemma 3.4 in \cite{Spener}. Moreover, since there holds 
	\begin{equation}  \label{A0.squared}
	\mid A^0_f \mid^2 = g_f^{ik} \,g_f^{jl}\, \langle (A^0_f)_{ij}, (A^0_f)_{kl} \rangle_{\rel^n}
	\end{equation}
	and $(A_f)_{ij} \equiv A_f(\partial_i,\partial_j) 
	= \partial_{ij} f - (\Gamma_f)^m_{ij} \,\partial_mf$ for every umbilic-free $W^{4,p}$-immersion $f$ by formulae 
	(\ref{double_nabla}) and (\ref{second.fund.i.j}), 
	one can easily verify by means of formula (\ref{Cramers.rule}) below, that there is a unique rational function $Q^G\in C^{\infty}(\Sigma')[v_1,\ldots,v_{2n},w_1,\ldots,w_{4n}]$
	in $6n$ real variables, whose coefficients are rational expressions involving the partial derivatives 
	$\partial_i(G^1)$, $\ldots$, $\partial_i(G^n)$, $\partial_{ij}(G^1),\ldots, \partial_{ij}(G^n)$ 
	of the components of $G$ up to second order, such that
	\begin{eqnarray}  \label{A0_hurra}
		\mid A^0_f \mid^4 = 
		Q^G(\nabla^{G}_{1}(f^1),\nabla^{G}_{1}(f^2),\ldots,
		\nabla^{G}_{2}(f^n), \nabla^{G}_{11}(f^1),\nabla^{G}_{11}(f^2),\ldots, \nabla^{G}_{22}(f^n)) \qquad
	\end{eqnarray}
	holds on $\Sigma'$, which is a ``scalar'', i.e. which does not depend on the choice of the chart $\psi$ of $\Sigma'$. Combining (\ref{4nabla_hurra}) and (\ref{A0_hurra}) and recalling that every immersion $\PP(t,0,U_0)\circ \phi_t^{U_0,F_0}$ of the considered flow line of equation (\ref{de_Turck_equation_2}) does not have any umbilic 
	points on $\Sigma$, we obtain the existence of a unique  
	function $R^{F_0} \in C^{\infty}(\Sigma')[v_1,\ldots,v_{2n}, w_1,\ldots,w_{4n},y_1,\ldots, y_{8n}]$, 
	whose components are rational functions in their $14n$ real variables, whose coefficients are rational functions of the partial derivatives of the components of $F_0$ up to third order and which is linear w.r.t. the last $8n$ variables $y_1,\ldots, y_{8n}$ on account of formula (\ref{P.third.derivatives}), such that 
	\begin{eqnarray} \label{nabla_umgestulpi}
		\frac{1}{2} \mid A^0_{f_t} \mid^{-4} \,
		g^{ij}_{f_t} \, g^{kl}_{f_t}
		\, \nabla_i^{f_t} \nabla_j^{f_t} \nabla_k^{F_0} \nabla_l^{F_0}(f_t)    \qquad        \\
		= \frac{1}{2} \mid A^0_{f_t} \mid^{-4} \,
		g^{ij}_{f_t} \, g^{kl}_{f_t} \, \nabla_i^{F_0} \nabla_j^{F_0} \nabla_k^{F_0} \nabla_l^{F_0}(f_t) \qquad \nonumber \\
		+R^{F_0}(\nabla^{F_0}_{1}(f_t^1),\ldots,
		\nabla^{F_0}_{2}(f_t^n),\nabla^{F_0}_{11}(f_t^1),
		\ldots,\nabla^{F_0}_{22}(f_t^n),
		\nabla^{F_0}_{111}(f_t^1),\ldots, \nabla^{F_0}_{222}(f_t^n))      \qquad \nonumber
	\end{eqnarray}
	holds on $\Sigma' \times [0,T]$ for any family of immersions $\{f_t\}$ in a sufficiently small neighbourhood $W_{U_0,T,p}$ of the fixed flow line $\{\PP(t,0,U_0)\circ 
	\phi_t^{U_0,F_0}\}_{t\in [0,T]}$ in 
	the space $X_{T,p}$ from line (\ref{X.T}). 
	At this point we have already silently used the
	composition of embedding
	(\ref{trace.time.embedding}) with embedding (\ref{embedding.Sobol.Hoelder}), yielding the 
	continuous embedding  
	\begin{equation}  \label{C0.Hoelder.embedding.XT} 
	X_{T,p} \hookrightarrow C^{0}([0,T];C^{2,\alpha}(\Sigma,\rel^n))
	\end{equation} 
	for any $\alpha \in (0,1)$ with 
	$\alpha \leq 2-\frac{6}{p}$, provided $p \in (3,\infty)$.
	We should note here, that the $n$ components of 
	$$
	R^{F_0}(\nabla^{F_0}_{1}(f_t^1),\ldots,
	\nabla^{F_0}_{2}(f_t^n),
	\nabla^{F_0}_{11}(f_t^1),\ldots, \nabla^{F_0}_{22}(f_t^n),
	\nabla^{F_0}_{111}(f_t^1),\ldots,
	\nabla^{F_0}_{222}(f_t^n))
	$$ 
	are ``scalars'' as well, i.e. they do not depend on the choice of the chart $\psi$ of $\Sigma'$. On account of formulae (\ref{M.F_0.2}) and (\ref{Psi}), we infer from (\ref{nabla_umgestulpi}) the following representation of the second component of the operator $\Psi^{F_0,T}$:   
	\begin{eqnarray} \label{phi_revisited} 
	\partial_t f_t(x) - \Mill_{F_0}(f_t)(x)      
	\nonumber \\
	= \partial_t f_t(x) + \frac{1}{2} \mid A^0_{f_t}(x) \mid^{-4} \,g^{ij}_{f_t} \, g^{kl}_{f_t} \, \nabla_i^{F_0} \nabla_j^{F_0} \nabla_k^{F_0} \nabla_l^{F_0}(f_t)(x) \nonumber \\
	+R^{F_0}(\nabla^{F_0}_{1}(f_t^1),\ldots,
	\nabla^{F_0}_{2}(f_t^n),\ldots, \nabla^{F_0}_{111}(f_t^1),\nabla^{F_0}_{111}(f_t^2), \ldots,\nabla^{F_0}_{222}(f_t^n))(x) \nonumber \\
	- \NN^{F_0}(x,D_xf_t(x),D_x^2f_t(x)) 
	\cdot D^3_xf_t(x)                    \nonumber  \\
	- \CC^{F_0}(x,D_xf_t(x),D_x^2f_t(x)) \cdot D_x^2f_t(x)                           \nonumber  \\
	- \dom^{F_0}(x,D_xf_t(x),D_x^2f_t(x)) \cdot D_xf_t(x)                  \nonumber \\
	=: \partial_t f_t(x) + \frac{1}{2} \mid A^0_{f_t}(x) \mid^{-4} \,g^{ij}_{f_t} \, g^{kl}_{f_t} \, \nabla_i^{F_0} \nabla_j^{F_0} 
	\nabla_k^{F_0} \nabla_l^{F_0}(f_t)(x) \qquad \\
	+\F^{F_0}(x,D_xf_t(x),D_x^2f_t(x),D_x^3f_t(x)) \qquad \nonumber
	\end{eqnarray}
	for any family of immersions $\{f_t\} \in W_{U_0,T,p}$ and 
	for $(x,t) \in \Sigma' \times [0,T]$, if the neighbourhood 
	$W_{U_0,T,p}$ of the flow line $\{\PP(t,0,U_0)\circ \phi_t^{U_0,F_0}\}_{t\in [0,T]}$ in $X_{T,p}$ is chosen 
	sufficiently small. Here and in the sequel, the symbols 
	$D_xf, D_x^2f, D_x^3f, \ldots$ 
	abbreviate the matrix-valued functions $(\partial_{1}f,\partial_{2}f)$,    \\ 
	$(\nabla^{F_0}_{ij}f)_{i,j \in \{1,2\}}$,
	$(\nabla^{F_0}_{ijk}f)_{i,j,k \in \{1,2\}}$, $\ldots$. 
	On account of the arbitrariness of the choice of the coordinate patch $\Sigma'$ and its chart $\psi$ and by the compactness of $\Sigma$ equation (\ref{phi_revisited}) gives rise to a unique and well-defined function 
	$\F^{F_0}: \Sigma \times \rel^{2n} \times \rel^{4n} \times \rel^{8n} \longrightarrow \rel^n$ whose $n$ components are rational functions in their $14n$ real variables, such that (\ref{phi_revisited}) holds on $\Sigma \times [0,T]$, and not only on $\Sigma' \times [0,T]$. 
	It is important to note, that by (\ref{M.F_0.2}) and (\ref{nabla_umgestulpi}) the function $f \mapsto \F^{F_0}(\,\cdot\,,D_xf,D_x^2f,D_x^3f)$ in 
	(\ref{phi_revisited}) has the same algebraic structure 
	as the function $f \mapsto R^{F_0}(\nabla^{F_0}_{1}(f^1),\ldots,
	\nabla^{F_0}_{222}(f^n))$, in particular $f\mapsto \F^{F_0}(\,\cdot\,,D_xf,D_x^2f,D_x^3f)$ is
	affine linear w.r.t. the derivatives $(\nabla^{F_0}_{ijk}f)_{i,j,k \in \{1,2\}}$ of third order, similarly to the illustration in formula (\ref{P.third.derivatives}) on $\Sigma'$. Hence, 
	there are functions 
	$M^{F_0}:\Sigma \times \rel^{2n} \times \rel^{4n}  
	\longrightarrow \textnormal{Mat}_{n,8n}(\rel)$, 
	$$
	[(x,v_1,\ldots,v_{2n},w_1,\ldots,w_{4n}) \mapsto M^{F_0}(x,v_1,\ldots,v_{2n},w_1,\ldots,w_{4n})],
	$$
	and $\tilde \F^{F^0}: \Sigma \times \rel^{2n} \times \rel^{4n} \longrightarrow \rel^n$ depending smoothly on $x$ and rationally on the $6n$ variables $v_i$ and $w_j$, such that
	\begin{eqnarray}    \label{F.third.derivatives}
		\F^{F_0}(x,v_1,\ldots,v_{2n},
		w_1,\ldots,w_{4n},y_1,\ldots,y_{8n})                  \nonumber            \\
		= M^{F_0}(x,v_1,\ldots,v_{2n},w_1,\ldots,w_{4n}) \cdot  (y_1, y_2, \ldots, y_{8n-1},y_{8n})^{T}    \\
		+ \tilde \F^{F_0}(x,v_1,\ldots,v_{2n}, w_1,\ldots,w_{4n}).   \nonumber
	\end{eqnarray}  
	Now, since the entire flow line $\{\PP(t,0,U_0)\circ \phi_t^{U_0,F_0}\}$, starting in the fixed umbilic-free initial immersion $U_0$, solves equation (\ref{de_Turck_equation_2}) within the space $C^{\infty}(\Sigma \times [0,T],\rel^n)$ and thus consists of umbilic-free immersions, all summands in line (\ref{phi_revisited}) remain uniformly bounded on $\Sigma \times [0,T]$. Hence, there are open connected neighbourhoods $\OO_t$ of the graphs 
	\begin{eqnarray*}
	G_t:=\{(x,D_x(\PP(t,0,U_0)\circ \phi_t^{U_0,F_0})(x), 
	D_x^2(\PP(t,0,U_0)\circ \phi_t^{U_0,F_0})(x))|\, x\in \Sigma\} \\
	\subset \Sigma \times \rel^{6n}
    \end{eqnarray*}
	for $t\in [0,T]$, such that $M^{F_0}$ and $\tilde \F^{F_0}$
	are $C^{\infty}$-smooth on the open subset 
	$\bigcup_{t\in [0,T]} \OO_t$ of $\Sigma \times \rel^{6n}$. Now, on account of the compactness of $\Sigma$ and of $[0,T]$ there is some $\delta>0$, such that there holds:  
	\begin{equation}  \label{thickness} 
		\textnormal{dist}_{\Sigma \times \rel^{6n}}(G_t,\partial \OO_t) > 2 \delta   \qquad   \forall \, t\in [0,T].
	\end{equation}    
	Similarly to the proof of Theorem 2 in \cite{Jakob_Moebius_2016} and in view of inequality (\ref{thickness}) and of embedding (\ref{C0.Hoelder.embedding.XT}) we define an open neighbourhood $\V_{\delta}$ of the function
	$$
	[(x,t) \mapsto (D_x(\PP(t,0,U_0)\circ \phi_t^{U_0,F_0})(x), 
	D_x^2(\PP(t,0,U_0)\circ \phi_t^{U_0,F_0})(x))] 
	$$
	within the space $C^{0}([0,T],C^{0,\alpha}(\Sigma,\rel^{6n}))$, with the exponent $\alpha=\alpha(p) \in (0,1)$ from line (\ref{C0.Hoelder.embedding.XT}), by setting 
	\begin{eqnarray*}
	\bar h \in \V_{\delta} :\Longleftrightarrow 
	\parallel (D_x(\PP(t,0,U_0)\circ \phi_t^{U_0,F_0}), 
	D_x^2(\PP(t,0,U_0) \circ \phi_t^{U_0,F_0})) 
	- \bar h_t \parallel_{L^{\infty}(\Sigma,\rel^{6n})} \\
	< \delta  \,\,\,   \forall \, t\in [0,T].  \qquad
	\end{eqnarray*}
	It is an immediate consequence of this definition and of inequality (\ref{thickness}), that there holds: 
	\begin{equation}  \label{if.delta}
		\bar h \in \V_{\delta} \Longrightarrow 
		(x,\bar h_t(x)) \in \OO_t \subset \Sigma \times \rel^{6n}  
		\quad \forall \, (x,t) \in \Sigma \times [0,T].
	\end{equation}      
	Now, in order to show that the non-linear operator 
	\begin{eqnarray*} 
		W_{U_0,T,p} \ni [(x,t) \mapsto f_t(x)]  \mapsto 
		[(x,t) \mapsto \F^{F_0}(x,D_xf_t(x),D^2_xf_t(x),D^3_xf_t(x))] \in   \\
		\in L^p([0,T],L^p(\Sigma,\rel^n))
	\end{eqnarray*}
	is of class $C^1$ on a sufficiently small open neighborhood 
	$W_{U_0,T,p}$ of the flow line $\{\PP(t,0,U_0)\circ \phi_t^{U_0,F_0}\}$ in 
	$X_{T,p}$, we shall proceed similarly to the proof of Theorem 2 in \cite{Jakob_Moebius_2016} and 
	consider here the Banach space 
	\begin{equation}  \label{B.space}
	B:=C^{0}([0,T];C^{0,\alpha}(\Sigma,\rel^{6n})) \times L^{p}([0,T];W^{1,p}(\Sigma,\rel^{8n})),
	\end{equation}  
	equipped with the canonical product norm
	$$
	\parallel h \parallel_{B}:= \parallel \bar h \parallel_{C^{0}([0,T];
	C^{0,\alpha}(\Sigma,\rel^{6n}))}  + 
    \parallel h^* \parallel_{L^{p}([0,T];W^{1,p}(\Sigma,\rel^{8n}))}
	$$ 
	for $h=(\bar h,h^*)\in B$. We firstly derive from equation (\ref{F.third.derivatives}),
	that the classical derivative of the function 
	$\F^{F_0}$ in any point $(v,w,y)\in \rel^{2n} \times \rel^{4n} \times \rel^{8n}$, in which $\F^{F_0}$ is 
	smooth, and in direction of any fixed vector 
	$(\vec v,\vec w,\vec y) \in \rel^{2n} \times \rel^{4n} \times \rel^{8n}$ reads:   
	\begin{eqnarray}  \label{Derivative} 
		D_{(v,w,y)}\F^{F_0}(x,(v,w,y))\cdot 
		(\vec v,\vec w,\vec y)                      \\ 
		= D_{(v,w)}(M^{F_0}(x,v,w) \cdot y) \cdot 
		(\vec v,\vec w)
		+ M^{F_0}(x,v,w) \cdot \vec y    \nonumber  \\
		+ D_{(v,w)}\tilde\F^{F_0}(x,v,w)\cdot 
		(\vec v,\vec w).                \nonumber
	\end{eqnarray}
	Now we fix some $h \in B$, whose first component 
	$\bar h \in C^{0}([0,T];C^{0,\alpha}(\Sigma,\rel^{6n}))$ 
	is contained in $\V_{\delta}$, and we use statement (\ref{if.delta}), 
	equation (\ref{Derivative}) and the mean value theorem, 
	in order to compute for any $\eta \in B$ whose first component $\bar \eta \in C^{0}([0,T];C^{0,\alpha}(\Sigma,\rel^{6n}))$ satisfies $\parallel \bar \eta \parallel_{C^{0}([0,T];
	C^{0,\alpha}(\Sigma,\rel^{6n}))} \leq \epsilon$:
	\begin{eqnarray}  \label{pointwise.estim.F.1}
		|M^{F_0}(x,\bar h_t(x) + \bar \eta_t(x)) \cdot \eta_t^*(x) 
		- M^{F_0}(x,\bar h_t(x)) \cdot \eta_t^*(x) |    \\
		\leq  C(\Sigma,F_0,\bar h,\epsilon) \,|\bar \eta_t(x)| \,|\eta^*_t(x)|   \nonumber
	\end{eqnarray} 
	$\forall \, (x,t)\in \Sigma \times [0,T]$ and also:
	\begin{eqnarray}   \label{pointwise.estim.F.2}
		|D_{(v,w)}(M^{F_0}(x,\bar h_t(x) + \bar \eta_t(x)) \cdot h_t^*(x)) 
		- D_{(v,w)}(M^{F_0}(x,\bar h_t(x)) \cdot h^*_t(x))|           \qquad          \\ 
		\leq  C(\Sigma,F_0,\bar h,\epsilon) \,|\bar \eta_t(x)| \,|h^*_t(x)|  \qquad \nonumber
	\end{eqnarray}
	$\forall \, (x,t)\in \Sigma \times [0,T]$, provided  
	$\epsilon =\epsilon(F_0,\bar h)>0$ is sufficiently small. 
	Now we obtain by means of an integration of estimates 
	(\ref{pointwise.estim.F.1}) and (\ref{pointwise.estim.F.2}) 
	the estimates 
	\begin{eqnarray}  \label{Lp.estim.F.1}
		\parallel M^{F_0}(\,\cdot \,,\bar h + \bar \eta)\cdot \eta^* - M^{F_0}(\,\cdot\,,\bar h)\cdot \eta^* \parallel_{L^{p}(\Sigma \times [0,T],\rel^{n})}   \\ 
		\leq  C(\Sigma,F_0,\bar h,\epsilon) \,
		\parallel \eta^* \parallel_{L^{p}(\Sigma \times [0,T],\rel^{8n})}  \, \parallel \bar \eta \parallel_{L^{\infty}(\Sigma \times [0,T],\rel^{6n})}                 \nonumber
	\end{eqnarray} 
	and
	\begin{eqnarray}  \label{Lp.estim.F.2}
		\parallel D_{(v,w)}(M^{F_0}(\,\cdot\,,\bar h + \bar \eta) \cdot h^*) - D_{(v,w)}(M^{F_0}(\,\cdot\,,\bar h) \cdot h^*)
		\parallel_{L^{p}(\Sigma \times [0,T],\rel^{6n^2})}     \qquad  \\
		\leq  C(\Sigma,F_0,\bar h,\epsilon) \,\parallel h^* \parallel_{L^{p}(\Sigma \times [0,T],\rel^{8n})} \,\parallel \bar \eta \parallel_{L^{\infty}(\Sigma \times [0,T],\rel^{6n})},       \nonumber
	\end{eqnarray}  
	for any $\eta \in B$ whose first component 
	$\bar \eta$ satisfies $\parallel \bar \eta \parallel_{C^{0}([0,T];C^{0,\alpha}(\Sigma,
	\rel^{6n}))} \leq \epsilon$. 
    Moreover, we infer from the fact that 
	$\bar h \in \V_{\delta}$ and from the definition of $\V_{\delta}$:       
	\begin{eqnarray}  \label{Lp.estim.F.3}
		\parallel D_{(v,w)}(M^{F_0}(\,\cdot\,,\bar h
		+\bar \eta) \cdot \eta^*) 
		\parallel_{L^{p}(\Sigma \times [0,T],\rel^{6n^2})}       \\
		\leq C(\Sigma,F_0,\bar h,\epsilon) \,
		\parallel \eta^* \parallel_{L^{p}(\Sigma \times [0,T],\rel^{8n})}.            \nonumber
	\end{eqnarray}   
	Now we derive from the linearity of the map 
	$[(v,w,y) \mapsto D_{(v,w)}(M^{F_0}(x,v,w) \cdot y)]$ w.r.t. $y$ the identity
	\begin{eqnarray}  \label{tricky} 
		D_{(v,w)}(M^{F_0}(x,\bar h + s\,\bar \eta) \cdot 
		(h^* +s \eta^*)) 
		-D_{(v,w)}(M^{F_0}(x,\bar h) \cdot h^*)    \qquad            \\
		= D_{(v,w)}(M^{F_0}(x,\bar h + s\,\bar \eta) \cdot h^*) 
		-D_{(v,w)}(M^{F_0}(x,\bar h) \cdot h^*)    \nonumber   \\
		+ s\, D_{(v,w)}(M^{F_0}(x,\bar h + s\,\bar \eta) \cdot \eta^*),                  \nonumber         
	\end{eqnarray}
	and we infer from equation (\ref{Derivative}), from estimates (\ref{Lp.estim.F.1})--(\ref{Lp.estim.F.3}), from equation (\ref{tricky}) and from the mean value theorem that the operator 
	\begin{equation}  \label{F.sharp}
	\F^{F_0}_{\sharp}:\V_{\delta} \times L^{p}([0,T];W^{1,p}(\Sigma,\rel^{8n}))
	\subset B \longrightarrow  L^p([0,T];L^p(\Sigma,\rel^n))
	\end{equation} 
	mapping $h=(\bar h,h^*) \mapsto 
	\F^{F_0}(\,\cdot \,,\bar h,h^*)$, satisfies:
	\begin{eqnarray*} 
		\parallel  \F^{F_0}_{\sharp}(h+\eta) - \F^{F_0}_{\sharp}(h) 
		- D_{(v,w,y)}\F^{F_0}(\,\cdot\,,h) \cdot \eta  \parallel_{L^{p}([0,T];L^{p}(\Sigma,\rel^{n}))}  \\	
		= \parallel \int_0^1 \big{(} 
		D_{(v,w,y)}\F^{F_0}(\,\cdot\,,h + s\,\eta) - 
		D_{(v,w,y)}\F^{F_0}(\,\cdot\,,h) \big{)} \, ds \cdot \eta \parallel_{L^{p}([0,T];L^{p}(\Sigma,\rel^{n}))} \\
		\leq \parallel \int_0^1
		\Big{[} D_{(v,w)}(M^{F_0}(\,\cdot\,,\bar h + s\,\bar \eta) \cdot h^*) - D_{(v,w)}(M^{F_0}(\,\cdot\,,\bar h) \cdot h^*)     \\ 
		+ s\, D_{(v,w)}(M^{F_0}(\,\cdot\,,\bar h + s\,\bar \eta) \cdot \eta^*) \Big{]} \, ds 
		\cdot \bar \eta \parallel_{L^{p}([0,T];L^{p}(\Sigma,\rel^{n}))}       \\
		+ \parallel \int_0^1 \Big{[} M^{F_0}(\,\cdot\,,\bar h+s \bar \eta)                 
		- M^{F_0}(\,\cdot\,,\bar h)  \Big{]} ds \cdot \eta^*  \, ds
		\parallel_{L^{p}([0,T];L^{p}(\Sigma,\rel^{n}))}         \\
		+  \parallel \int_0^1 \Big{[} D_{(v,w)}\tilde \F^{F_0}(\,\cdot\,,\bar h+ s \,\bar \eta) 
		-  D_{(v,w)}\tilde \F^{F_0}(\,\cdot\,,\bar h) \Big{]} \,ds \cdot \bar \eta \parallel_{L^{p}([0,T];L^{p}(\Sigma,\rel^{n}))}     \\
		\leq C(\Sigma,F_0,\bar h,\epsilon) \, 
		\Big{[} \parallel h^* \parallel_{L^{p}([0,T];L^{p}(\Sigma,\rel^{8n}))} \,    \parallel \bar \eta \parallel_{L^{\infty}(\Sigma 
		\times [0,T],\rel^{6n})}^2                  \\
		+ \parallel \bar \eta \parallel_{L^{\infty}(\Sigma \times [0,T],\rel^{6n})} \,
		\parallel \eta^* \parallel_{L^{p}([0,T];L^{p}(\Sigma,\rel^{8n}))}    \\
		+ \parallel \eta^* \parallel_{L^{p}([0,T];L^{p}(\Sigma,\rel^{8n}))}  \,         \parallel \bar \eta \parallel_{L^{\infty}(\Sigma \times [0,T],\rel^{6n})}      \\       
		+ \parallel \bar \eta \parallel_{L^{p}([0,T];L^{p}(\Sigma,\rel^{6n}))} 
        \, \parallel \bar \eta \parallel_{L^{\infty}(\Sigma \times [0,T],\rel^{6n})}\Big{]}    \\
		\leq  \hat C(\Sigma,T,F_0,h,\epsilon,p) \, 
		\parallel \eta \parallel_B^2,
	\end{eqnarray*} 
	for any $\eta \in B$ whose first component 
	$\bar \eta \in C^{0}([0,T];C^{0,\alpha}(\Sigma,\rel^{6n}))$ satisfies $\parallel \bar \eta \parallel_{C^{0}([0,T];C^{0,\alpha}(\Sigma,\rel^{6n}))}\leq \epsilon$. This estimate proves that the operator $\F^{F_0}_{\sharp}$ 
	in line (\ref{F.sharp}) is Fr\'echet differentiable in any point $h \in \V_{\delta} \times L^{p}([0,T];W^{1,p}(\Sigma,\rel^{8n}))$ and that the Fr\'echet derivative of $\F^{F_0}_{\sharp}$ in any such $h \in B$ is concretely given by the formula:
	\begin{eqnarray}  \label{Frechet.derivative.F} 
		D\F^{F_0}_{\sharp}(h).\eta =  D_{(v,w,y)}\F^{F_0}(\,\cdot\,,h) \cdot \eta \quad  \\ 
		= D_{(v,w)}(M^{F_0}(\,\cdot\,,\bar h) \cdot h^*) \cdot \bar \eta
		+ M^{F_0}(\,\cdot\,,\bar h) \cdot \eta^*
		+ D_{(v,w)}\tilde\F^{F_0}(\,\cdot\,,\bar h)\cdot \bar \eta,  \nonumber
	\end{eqnarray} 
	where we have used equation (\ref{Derivative}) in the second equality. Moreover, one can easily deduce from formulae (\ref{pointwise.estim.F.1}),
	(\ref{pointwise.estim.F.2}) and (\ref{Frechet.derivative.F}), that the Fr\'echet derivative 
	$D\F^{F_0}_{\sharp}(h): B \longrightarrow  L^{p}([0,T];L^p(\Sigma,\rel^{n}))$ 
	depends continuously on $h \in B$. Furthermore, since we had taken the exponent $\alpha=\alpha(p) \in (0,1)$ from embedding (\ref{C0.Hoelder.embedding.XT}), we can now use this embedding, in order to see that the linear map 
	$$ 
	X_{T,p} \ni [(x,t) \mapsto f_t(x)] \stackrel{\Lift}\mapsto 
	[(x,t) \mapsto (D_xf_t(x),D_x^2f_t(x),D^3_xf_t(x))] \in B 
	$$
	is continuous and thus maps a sufficiently small open neighborhood $W_{U_0,T,p}$ of the flow line 
	$\{\PP(t,0,U_0)\circ \phi_t^{U_0,F_0}\}_{t\in [0,T]}$ in $X_{T,p}$ into the product $\V_{\delta} \times L^{p}([0,T];W^{1,p}(\Sigma,\rel^{8n}))$. 
	Hence, the composition 
	\begin{eqnarray*}
		W_{U_0,T,p} \ni [(x,t) \mapsto f_t(x)] \mapsto 
		(\F^{F_0}_{\sharp}\circ \Lift)([(x,t) \mapsto f_t(x)])   \\ 
		\equiv [(x,t) \mapsto \F^{F_0}(x,D_xf_t(x),D_x^2f_t(x),D^3_xf_t(x))]   \in L^p([0,T];L^p(\Sigma,\rel^n))
	\end{eqnarray*}
	has just turned out to be a non-linear $C^1$-operator, and we infer from formula (\ref{Frechet.derivative.F}) that 
	\begin{eqnarray}   \label{Frechet.F}
		D\F^{F_0}(\,\cdot\,,D_xf_t,D_x^2f_t,D^3_xf_t).
		(D_x\eta_t,D_x^2\eta_t,D^3_x\eta_t) \quad \,\,\,\\
		= D_{(v,w)}(M^{F_0}(\,\cdot\,,D_xf_t,D_x^2f_t) \cdot D^3_xf_t) \cdot (D_x\eta_t,D_x^2\eta_t)   \nonumber \\
		+ M^{F_0}(\,\cdot\,,D_xf_t,D_x^2f_t) \cdot D^3_x\eta_t                           
		+D_{(v,w)}\tilde\F^{F_0}(\,\cdot\,,D_xf_t,D_x^2f_t)\cdot (D_x\eta_t,D_x^2\eta_t),  \nonumber
	\end{eqnarray} 
	for every $\{f_t\}_{t \in [0,T]} \in W_{U_0,T,p}$ and $\{\eta_t\} \in T_{f}W_{U_0,T,p} = X_{T,p}$. Repeating the above reasoning starting in line (\ref{B.space}) - using now equations (\ref{A0.squared}) and (\ref{A0_hurra}) instead of equation (\ref{F.third.derivatives}) and again embedding (\ref{C0.Hoelder.embedding.XT}) - and adjusting the definition of the auxiliary Banach space $B$ in (\ref{B.space}) appropriately, one can also show that the operator 
	\begin{eqnarray}  \label{A.def}
		X_{T,p} \supset W_{U_0,T,p} \ni  [(x,t) \mapsto f_t(x)]  \mapsto 
		[(x,t) \mapsto \Area(x,D_xf_t(x),D_x^2f_t(x))]   \nonumber   \\  
		:= [(x,t) \mapsto \frac{1}{2} \, |A^0_{f_t}|^{-4}(x)]     
		\in C^0(\Sigma \times [0,T],\rel)     \qquad
	\end{eqnarray} 
	and also the operator 
	\begin{eqnarray}    \label{Q.def}
		X_{T,p} \supset W_{U_0,T,p} \ni [(x,t) \mapsto f_t(x)]   \mapsto 
		[(x,t) \mapsto  Q(D_xf_t(x),D_x^4f_t(x))]     \qquad                                       \\
		:=[(x,t) \mapsto g^{ij}_{f_t}(x) \, g^{kl}_{f_t}(x) \, \nabla_i^{F_0} \nabla_j^{F_0} 
		\nabla_k^{F_0} \nabla_l^{F_0}(f_t)(x)]  \in  L^p([0,T];L^p(\Sigma,\rel^n))        \nonumber
	\end{eqnarray}   
	are of class $C^1$ on a sufficiently small open neighborhood $W_{U_0,T,p}$ of the flow line $\{\PP(t,0,U_0)\circ \phi_t^{U_0,F_0}\}_{t\in [0,T]}$ in $X_{T,p}$. 
	Hence, both the first and the second component of the operator $\Psi^{F_0,T}:W_{F_0,T,p} \longrightarrow Y_{T,p}$ have turned out to be of class $C^1$. We should mention here already for later use, that the above function 
	\begin{equation}  \label{Q.def.2}
	Q:\rel^{2n} \times \rel^{16n} \longrightarrow \rel^n,   \qquad  (v,p)\mapsto Q(v,p), 
	\end{equation} 
	in (\ref{Q.def}) is rational in its first $2n$ variables $v_1,\ldots,v_{2n}$ and linear in its last $16n$ variables $p_1,\ldots p_{16n}$. This follows immediately from the definition in (\ref{Q.def}) and from the identity 
	\begin{equation}   \label{Cramers.rule}
	(g_f)^{ij} \equiv ((g_{f})_{ij})^{-1} 
	= \frac{1}{\det (g_{f})} \,((g_{f})_{ij})^{\sharp}
	\end{equation} 
	by Cramer's rule, for any $C^1$-immersion $f:\Sigma \longrightarrow \rel^n$. See here also \cite{Jakob_Moebius_2016}, p. 1165.    \\
	\noindent 
	$\beta$) Now we go on proving higher regularity of the operator 
	$\Psi^{F_0,T}:W_{U_0,T,p} \longrightarrow Y_{T,p}$, 
	following the lines of the proof of Proposition 3.5 in 
	Appendix A.1 of \cite{Spener}. To this end, we first recall the special algebraic structure of the function 
	\begin{eqnarray}  \label{map.the.F}
	X_{T,p} \supset W_{U_0,T,p} \ni  [(x,t) \mapsto f_t(x)]  \qquad \nonumber  \\   
    \mapsto [(x,t) \mapsto   \F^{F_0}(x,D_xf_t(x),D_x^2f_t(x),D_x^3f_t(x))]\in L^p([0,T];L^p(\Sigma,\rel^n))        \qquad                      
	\end{eqnarray} 
	in (\ref{phi_revisited}), which is expressed by the affine 
	linear structure of the function 
	$\F^{F_0}(x,v_1,\ldots,v_{2n},w_1,\ldots,w_{4n},
	y_1, \ldots,y_{8n})$ in formula (\ref{F.third.derivatives}) w.r.t. the variable $y \in \rel^{8n}$, 
	representing the third covariant derivatives 
	$\nabla^{F_0}_{111}(f_t)$, $\ldots$, 
	$\nabla^{F_0}_{222}(f_t)$ in 
	equations (\ref{de_Turck_equation_2}) and (\ref{phi_revisited}).
	Now, the function $M^{F_0}:\Sigma \times \rel^{2n} \times \rel^{4n} \longrightarrow \textnormal{Mat}_{n,8n}(\rel)$, i.e.  
	\begin{equation} \label{M.F.0}
		[(x,v_1,\ldots,v_{2n},w_1,\ldots,w_{4n}) \mapsto 
		M^{F_0}(x,v_1,\ldots,v_{2n},w_1,\ldots,w_{4n})],
	\end{equation} 
	and the function $\tilde \F^{F^0}: \Sigma \times \rel^{2n} \times \rel^{4n} \longrightarrow \rel^n$ depend rationally 
	on the $6n$ variables $v_i$ and $w_j$, representing the derivatives 
	$\nabla^{F_0}_{1}(f_t)$,$\nabla^{F_0}_{2}(f_t)$, 
	$\nabla^{F_0}_{11}(f_t)$, $\ldots$, 
	$\nabla^{F_0}_{22}(f_t)$ in equations (\ref{de_Turck_equation_2}) and (\ref{phi_revisited}). 
	We recall here that the embedding (\ref{C0.Hoelder.embedding.XT}) guarantees, that we have 
	for the covariant derivatives up to second order: 
	\begin{equation}    \label{derivatives.first.second.order}
		\nabla^{F_0}_{1}(f_t), \nabla^{F_0}_{2}(f_t), 
		\nabla^{F_0}_{11}(f_t), \ldots, \nabla^{F_0}_{22}(f_t) 
		\in C^{0}(\Sigma \times [0,T],\rel^n),
	\end{equation} 
	whereas there holds trivially: 
	\begin{equation}    \label{derivatives.third.fourth.order}
		\nabla^{F_0}_{111}(f_t), \ldots ,\nabla^{F_0}_{222}(f_t), 
		\nabla^{F_0}_{1111}(f_t), \ldots, \nabla^{F_0}_{2222}(f_t)
		\in L^p([0,T];L^p(\Sigma,\rel^n))
	\end{equation}  
	for any $\{f_t\} \in W_{U_0,T,p}$. 
	Therefore, we may proceed as in the proof of Proposition 3.5 in \cite{Spener}, making use of the fact that
	by Proposition 2.5 in \cite{Shao.2015} the bilinear maps 
	\begin{eqnarray}  \label{bilinear.1}
		C^{0}(\Sigma \times [0,T],\rel) \bigotimes C^{0}(\Sigma \times [0,T],\rel) \longrightarrow 
		C^{0}(\Sigma \times [0,T],\rel),       \\  
		u_1 \otimes u_2 \mapsto  
		[(x,t) \mapsto u_1(x,t) \cdot  u_2(x,t)]      \nonumber
	\end{eqnarray} 
	and 
	\begin{eqnarray}  \label{bilinear.2}
		C^{0}(\Sigma \times [0,T],\textnormal{Mat}_{n,m}(\rel)) 
		\bigotimes L^p([0,T];L^p(\Sigma,\rel^m))
		\longrightarrow L^p([0,T];L^p(\Sigma,\rel^n)), \qquad                               \\  
		(M,Y)   \mapsto   
		[(x,t) \mapsto M(x,t) \cdot Y(x,t)]    \qquad \quad  \nonumber
	\end{eqnarray} 
	are continuous and thus real analytic, for any fixed 
	pair of integers $m,n \in \nat$, and that 
	- similarly to Proposition 6.4 in \cite{Shao.Simonett} - the non-linear operator 
	\begin{eqnarray}  \label{analytic.inversion}
		\{\, u \in C^{0}(\Sigma \times [0,T],\rel) | 
		\inf_{\Sigma \times [0,T]}|u| >\delta \,\}
		\longrightarrow  C^{0}(\Sigma \times [0,T],\rel), 
		\,\,\, u \mapsto \frac{1}{u},  \qquad
	\end{eqnarray} 
	is real analytic, for any fixed $\delta>0$. 
	Now, due to formula (\ref{F.third.derivatives}) the operator in line (\ref{map.the.F}) can be expressed as the 
	sum of the two following operators: \\
	a) The composition of the linear continuous operator
	\begin{eqnarray*}  
	W_{U_0,T,p} \ni [(x,t) \mapsto f_t(x)] 
	\mapsto  [(x,t) \mapsto (x,D_xf_t(x),D_x^2f_t(x))]    \\
	\in C^{0}(\Sigma \times [0,T],\Sigma \times \rel^{6n}) 
	\end{eqnarray*}  
	with the function $[(x,v,w) \mapsto M^{F_0}(x,v_1,\ldots,v_{2n},w_1,\ldots,w_{4n})]$,
	matrix-multiplied by the $\rel^{8n}$-valued function 
	$[(x,t)\mapsto D_x^3f_t(x)] \in 
	L^p([0,T];L^p(\Sigma,\rel^{8n}))$ of all covariant derivatives of third order of $\{f_t\} \in W_{U_0,T,p}$,\\
	b) the composition of the above linear function 
	$[(x,t) \mapsto f_t(x)] \mapsto 
	[(x,t) \mapsto (x,D_xf_t(x),D_x^2f_t(x))]$ 
	with the function
	\begin{equation}  \label{F.0.tilde}
	[(x,v,w) \mapsto \tilde 
	\F^{F_0}(x,v_1,\ldots,v_{2n},w_1,\ldots,w_{4n})].
	\end{equation} 
	Now, taking into account that the functions $M^{F_0}$ and 
	$\tilde \F^{F_0}$ in (\ref{M.F.0}) and (\ref{F.0.tilde}) are rational expressions in the variables $v_1,\ldots,v_{2n}, w_1,\ldots,w_{4n}$ and in the partial derivatives of the smooth immersion $F_0$ up to third order, 
	and recalling the regularity of the covariant derivatives 
	of any fixed $\{f_t\} \in W_{U_0,T,p}$ up to second order 
	respectively up to fourth order in  
	(\ref{derivatives.first.second.order}) and 
	(\ref{derivatives.third.fourth.order}), we can successively apply statements (\ref{bilinear.1})--(\ref{analytic.inversion}), 
	in order to argue as in the proof of Proposition 3.5 in \cite{Spener}, that the non-linear operators composed above in (a) and (b) are real analytic maps from $W_{U_0,T,p}$ to $L^{p}([0,T];L^p(\Sigma,\rel^{n}))$, implying that the non-linear operator in line (\ref{map.the.F}) is actually real analytic from $W_{U_0,T,p}$ to $L^p([0,T],L^p(\Sigma,\rel^n))$, as well. 
	Moreover, we can see that the operator 
	$\{f_t\}  \mapsto \{Q(D_xf_t,D_x^4f_t)\}$
	in formula (\ref{Q.def}) is the composition of the 
	linear continuous operator
	\begin{eqnarray*} 
	W_{U_0,T,p} \ni [(x,t) \mapsto f_t(x)]  
	\mapsto [(x,t) \mapsto (D_xf_t(x),D_x^4f_t(x))]     \\
	\in C^{0}([0,T];C^{1,\alpha}(\Sigma,\rel^{2n})) \times L^p([0,T];L^p(\Sigma,\rel^{16n}))
	\end{eqnarray*} 
	with the function $Q:\rel^{2n} \times \rel^{16n} \to \rel$ from line (\ref{Q.def.2}), which is rational 
	in the variables $v_1,\ldots,v_{2n}$ and linear in the variables $p_1,\ldots,p_{16n}$. We may therefore argue again by means of statements (\ref{bilinear.1})--(\ref{analytic.inversion}), that the non-linear operator in line (\ref{Q.def}) is actually real analytic. Similarly, one can show by means of formulae (\ref{second.fund.i.j}), (\ref{mean_curvat_laplacian}), 
	(\ref{A0.squared}), (\ref{bilinear.1}) and (\ref{analytic.inversion}) that the non-linear operator  
	\begin{eqnarray*} 
		W_{U_0,T,p}  \ni [(x,t) \mapsto f_t(x)] 
		\mapsto [(x,t) \mapsto \Area(x,D_xf_t(x),D_x^2f_t(x))]           \\ 
		\equiv [(x,t) \mapsto \mid A^0_{f_t}(x) \mid^{-4}] \in  C^{0}(\Sigma \times [0,T],\rel) 
		\nonumber
	\end{eqnarray*}  
	from line (\ref{A.def}) is real analytic. Hence, combining these three results we can finally conclude by means of formula (\ref{phi_revisited}) - and again using statement 
	(\ref{bilinear.2}) - that the ``DeTurck-modification'' 
	$$
	\Mill_{F_0}: W_{U_0,T,p} \longrightarrow  L^p([0,T];L^p(\Sigma,\rel^n))
	$$ 
	of the original differential operator 
	$[\{f_t\} \mapsto  \{|A^0_{f_t}|^{-4} \, \nabla_{L^2}\Will(f_t)\}]$ 
	in formulae (\ref{de_Turck_equation}) and (\ref{D_F_0}) is a real analytic operator from $W_{U_0,T,p}$ to 
	$L^p([0,T];L^p(\Sigma,\rel^n))$. Since the operator 
	$\partial_t:X_{T,p} \longrightarrow L^p([0,T];L^p(\Sigma,\rel^n))$ 
	is linear and continuous, we have thus proved the assertion of the third part of the theorem. 
	\item[3)] Just as in the proof of the second part of Theorem 2 in \cite{Jakob_Moebius_2016}, one can infer from a combination of the result of the third part of this theorem with formulae (\ref{phi_revisited}) and (\ref{Frechet.F})--(\ref{Q.def}) 
	and with the chain- and product rule for non-linear $C^1$-operators, that the Fr\'echet derivative 
	of the second component of $\Psi^{U_0,T}$ in any fixed 
	$\{f_t\} \in W_{U_0,T,p}$ is given by the formula:  
	\begin{eqnarray}  \label{D_Psi_f} 
		(D(\Psi^{F_0,T})_2(\{f_t\})).(\{\eta_t\})  
		\equiv  \partial_t \eta_t  -  D(\Mill_{F_0})(\{f_t\}).(\{\eta_t\})   \qquad         \\             
		= \partial_t \eta_t + 
		\frac{1}{2} \mid A^0_{f_t} \mid^{-4} \,
		\big{(} D_{v}Q(D_xf_t,D_x^4f_t) \cdot D_x(\eta_t)   \nonumber \\
		+ g^{ij}_{f_t} \, g^{kl}_{f_t} \, \nabla_i^{F_0} \nabla_j^{F_0} 
		\nabla_k^{F_0} \nabla_l^{F_0}(\eta_t) \big{)}                          \nonumber      \\
		+ D_{(v,w)}\Area(\,\cdot\,,D_xf_t,D_x^2f_t)  \,\cdot \,
		\big{(}D_x\eta_t,D_x^2\eta_t\big{)} \, \, g^{ij}_{f_t} \, g^{kl}_{f_t} \, 
		\nabla_i^{F_0} \nabla_j^{F_0} \nabla_k^{F_0} \nabla_l^{F_0}(f_t)                  \nonumber \\
		+ D_{(v,w,y)}\F^{F_0}(\,\cdot\,,D_xf_t,D_x^2f_t,D_x^3f_t) \cdot 
		\big{(} D_x\eta_t,D_x^2 \eta_t,D_x^3\eta_t \big{)}                                \nonumber \\
		= \partial_t \eta_t +
		\frac{1}{2} \mid A^0_{f_t} \mid^{-4} \,
		g^{ij}_{f_t} \, g^{kl}_{f_t} \, \nabla_i^{F_0} \nabla_j^{F_0} \nabla_k^{F_0} \nabla_l^{F_0}(\eta_t)     
		+ B_{3}^{ijk}(\,\cdot\,,D_xf_t,D_x^2f_t) \cdot \nabla^{F_0}_{ijk}(\eta_t)            \nonumber \\ 
		+B_2^{ij}(\,\cdot\,,D_xf_t,D_x^2f_t,D_x^3f_t,D_x^4f_t) \cdot \nabla^{F_0}_{ij}(\eta_t) \nonumber \\
		+ B_1^i(\,\cdot\,,D_xf_t,D_x^2f_t,D_x^3f_t,D_x^4f_t) \cdot \nabla^{F_0}_{i}(\eta_t) 
		\nonumber
	\end{eqnarray}
	on $\Sigma \times [0,T]$, for $\eta=\{\eta_t\} \in T_{f}W_{U_0,T,p}=X_{T,p}$, 
	where $B_{3}^{ijk}(\,\cdot \,,D_xf_t,D_x^2f_t)$, 
	$B_2^{ij}(\,\cdot\,,D_xf_t,D_x^2f_t,D_x^3f_t,D_x^4f_t)$,
	$B_1^i(\,\cdot\,,D_xf_t,D_x^2f_t,D_x^3f_t,D_x^4f_t)$ are coefficients of $\textnormal{Mat}_{n,n}(\rel)$-valued, contravariant tensor fields of degrees $3,2$ and $1$. It is important to note here, that a quick inspection of formula (\ref{Q.def}) shows, that the gradient $D_{v}Q(v,p)$ of the function $Q$ in line (\ref{Q.def.2}) w.r.t. $v$ depends linearly on the variables $p_1,\ldots,p_{16n}$, just as the function 
	$[(v,p) \mapsto Q(v,p)]$ itself does. See here again 
	\cite{Jakob_Moebius_2016}, p. 1165. Combining this information with formulae (\ref{Frechet.F}), (\ref{A.def}) and (\ref{D_Psi_f}), we see that the coefficients
	$B_1^i(\,\cdot\,,D_xf_t,D_x^2f_t,D_x^3f_t,D_x^4f_t)$ and  $B_2^{ij}(\,\cdot\,,D_xf_t,D_x^2f_t,D_x^3f_t,D_x^4f_t)$ 
	depend affine linearly on the derivatives of 
	third order $\nabla^{F_0}_{ijk}f_t$ and of fourth order $\nabla^{F_0}_{ijkl}f_t$ of any fixed element 
	$\{f_t\} \in W_{U_0,T,p}$. Now, since the derivatives of third and fourth order $\nabla^{F_0}_{ijk}f_t$ and 
	$\nabla^{F_0}_{ijkl}f_t$ are of class $L^p([0,T];L^p(\Sigma,\rel^n))$ by  (\ref{derivatives.third.fourth.order}), and since all the remaining factors respectively summands
	in $B_2^{ij}(\,\cdot\,,D_xf_t,D_x^2f_t,D_x^3f_t,D_x^4f_t)$ and $B_1^{i}(\,\cdot\,,D_xf_t,D_x^2f_t,D_x^3f_t,D_x^4f_t)$ are of class $C^{0}([0,T];C^{0,\alpha}(\Sigma))$ on account of formulae (\ref{F.third.derivatives}), (\ref{Frechet.F}) and (\ref{D_Psi_f}) and on account of embedding (\ref{C0.Hoelder.embedding.XT}), 
	these two types of coefficients are of class $L^{p}([0,T];L^{p}(\Sigma,\textnormal{Mat}_{n,n}(\rel)))$, provided the open neighborhood $W_{U_0,T,p}$ of the fixed flow line $\{\PP(t,0,U_0)\circ \phi_t^{U_0,F_0}\}_{t\in [0,T]}$ in $X_{T,p}$ is sufficiently small. Compare here also with the proof of Lemma 3.4 in \cite{Spener}.
    Moreover, using the fact that the map 
	$M^{F_0}:\Sigma \times \rel^{2n} \times \rel^{4n}  
	\longrightarrow \textnormal{Mat}_{n,8n}(\rel)$
	in (\ref{M.F.0}) depends smoothly on 
	$x$ and rationally on the $6n$ variables $v_i$ and $w_j$, 
	we conclude that the coefficients 
	$B_{3}^{ijk}(\,\cdot \,,D_xf_t,D_x^2f_t)$ are of class $C^{0,\alpha}([0,T],C^{0,\alpha}(\Sigma,
	\textnormal{Mat}_{n,n}(\rel)))$, for some appropriate $\alpha=\alpha(p)\in (0,1)$, on account of formulae (\ref{Frechet.F}) and (\ref{D_Psi_f}), 
	and on account of the continuous embedding
	\begin{equation}  \label{Hoelder.embedding.XT} 
	X_{T,p} \hookrightarrow C^{0,\alpha}([0,T];C^{2,\alpha}(\Sigma,\rel^n)),
	\end{equation} 
    for some small $\alpha\in (0,1)$ only depending on  
    $p$, provided $p\in (3,\infty)$.
        \footnote{Embedding (\ref{Hoelder.embedding.XT}) 
        obviously improves embedding
        (\ref{C0.Hoelder.embedding.XT}) 
        - at least if $T$ remains fixed - and it can be proved exactly as the statement of the third part of Lemma 3.3 in \cite{Spener} by means of identification of fractional Sobolev spaces with appropriate Besov spaces and identification of certain Besov spaces 
        with appropriate H\"older-Zygmund spaces, 
    	taken from \cite{Amann.2} and \cite{Triebel.2}, 
    	furthermore by means of continuous embeddings 
    	between appropriate Besov spaces in \cite{Amann.2} and by means of appropriate interpolation theorems, taken from \cite{Amann.2} and \cite{Lunardi.interpolation}.} 
	Using again embedding (\ref{Hoelder.embedding.XT}), one can see immediately in formula (\ref{D_Psi_f}), that the leading term
	\begin{equation}  \label{leading.operator} 
		\{\eta_t\}  \mapsto \Big{\{}\partial_t \eta_t + \frac{1}{2} \mid A^0_{f_t} \mid^{-4} \,
		g^{ij}_{f_t} \, g^{kl}_{f_t} \, \nabla_i^{F_0} \nabla_j^{F_0} 
		\nabla_k^{F_0} \nabla_l^{F_0}(\eta_t) \Big{\}}
	\end{equation} 
	of the second component of the Fr\'echet derivative $D\Psi^{F_0,T}(\{f_t\})$ 
	is a uniformly parabolic linear differential operator of fourth order in the sense of Proposition \ref{A.priori.estimate} below, mapping $X_{T,p}$ into $L^p([0,T];L^p(\Sigma,\rel^n))$, for any fixed $\{f_t\} \in W_{U_0,T,p}$, whose coefficients are of class $C^{0,\alpha}([0,T],C^{0,\alpha}(\Sigma,\rel))$, 
	and it obviously acts on each component of $\{\eta_t\} \in X_{T,p}$ separately. This proves all assertions of the 
	fourth part of the theorem. 
	\item[4)] First of all, the coefficients of the leading term (\ref{leading.operator}) of the second component of $D\Psi^{F_0,T}(\{f_t\})$ are of class $C^{2}(\Sigma \times [0,T],\rel)$, if $\{f_t\}\in C^{4}(\Sigma \times [0,T],\rel^n) \cap W_{U_0,T,p}$. Hence, they satisfy the regularity requirements of Proposition \ref{A.priori.estimate} below. Moreover, the leading operator in (\ref{leading.operator}) is obviously 
	uniformly elliptic and diagonal on $\Sigma \times [0,T]$ 
	in the sense of Proposition \ref{A.priori.estimate} below for $g:=g_{F_0}$, since the functions $f_t:\Sigma \longrightarrow \rel^n$ have to stay umbilic-free immersions for every $t\in [0,T]$, if the open neighborhood $W_{U_0,T,p}$ of the flow line 
	$\{\PP(t,0,U_0)\circ \phi_t^{U_0,F_0}\}_{t\in [0,T]}$ in $X_T$ is sufficiently small, taking here embedding (\ref{Hoelder.embedding.XT}) again into account. Moreover, the remaining coefficients 
	$B_{3}^{ijk}(\,\cdot \,,D_xf_t,D_x^2f_t)$, 
	$B_2^{ij}(\,\cdot\,,D_xf_t,D_x^2f_t,D_x^3f_t,D_x^4f_t)$ and
	$B_1^{i}(\,\cdot\,,D_xf_t,D_x^2f_t,D_x^3f_t,D_x^4f_t)$ of 
	the second component of $D\Psi^{F_0,T}(\{f_t\})$ in 
	(\ref{D_Psi_f}) are of class $C^{0}(\Sigma \times [0,T],\textnormal{Mat}_{n,n}(\rel))$, and thus their restrictions to $\Omega_{\kappa} \times [0,T]$, 
	for any fixed coordinate patch $\Omega_{\kappa} \subset \Sigma$ of a finite atlas of $\Sigma$ as in Section 10 of \cite{Amann.Hieber.Simonett}, satisfy the conditions of Proposition \ref{A.priori.estimate} below. Hence, we infer from this proposition, that both the linear operator 
	$$
	\partial_t - D(\Mill_{F_0})(\{f_t\}): 
	X_T^0:=\{\{u_t\} \in X_{T,p} \,|\, u_0 \equiv 0 \,\}
	\stackrel{\cong}\longrightarrow L^p([0,T];L^p(\Sigma,\rel^n))
	$$ 
	and the Fr\'echet derivative 
	$D\Psi^{F_0,T}(\{f_t\}):X_{T,p} \stackrel{\cong}\longrightarrow Y_{T,p}$ 
	of the non-linear operator $\Psi^{F_0,T}$ 
	in (\ref{Psi}) are isomorphisms, in any 
	$\{f_t\} \in C^{4}(\Sigma \times [0,T],\rel^n) \cap W_{U_0,T,p}$. 
	\qed
\end{itemize} 
\noindent
In order to justify the proof of the fourth part of Theorem 
\ref{Psi.of.class.C_1}, we combine results due to Amann, Duong, Hieber and Simonett in \cite{Amann.1995}, 
\cite{Amann.2004}, \cite{Amann.Hieber.Simonett} and \cite{Duong.Simonett.1997}, yielding the following general statement.
\begin{proposition}    \label{A.priori.estimate}
	Let $\Sigma$ be a smooth, compact, orientable surface without boundary, $g \in \Gamma(T^{0,2}\Sigma)$ a $C^{\infty}$-smooth 
	Riemannian metric on $\Sigma$, $\nabla^g$ the 
	corresponding Levi-Civita-connection on $\Sigma$, $n\in \nat$, $p \in (1,\infty)$, and 
	$T>0$ arbitrarily fixed, and let $\psi:O \stackrel{\cong} \longrightarrow \Sigma'$ be a smooth chart of an arbitrary coordinate patch $\Sigma'$ of $\Sigma$, yielding partial derivatives $\partial_m$, $m=1,2$, and the coefficients 
	$g_{ij}:= g(\partial_i,\partial_j)$ of the metric $g$, 
	restricted to $\Sigma'$. Moreover, let  
	$$
	\partial_t + L:X_{T,p} \longrightarrow L^p([0,T];L^p(\Sigma,\rel^n))
	$$
	be a linear differential operator of order $4$, whose leading operator of fourth order is diagonal, i.e. 
	acts on each component of $f$ separately:
	\begin{eqnarray}  \label{L}
		(\partial_t + L)(f)(x,t) := \partial_t(f_t)(x)  
		+ A_{4}^{ijkl}(x,t)\, \nabla^g_{ijkl}(f_t)(x) \qquad  \\
		+ A_{3}^{ijk}(x,t)\, \nabla^g_{ijk}(f_t)(x) 
		+ A_2^{ij}(x,t)\, \nabla^g_{ij}(f_t)(x)  \qquad \nonumber \\
		+ A_1^i(x,t) \, \nabla^g_{i}(f_t)(x) + 
		A_0(x,t) \,(f_t)(x),  \qquad \nonumber 
	\end{eqnarray}
	in every pair $(x,t) \in \Sigma' \times [0,T]$, and $L$ has to satisfy the following ``structural hypotheses'':
	\begin{itemize} 
		\item[1)] $A_{4}^{ijkl}$, $A_{3}^{ijk}$, 
		$A_2^{ij}$, $A_1^i$, $A_0$ are the coefficients of contravariant and continuous tensor fields $A_4$, $A_3$, $A_2$, $A_1$, $A_0$ on $\Sigma \times [0,T]$ of degrees $4,3,\ldots,0$. Moreover, the tensor $A_{4}$ has to be the square $E \otimes E$ of a contravariant 
		real-valued symmetric tensor field $E$ of order $2$, i.e. $A_{4}^{ijkl}(x,t) = E^{ij}(x,t) \, E^{kl}(x,t)$ with $E^{ij}(x,t)=E^{ji}(x,t) \in \rel$, for every $(x,t) \in \Sigma' \times [0,T]$. 
		\item[2)] The tensor field $E$ is required to be positive definite on $\Sigma \times [0,T]$, i.e. there has to be a number $\Lambda \geq 1$ such that there holds:
		\begin{equation} \label{uniformly_parabolic} 
		E^{ij}(x,t)\,  \xi_i \xi_j  
		\geq \Lambda^{-1/2}\, g^{ij}(x)\, \xi_i \xi_j  
		\end{equation} 
		for any vector $\xi =(\xi_1,\xi_2) \in \rel^2$ and for any pair $(x,t) \in \Sigma' \times [0,T]$.
		\item[3)] There holds 
		$\sup_{\Sigma \times [0,T]} |A_{r}| 
		\leq \Lambda$, for $r=0,1,2,3,4$.
	\end{itemize}
	Then the following three statements hold: 
	\begin{itemize} 
	\item [1)] For every $\theta \in (0,\pi/2)$ there are numbers $\omega=\omega(\Sigma,\Lambda,p,\theta,n)>0$ and $M=M(\Sigma,\Lambda,p,\theta,n) \geq 1$, such that for each fixed $t^*\in [0,T]$ the elliptic operator   
	\begin{eqnarray} 
	    \omega\, \textnormal{Id}_{L^{p}(\Sigma,\rel^n)} + L_{t^*}:= 
		\omega \,\textnormal{Id} _{L^{p}(\Sigma,\rel^n)} 
		+ \big{(} A_{4}^{ijkl}(\,\cdot\,,t^*)\, \nabla^g_{ijkl} + A_{3}^{ijk}(\,\cdot\,,t^*)\, \nabla^g_{ijk}    \nonumber \\ 
		+ A_2^{ij}(\,\cdot\,,t^*)\, \nabla^g_{ij}  +      
		A_1^i(\,\cdot\,,t^*)\, \, \nabla^g_{i} +  A_0(\,\cdot\,,t^*)\, \big{)}    \qquad
		\end{eqnarray}  
		in $L^{p}(\Sigma,\rel^n)$ maps 
		$W^{4,p}(\Sigma,\rel^n)$ isomorphically onto $L^{p}(\Sigma,\rel^n)$ and is of type $(M,\pi-\theta)$ in the sense of formula (1.1) in \cite{Amann.Hieber.Simonett} respectively of formula 
		(4.7.8) in Chapter III of \cite{Amann.1995}.  
		\item [2)] The operator 
		$$
		(\gamma_0,\partial_t + L): X_{T,p} \stackrel{\cong} \longrightarrow Y_{T,p} \equiv W^{4-\frac{4}{p},p}(\Sigma,\rel^n) \times L^p([0,T],L^p(\Sigma,\rel^n))   
		$$
		is a topological isomorphism, where 
		$\gamma_0: X_{T,p} \longrightarrow W^{4-\frac{4}{p},p}(\Sigma,\rel^n)$ 
		is the trace operator at time $t=0$, mapping 
		$\{u_t\} \mapsto u_0$; see also formula (\ref{Trace.X.T}).   
		\item [3)] The parabolic linear operator   
		\begin{equation}\label{global.initial.value.problem}   \partial_t + L: 
		X_T^0:=\{\{u_t\} \in X_{T,p} \,|\, \gamma_0(\{u_t\}) \equiv 0 \,\}   \stackrel{\cong}\longrightarrow L^p([0,T],L^p(\Sigma,\rel^n))
		\end{equation} 
		is a topological isomorphism. 	 
	\end{itemize}
\end{proposition}  
\proof
We fix an arbitrary uniformly parabolic operator $L$ as in 
(\ref{L}) satisfying the structural hypotheses (1), (2) and 
(3) of this proposition, and we choose a finite atlas $\{(\Omega_{\kappa},\psi_{\kappa})\}_{\kappa=1,\ldots,K}$ 
of smooth charts for $\Sigma$ as in Section 10 of \cite{Amann.Hieber.Simonett}.
Since the coefficients of $L$ are assumed to be continuous, 
we can consider - instead of the original operator 
$\partial_t + L$ in (\ref{L}) - for every fixed $t^*\in [0,T]$ the linear, uniformly elliptic differential operator 
\begin{eqnarray}   \label{L.t}
	\eta \mapsto L_{t^*}(\eta)
	:= \big{(} A_{4}^{ijkl}(\,\cdot\,,t^*)\, \nabla^g_{ijkl} + A_{3}^{ijk}(\,\cdot\,,t^*)\, \nabla^g_{ijk}               +    \nonumber    \\
	+ A_2^{ij}(\,\cdot\,,t^*)\, \nabla^g_{ij} + A_1^i(\,\cdot\,,t^*)\, \, \nabla^g_{i} +  A_0(\,\cdot\,,t^*)\, \big{)}(\eta)
\end{eqnarray}  
in $L^p(\Sigma,\rel^n)$ with domain $W^{4,p}(\Sigma,\rel^n)$. 
Now, by the second assumption of this proposition 
the highest order term 
\begin{equation}  \label{highest.order} 
	W^{4,p}(\Sigma,\rel^n)  \ni \eta  \mapsto 
	\Area_{t^*}(\eta):= E^{ij}(\,\cdot\,,t^*)\, E^{kl}(\,\cdot\,,t^*)\, \nabla^g_{ijkl}(\eta) \in L^{p}(\Sigma,\rel^n)    
\end{equation} 
of the linear operator $L_{t^*}$ in line (\ref{L.t}) 
is uniformly elliptic and diagonal, in every fixed $t^*\in [0,T]$. Hence, the spectrum of its ``symbol'', which is 
given by the family of linear maps 
\begin{equation}   \label{symbol.A}
	(\Area_{t^*})_{\sharp}(x,\xi):= \textnormal{Id}_{\com^n} \cdot
	\Big{(}\sum_{i,j,k,l} E^{ij}(x,t^*)\, E^{kl}(x,t^*) \, 
	\, \xi_i \,\xi_j \,\xi_k \,\xi_l \Big{)}:\com^n \longrightarrow \com^n,
\end{equation}  
for $x \in \Omega_{\kappa}$, for every $\xi \in \rel^2$ 
with $|\xi|= 1$ and $\kappa \in \{1,\ldots,K\}$, is 
precisely given by:
$$
	\sigma((\Area_{t^*})_{\sharp}(x,\xi)) = 
	\Big{\{} \sum_{i,j,k,l} E^{ij}(x,t^*) \, E^{kl}(x,t^*) \, 
	\xi_i \,\xi_j \,\xi_k \,\xi_l  \Big{\}}
$$  
for every $x\in \Omega_{\kappa}$, for
$\xi \in \rel^2$ with $|\xi|= 1$ and for $\kappa \in \{1,\ldots,K\}$, and these are real positive numbers, being uniformly bounded from below in terms of $\Lambda^{-1}$, on account of the required uniform positive definiteness of the tensor field $E$ on $\Sigma \times [0,T]$ by assumption (\ref{uniformly_parabolic}). In particular, the ``principal symbol'' of the operator $L_{t^*}$ in line (\ref{L.t}) is 
``parameter elliptic'' in every fixed point 
$x \in \Sigma$ with ``angle of ellipticity''$=0$ in the 
sense of Section 10 in \cite{Amann.Hieber.Simonett} 
respectively Section 7 in \cite{Duong.Simonett.1997}, 
for every fixed $t^*\in [0,T]$; 
see here also Section 5.1 in \cite{Denk.Hieber.Pruess.1}. 
Since the coefficients of the linear operators 
$\{L_{t}\}_{t \in [0,T]}$ are also assumed to be 
uniformly bounded in $L^{\infty}(\Sigma)$ 
in terms of $\Lambda$, we can thus infer from 
Theorem 10.1 in \cite{Amann.Hieber.Simonett} 
and from Theorem 7.1 in \cite{Duong.Simonett.1997}, 
that for any fixed $\theta \in (0,\pi/2)$ there are numbers $\omega=\omega(\Sigma,\Lambda,p,\theta,n)>0$ and 
$M=M(\Sigma,\Lambda,p,\theta,n) \geq 1$, such that 
\begin{itemize}
	\item [1)] the linear operator 
	$\omega\, \textnormal{Id}_{L^{p}(\Sigma,\rel^n)} + L_{t^*}$ in $L^p(\Sigma,\rel^n)$ yields an isomorphism  
	$$ 
	\omega\, \textnormal{Id}_{L^{p}(\Sigma,\rel^n)} + L_{t^*}: W^{4,p}(\Sigma,\rel^n) \stackrel{\cong}\longrightarrow L^p(\Sigma,\rel^n),
	$$ 
	\item [2)] the linear operator 
	$\omega\, \textnormal{Id} _{L^{p}(\Sigma,\rel^n)} + L_{t^*}$ is of positive type $(M,\pi-\theta)$, which means 
	precisely that the inverse
	$$
	\big{(} (\omega  + \lambda) \,
	\textnormal{Id}_{L^{p}(\Sigma,\rel^n)} 
	+ L_{t^*}  \big{)}^{-1}: L^{p}(\Sigma,\rel^n) 
	\longrightarrow L^{p}(\Sigma,\rel^n)
	$$
	exists and is bounded for every complex number $\lambda$, 
	which is contained in the sector 
	$$
	S_{\vartheta}:=\{\,\lambda \in \com^* \,| \,
	|\textnormal{arg}(\lambda)| \leq \vartheta \, \} 
	\cup \{0\}
	$$
	with apex in $0$ and with half opening angle 
	$\vartheta:=\pi-\theta \in \big{(}\frac{\pi}{2},\pi \big{)}$, and that there holds the estimate
	\begin{equation}  \label{sectorial.estimate}
	\parallel \big{(} (\omega + \lambda)\, 
	\textnormal{Id}_{L^{p}(\Sigma,\rel^n)} + L_{t^*} \big{)}^{-1} \parallel_{\Lift(L^{p}(\Sigma,\rel^n),
	L^{p}(\Sigma,\rel^n))}
	\leq \frac{M}{1+|\lambda|}
	\end{equation}
	for every $\lambda \in S_{\vartheta}$, and for every 
	fixed $t^*\in [0,T]$, 
	\item[3)] the linear operator 
	$\omega \,\textnormal{Id}_{L^{p}(\Sigma,\rel^n)} + L_{t^*}$ admits a bounded $\HH_{\infty}$-calculus 
	in $L^p(\Sigma,\rel^n)$ with angle $\theta$,
	which means precisely that for any bounded holomorphic function $f :\mathring{S_{\theta}} \longrightarrow \com$ there holds 
	$$ 
	\parallel f(\omega \,\textnormal{Id} _{L^{p}(\Sigma,\rel^n)} + L_{t^*})  \parallel_{\Lift(L^p(\Sigma,\rel^n),L^p(\Sigma,\rel^n))} \leq M \, \parallel f \parallel_{L^{\infty}(\mathring{S_{\theta}})}, 
	$$  
	in the sense of Sections 2 and 7 in \cite{Duong.Simonett.1997},
	again for some sufficiently large constant $M=M(\Sigma,\Lambda,p,\theta,n)>0$. 
	By formula (2.4) in \cite{Amann.Hieber.Simonett} respectively Corollary 7.2 in \cite{Duong.Simonett.1997} the third conlcusion implies particularly, that 
	\item[4)] the operator $\omega \,\textnormal{Id}_{L^{p}(\Sigma,\rel^n)} + L_{t^*}$ has bounded imaginary powers in $L^p(\Sigma,\rel^n)$ with angle 
	$\theta$, i.e. that there holds for every fixed 
	$t^* \in [0,T]$: 
	$$ 
	\parallel (\omega \,\textnormal{Id}_{L^{p}(\Sigma,\rel^n)} + L_{t^*})^{i\,s}  \parallel_{\Lift(L^p(\Sigma,\rel^n),
	L^p(\Sigma,\rel^n))}\leq M \, e^{\theta \,s}, 
	$$   
	for every $s \in \rel$, and for every fixed 
	$\theta \in (0,\pi/2)$.   
\end{itemize}
On account of our assumption, that the coefficients 
$A_{4}^{ijkl}, \ldots, A_{0}$ 
of the family of linear operators $L=\{L_{t}\}_{t \in [0,T]}$
are continuous on $\Sigma \times [0,T]$ there holds
trivially  
$$ 
[t \mapsto L_t] \in 
C^0([0,T],\Lift(W^{4,p}(\Sigma,\rel^n), 
L^{p}(\Sigma,\rel^n)). 
$$
Hence, combining the above four statements (1)--(4) 
with the facts that 
$L^p(\Sigma,\rel^n)$ is a UMD-space by 
Theorem 4.5.2 in Chapter III of \cite{Amann.1995},
and that the trace operator $\gamma_0$ maps $X_{T,p}$ 
continuously onto $W^{4-\frac{4}{p},p}(\Sigma,\rel^n)$ 
by formula (\ref{Trace.X.T}), we can deduce the second 
and the third assertion of this proposition from a 
combination of Theorem 4.10.8 in Chapter III of \cite{Amann.1995} with Theorem 7.1 in \cite{Amann.2004}.                         
\qed    \\
\noindent    
\\
\underline{Proof of Theorem \ref{short.time.solution}}:
   \begin{itemize} 
   \item[1)] Formulae (\ref{K.F_0.f_t.f_t}) respectively (\ref{A.F_0.f.f}) motivate us, to consider the quasilinear parabolic Cauchy problem 
    \begin{equation}  \label{A.F_0.ft.ft}
    	\partial_tf_t = K^{F_0}(f_t).(f_t)  \quad 
    	\textnormal{on} \quad \Sigma \times [0,T],
    	\quad \textnormal{with} \,\,\, f_0=u_0 
    	\quad \textnormal{on} \quad \Sigma 
    \end{equation}    
     in the space $X_{T,p}$, instead of Cauchy problem (\ref{initial.value.problem}). Now, on account of embedding (\ref{embedding.Sobol.Hoelder}) 
    we can choose some small open neighbourhood  
    $V_{U_0}$ of the umbilic-free immersion $U_0$ in 
    $W^{4-\frac{4}{p},p}(\Sigma,\rel^n)$ which 
    is still contained in $\textnormal{Imm}_{\textnormal{uf}}(\Sigma,\rel^n)$,  
    and we can infer from Lemma \ref{preparation}, that the non-linear operator $K^{F_0}$ in 
    formula (\ref{A.F_0.ft.ft}) is of class 
    $C^{0,1}(V_{U_0},\Lift(W^{4,p}(\Sigma,\rel^n), 
    L^{p}(\Sigma,\rel^n))$, 
    and that especially the linear differential operator 
    $-K^{F_0}(U_0)$ is of maximal $L^p$-regularity on 
    $[0,T]$ w.r.t. $(W^{4,p}(\Sigma,\rel^n), 
    L^{p}(\Sigma,\rel^n))$. We may therefore apply  
    Theorem 2.1 - here with interpolation parameter $\mu=1$ - 
    and Remark 2.3 in \cite{Koehne.Pruess.Wilke.2010} to Cauchy 
    problem (\ref{A.F_0.ft.ft}), and we obtain from this theorem  
    some small positive time $T=T(U_0)>0$ 
    and an $\varepsilon = \varepsilon(U_0)>0$, such that 
    for any immersion $u_0 \in V_{U_0}$ with 
    \begin{equation}  \label{small.initial}
   	\parallel u_0 - U_0 
	\parallel_{W^{4-\frac{4}{p},p}(\Sigma,\rel^n)}<\varepsilon
	\end{equation}
	the quasilinear parabolic Cauchy problem 
	(\ref{A.F_0.ft.ft}) has exactly one solution 
	$[t \mapsto u(t,u_0)]$ in $X_{T(U_0),p}$.  
	Now, on account of the evolution equation in 
	(\ref{A.F_0.ft.ft}) and on account of embedding 
	(\ref{C0.Hoelder.embedding.XT}) the unique short-time 
	solution $[t \mapsto u(t,u_0)]$ of Cauchy problem (\ref{A.F_0.ft.ft}) in $X_{T(U_0),p}$ has to consist of 
	umbilic-free immersions only. Hence, formula 
	(\ref{A.F_0.f.f}) tells us, that the function
	$[t \mapsto u(t,u_0)]$ also solves the original Cauchy problem 
	(\ref{initial.value.problem}) in $X_{T(U_0),p}$, 
	for every initial immersion $u_0 \in V_{U_0}$ satisfying 
	condition (\ref{small.initial}). Now vice versa,
	any short-time solution of the original Cauchy problem 
	(\ref{initial.value.problem}) in $X_{T^*(u_0),p}$, 
	for some sufficiently small $T^*(u_0)>0$, has to consist 
	of only umbilic-free immersions - for the same reason as 
	above - and therefore any such short-time solution is  
	a short-time solution of Cauchy problem (\ref{A.F_0.ft.ft})
	in $X_{T^*(u_0),p}$ as well, for any fixed initial 
	immersion $u_0 \in V_{U_0}$ satisfying 
	condition (\ref{small.initial}). Therefore, the above
	unique short-time solution $[t \mapsto u(t,u_0)]$ 
	of Cauchy problem (\ref{A.F_0.ft.ft}) in $X_{T(U_0),p}$ 
	also solves Cauchy problem (\ref{initial.value.problem}) 
	\underline{uniquely} in $X_{T(U_0),p}$, 
	for every initial immersion $u_0 \in V_{U_0}$ satisfying 
	condition (\ref{small.initial}) with $\varepsilon=\varepsilon(U_0)>0$ 
	sufficiently small, which proves the assertion. 
	\item[2)] The proof of the second part of the 
	theorem now follows immediately from the result 
	of the first part of this theorem, combined with 
	Lemma \ref{preparation} of this article and 
	Corollary 2.2 and Remark 2.3 of \cite{Koehne.Pruess.Wilke.2010}. 
	\qed	
	\end{itemize} 
\noindent 
Now we finally arrive at the   
\underline{proof of Theorem \ref{real.analytic.flow}}: 
\begin{itemize} 
	\item[1)] Theorem 1 in \cite{Jakob_Moebius_2016} 
	yields the existence of a $C^{\infty}$-smooth short-time solution $\{\PP(t,0,U_0)\}_{t\in [0,\tau]}$ of the MIWF-equation (\ref{Moebius.flow}) starting in $U_0$ for some $\tau>0$, and its proof also yields a smooth family of smooth diffeomorphisms
	$\phi^{U_0,F_0}_t:\Sigma \longrightarrow \Sigma$ with 
	$\phi^{U_0,F_0}_0=\textnormal{id}_{\Sigma}$, for 
	$t\in [0,\tau]$, such that the composition 
	$\{\PP(t,0,U_0) \circ 
	\phi^{U_0,F_0}_t\}_{t\in [0,\tau]}$ is a 
	$C^{\infty}$-smooth solution 
	of the modified evolution equation 
	(\ref{de_Turck_equation}), i.e. it solves Cauchy problem
	(\ref{initial.value.problem}) with $u_0=U_0$ in 
	the class $C^{\infty}(\Sigma \times 
	[0,\tau],\rel^n)$. This short-time solution 
	can be continued to a unique maximal flow line of evolution equation (\ref{de_Turck_equation}) starting 
	in $U_0$, in the sense of Definition 
	\ref{Umbilic.free}, (d). On the other hand 
	this maximal flow line also solves 
	problem (\ref{initial.value.problem}) with $u_0=U_0$ 
	in the space $X_{T^*,p}$, for any 
	$T^* \in (0,T_{\textnormal{max}}(U_0))$,
	and thus coincides with the unique maximal solution 
	$\{\PP^*(t,0,U_0)\}_{t\in [0,t^+(U_0))}$  
	of Cauchy problem (\ref{initial.value.problem})   
	from the second part of Theorem \ref{short.time.solution}
	restricted to $\Sigma \times [0,T^*]$, which implies 
	especially that $T_{\textnormal{max}}(U_0) \leq t^+(U_0)$. 
	Now we show that actually 
	$T_{\textnormal{max}}(U_0)=t^+(U_0)$ holds. 
	To this end, we choose some $T\in (0,t^+(U_0))$, 
	and we recall that the unique solution 
	$\{\PP^*(t,0,U_0)\}_{t\in [0,T]}$  
	of problem (\ref{initial.value.problem}) in 
	$X_{T,p}$ has to be of class 
	$C^{0,\alpha}([0,T];C^{2,\alpha}(\Sigma,\rel^n))$, 
	on account of embedding (\ref{Hoelder.embedding.XT}), 
	for some small $\alpha=\alpha(p) \in (0,1)$.  
	Therefore Lemma \ref{preparation} guarantees us, that 
	evolution equation (\ref{de_Turck_equation}) 
	has the quasilinear structure appearing in formulae (\ref{A.F_0.f.f}) or (\ref{A.F_0.ft.ft}) 
	along its solution $\{\PP^*(t,0,U_0)\}_{t\in [0,T]}$,
    and this fact inspires us, to consider the linear operator 
    \begin{eqnarray}  \label{L}
    L^{F_0,U_0}:=\partial_t - K^{F_0}(\PP^*(\,\cdot\,,0,U_0)):         \\ 
    C^{4+\alpha,1+\frac{\alpha}{4}}(\Sigma \times [0,T],\rel^n)
    \longrightarrow  
    C^{\alpha,\frac{\alpha}{4}}(\Sigma \times [0,T],\rel^n), 
    \nonumber
    \end{eqnarray}  
    which is uniformly parabolic and has coefficients of class 
    $C^{\alpha,\frac{\alpha}{4}}(\Sigma \times [0,T],\rel^n)$, thus satisfies all conditions 
    of Corollary 3 in \cite{Jakob_Moebius_2016} and yields 
    therefore an isomorphism 
    $$ 
    L^{F_0,U_0}: 
    \{\{\eta_t\} \in 
    C^{4+\alpha,1+\frac{\alpha}{4}}(\Sigma \times [0,T],\rel^n) \,|\, 
    \eta_0=U_0 \} \stackrel{\cong}\longrightarrow 
    C^{\alpha,\frac{\alpha}{4}}(\Sigma \times [0,T],\rel^n)    
    $$ 
    between affine linear spaces, again with the exponent $\alpha$ from embedding (\ref{Hoelder.embedding.XT}). 
    Moreover, there holds the linear parabolic differential 
    equation  
    \begin{eqnarray}  \label{linear.quasilinear} 
    L^{F_0,U_0}(\PP^*(t,0,U_0)) 
    = \partial_t(\PP^*(t,0,U_0)) - K^{F_0}(\PP^*(t,0,U_0)).(\PP^*(t,0,U_0)) \nonumber \\
    = 0  \quad \textnormal{for a.e.} \,\,t\in [0,T]  \qquad 
    \end{eqnarray}
    on account of equations (\ref{initial.value.problem}), (\ref{A.F_0.f.f}) and (\ref{L}), 
    and also $\PP^*(0,0,U_0)=U_0$. 
    Hence, since obviously 
    $C^{4+\alpha,1+\frac{\alpha}{4}}(\Sigma \times [0,T],\rel^n)$ is a Banach subspace of $X_{T,p}$,  
    we can therefore conclude as in the 
    proof of Theorem 3 in \cite{Jakob_Moebius_2016}, that 
    the family of immersions 
    $\{\PP^*(t,0,U_0)\}_{t\in [0,T]}$ 
    solves Cauchy problem (\ref{initial.value.problem}) 
    not only in $X_{T,p}$ but actually in the smaller 
    Banach space $C^{4+\alpha,1+\frac{\alpha}{4}}(\Sigma \times [0,T],\rel^n)$. Now we can combine this new 
    information again with the linear parabolic equation 
    (\ref{linear.quasilinear}) and with the smoothness 
    of $U_0$, in order to infer by means of Proposition 3 in \cite{Jakob_Moebius_2016} higher regularity of class 
    $C^{4+k+\alpha,1+\frac{k+\alpha}{4}}
    (\Sigma \times [0,T],\rel^n)$, inductively for 
    every $k \in \nat_0$, just as in the proof of 
    Theorem 3 in \cite{Jakob_Moebius_2016}. 
    Hence, the unique solution 
    $\{\PP^*(t,0,U_0)\}_{t\in [0,T]}$
    of Cauchy problem (\ref{initial.value.problem}) in $X_{T,p}$ 
    has turned out to be of class $C^{\infty}(\Sigma \times [0,T],\rel^n)$, for every $T \in (0,t^+(U_0))$, 
    proving that the maximal solution 
    $\{\PP^*(t,0,U_0)\}_{t\in [0,t^+(U_0))}$ 
    of problem (\ref{initial.value.problem}) from the second 
    part of Theorem \ref{short.time.solution}
    is a maximal flow line of evolution equation (\ref{de_Turck_equation}) 
    in the sense of Definition \ref{Umbilic.free}, (d), 
    and that also $T_{\textnormal{max}}(U_0) \geq t^+(U_0)$ 
    has to hold. 
	\item[2)] Now we consider the unique maximal solution
	$\{\PP^*(t,0,U_0)\}_{t\in [0,T_{\textnormal{max}}(U_0))}$ 
	of problem (\ref{initial.value.problem}) from the first part of this theorem - starting in a smooth and umbilic-free immersion $U_0$ - and we assume here additionally, that $\Sigma$ is a complex compact torus 
	and that the immersion $F_0$ is not only smooth 
	but real analytic on $\Sigma$ w.r.t. the prescribed 
	complex structure. On account of the 
	complicated, quasilinear 
	structure of the differential operator 
	$[f \mapsto \Mill_{F_0}(f)]$ in formula
	(\ref{M.F_0.2}) respectively (\ref{phi_revisited}), 
	we cannot directly apply Simonett's and Shao's 
	regularity theory for linear parabolic equations 
	on uniformly regular Riemannian manifolds \cite{Shao.2015}, \cite{Shao.Simonett} to smooth solutions of our equation (\ref{initial.value.problem}), 
	in order to prove their real analyticity, locally 
	in space and time. 
	But we still succeed here as Shao did in \cite{Shao.2013} or in Section 4 of \cite{Shao.2015} - discussing the real analyticity of the Willmore flow 
    respectively of the Ricci-DeTurck flow - 
 	localizing the quasilinear differential operator 
	in (\ref{A.F_0.f.f}) respectively (\ref{A.F_0.ft.ft}) 
	in terms of a particular 
	partition of unity, then adopting the argument 
	of the proof of Proposition 3.7 in \cite{Shao.2015} 
	via Theorem 4.2 in \cite{Escher.Pruess.Simonett}
	for each localized differential operator,
	and then conclude by means of Proposition 2.5 and 
	Lemma 3.8 in \cite{Shao.2015} and Proposition 6.4 
	in \cite{Shao.Simonett}. Moreover, we will have to work 
	here in slightly different function spaces, namely  
	in little H\"older spaces $h^{s}(\Sigma,\rel^n)$
	instead of $L^p$-spaces, thus exactly following Shao's work in \cite{Shao.2013} or Sections 4--6 of his paper \cite{Shao.2015}. Our argument therefore clearly 
	deviates from the proof of Theorem 3.14 in Appendix A.3 
	in \cite{Spener}, discussing the real analyticity 
	of short-time solutions of the elastic energy flow 
	in $\rel^{m+1}$.   
Henceforth, we will work with a finite conformal atlas $\Area:=\{(O_j,\varphi_j)\}_{j=1,\ldots,N}$ of $\Sigma$ with $\varphi_j:O_j \stackrel{\cong}\longrightarrow B_1^2(0)$, $j=1,\ldots,N$, and invers maps $\psi_j:=(\varphi_j)^{-1}$, 
and we shall consider some arbitrary point $p\in \Sigma$ 
and some coordinate patch $(O_{j_p},\varphi_{j_p})$ 
of $\Area$ containing $p$. 
As explained in Section 3.1 of \cite{Shao.2015}, 
we can add a conformal chart $\varphi_{\iota}:O_{\iota} \stackrel{\cong}\longrightarrow B_1^2(0)$ to the atlas $\Area$, satisfying $\varphi_{\iota}(p)=0$, 
resulting in another conformal atlas 
$\tilde \Area = \Area \cup (O_{\iota},\varphi_{\iota})$.  
Furthermore, we introduce as in Section 3.1 of \cite{Shao.2015} 
or in Section 4 of \cite{Shao.Simonett} open balls 
$B_i:=B^2_{i\varepsilon_0}(0)$, with 
$i=1,\ldots,5$ and $0<\varepsilon_0<\frac{1}{5}$, 
cut-off functions $\chi \in C^{\infty}_c(B_2,[0,1])$ 
and $\varsigma \in C^{\infty}_c(B_5,[0,1])$  
with $\chi \equiv 1$ on $B_1$ and 
$\varsigma \equiv 1$ on $B_4$, the truncated shift
\begin{equation}   \label{theta.mu} 
\theta_{\mu}(y):=y + \mu\, \chi(y), \quad  
\textnormal{for} \quad y \in B^2_1(0),        
\end{equation}  
of points $y$ in $B^2_1(0)$, with $\mu \in B^2_r(0)$ 
for some sufficiently small $r \in (0,\varepsilon_0)$ - see 
here Section 2 in \cite{Escher.Pruess.Simonett} - 
and the induced truncated shift about $p$ on $\Sigma$: 
\[ \Theta_{\mu}(x):=\left\{ \begin{array}{r@{\quad:\quad}l} 
	(\varphi_{\iota})^{-1}(\theta_{\mu}(\varphi_{\iota}(x))
	& x \in O_{\iota}                               \\ 
	x             &     x\not \in O_{\iota}. 
\end{array} \right.    \]
By Proposition 2.4 in \cite{Escher.Pruess.Simonett} 
- see also Lemma 3.1 in \cite{Shao.2015} -  
the truncated shift $\theta_{\mu}$ induces 
topological automorphisms of $h^{s}(U,\rel^n)$ 
for any open subset $U \subset \rel^2$ 
containing the closed ball $\overline{B_3}=\overline{B^2_{3\varepsilon_0}(0)}$ 
and for any H\"older exponent 
$s \in \rel_+\setminus \nat_0$, simply in 
terms of composition: 
$$ 
u \mapsto (\theta_{\mu})^{*}(u) = u \circ \theta_{\mu}
$$ 
with invers 
$$ 
u \mapsto (\theta^{\mu})_{*}(u)=u \circ (\theta_{\mu})^{-1}, 
$$
for functions $u:U \longrightarrow \rel^n$.
Similarly, for functions $u:\Sigma \longrightarrow \rel^n$ the truncated shift $\Theta_{\mu}$ 
induces a transformation 
$$
u \mapsto \Theta^*_{\mu}(u):=
\varphi_{\iota}^*\theta_{\mu}^*\psi^*_{\iota} 
(\varsigma u) + (1-\varsigma) \, u,  
$$
where we identified $\varsigma$ with its pullback
$\varphi_{\iota}^*\varsigma$ from $B_1^2(0)$ onto 
$O_{\iota}\subset \Sigma$, and this transformation 
is again an automorphism of 
$h^s(\Sigma,\rel^n)$ for any H\"older exponent 
$s \in \rel_+\setminus \nat_0$, on account of 
Proposition 3.3 in \cite{Shao.2015}.
Now, as in Section 3.4 in \cite{Shao.2015} we choose 
a time $T \in (0,T_{\textnormal{max}}(U_0))$, 
another time $t_0 \in (0,T)$, 
a small open interval $B_{3\varepsilon_0}(t_0)\subset (0,T)$ 
about $t_0$ and another bump function 
$\xi \in C^{\infty}_c((0,T),\rel)$ with 
$$
\textnormal{supp}(\xi)\subset   (t_0-2\varepsilon_0,t_0+2\varepsilon_0)
\quad \textnormal{and} \quad  \xi \equiv 1 \quad \textnormal{on} \quad [t_0-\varepsilon_0,t_0+\varepsilon_0],
$$
yielding the time-dependent time-translations 
\begin{equation}  \label{varrho}
\varrho_{\lambda}(t):= t + \lambda \,\xi(t), 
\quad \textnormal{for} \,\,\, t\in [0,T], 
\end{equation} 
for fixed parameters $\lambda \in \rel$.
Furthermore, following the procedure in Section 
3.4 of \cite{Shao.2015}, we consider the time-dependent 
shift $\theta_{\xi(t)\cdot\mu}(y) = y +\chi(y) \,\xi(t)\,\mu$, inducing 
$$
\tilde T_{\mu}(u)(y,t)
:=(\theta_{\xi(t)\cdot\mu})^*(u)(y,t), 
\quad \textnormal{for} \quad 
(y,t) \in U \times [0,T], 
$$ 
for functions $u:U \times [0,T] \longrightarrow 
\rel^n$, and analogously 
\begin{equation}   \label{T.mu} 
T_{\mu}(u)(x,t):=
(\Theta_{\xi(t)\cdot\mu})^*(u)(x,t), 
\quad \textnormal{for} \quad (x,t) \in 
\Sigma \times [0,T], 
\end{equation}  
for functions $u:\Sigma \times [0,T] \longrightarrow 
\rel^n$ and for $|\mu|<r$, with $r$ as in 
(\ref{theta.mu}). Finally, together with the time-shift in (\ref{varrho}) we consider the local modifications
	\begin{eqnarray}   \label{Phi.lambda.mu}
	u_{\lambda,\mu}(x,t):=(\Phi_{\lambda,\mu}^*u)(x,t)
	:= (T_{\mu} \circ \varrho_{\lambda}^*)(u)(x,t)
	\nonumber \\ 
	= u(x + \xi(t)\,\chi(x) \,\mu, t + \xi(t)\,\lambda),  
	\,\,  \textnormal{for} \,\, (x,t) \in  
	\Sigma \times [0,T],  
	\end{eqnarray}  
	of any given function $u :\Sigma \times [0,T] 
	\longrightarrow \rel^n$, with  
	$|\lambda|^2 + |\mu|^2<r^2$. By Proposition 3.10 
	in \cite{Shao.2015} and by the open mapping theorem  
	the local modifications in (\ref{Phi.lambda.mu}) 
	induce topological automorphisms of the Banach space 
	\begin{equation}  \label{E.1}
	{\bf E}_1:= C^0([0,T];h^{4+\beta}(\Sigma,\rel^n)) 
	\cap C^{1}([0,T];h^{\beta}(\Sigma,\rel^n))		
	\end{equation} 
	for some arbitrary $\beta \in (0,1)$, provided 
	$r$ in (\ref{theta.mu}) and (\ref{Phi.lambda.mu}) 
	is chosen sufficiently small. We note here for later 
	use, that there holds
	$$
	u_{\lambda,\mu}(x,0)= u(x+\xi(0)\,\chi(x) \,\mu, \xi(0)\,\lambda)=u(x,0), \quad 
	\textnormal{for} \quad x \in \Sigma,
	$$
	for any fixed function $u:\Sigma \times [0,T] \longrightarrow \rel^n$. Moreover, Proposition 3.10 in \cite{Shao.2015} yields a certain operator 
	\begin{equation}  \label{B.real.analytic}
		[(\lambda,\mu) \mapsto B_{(\lambda,\mu)}] \in 
		C^{\omega}(B_{r}^3(0),C^0([0,T],
		\Lift(h^{4+\beta}(\Sigma,\rel^n),
		h^{\beta}(\Sigma,\rel^n)))) 
	\end{equation}
	with $B_{(\lambda,0)} \equiv 0$ such that 
	\begin{equation}   \label{Prop.3.10}
	\partial_t u_{\lambda,\mu}(x,t) 
	= (1+\xi'(t)\,\lambda) \,(\Phi_{\lambda,\mu})^*(\partial_tu_t)(x)
	+ B_{(\lambda,\mu)}(u_{\lambda,\mu})(x,t)  
	\end{equation}
	holds $\forall \,(x,t) \in \Sigma \times [0,T]$ and for any fixed family of functions $\{u_t\}\in {\bf E}_1$.
	Now we apply formulae (\ref{Phi.lambda.mu}) and (\ref{Prop.3.10}) directly to the flow line $\{\PP^*(t,0,U_0)\}_{t\in [0,T]}$ - which we 
	shall henceforth abbreviate by $\{u^*_t\}_{t\in [0,T]}$ - 
	and we shall use the important fact, that 
	the right hand side of evolution equation (\ref{de_Turck_equation}) is given by the time-independent, quasilinear operator in formula (\ref{M.F_0.2}) respectively (\ref{A.F_0}) along smooth flow lines, 
	on account of the first part of Theorem \ref{Psi.of.class.C_1}. We see therefore, 
	that the function $u^*_{\lambda,\mu} \equiv 
	(\Phi_{\lambda,\mu})^*(u^*)$ satisfies
	$u^*_{\lambda,\mu}(\,\cdot \,,0)
	=u^*(\,\cdot \,,0)=U_0$ on $\Sigma$ and the equation
	\begin{eqnarray}  \label{prepaing.Psi} 
		\partial_t u^*_{\lambda,\mu}(x,t) 
		= (1+\xi'(t)\,\lambda) \,(\Phi_{\lambda,\mu})^*(\partial_tu^*_t)(x)
		+B_{(\lambda,\mu)}(u^*_{\lambda,\mu})(x,t)  \qquad  \\
		= (1+\xi'(t)\,\lambda) \,
		T_{\mu} \circ
		\varrho_{\lambda}^*(K^{F_0}(u^*_t).(u^*_t))(x)
		+B_{(\lambda,\mu)}(u^*_{\lambda,\mu})(x,t)    \nonumber                                \\
		= (1+\xi'(t)\,\lambda) \,
		T_{\mu} \big{(} 
		K^{F_0}(\varrho_{\lambda}^* u^*_t).(\varrho_{\lambda}^*u^*_t)\big{)}(x)
		+B_{(\lambda,\mu)}(u^*_{\lambda,\mu})(x,t)          \nonumber   \\
		=(1+\xi'(t)\,\lambda) \,
		T_{\mu} \big{(}K^{F_0}((T_{\mu})^{-1}(u^*_{\lambda,\mu})).
		((T_{\mu})^{-1}(u^*_{\lambda,\mu}))\big{)}(x,t) 
		+ B_{(\lambda,\mu)}(u^*_{\lambda,\mu})(x,t)  \nonumber   \\
		=(1+\xi'(t)\,\lambda) \,
		T_{\mu} \circ \Mill_{F_0} \circ (T_{\mu})^{-1}(u^*_{\lambda,\mu})(x,t) 
		+ B_{(\lambda,\mu)}(u^*_{\lambda,\mu})(x,t)  
		\nonumber  
	\end{eqnarray} 
	for $(x,t) \in \Sigma \times [0,T]$, where we have 
	used the time-independence of the operator 
	$K^{F_0}$ from formula (\ref{A.F_0}) in the third line. Following now the notation of Section 4 in \cite{Shao.2013} or Sections 4--6 in \cite{Shao.2015}, we shall abbreviate   
	$E_0:=h^{\beta}(\Sigma,\rel^n)$ and
	$E_1:=h^{4+\beta}(\Sigma,\rel^n)$, and moreover 
    $$ 
    {\bf E}_0:=C^0([0,T];h^{\beta}(\Sigma,\rel^n)),  
	$$ 
	for the same $\beta \in (0,1)$ as in the definition
	of the space ${\bf E}_1$ in (\ref{E.1}). 
	Hence, in this proof we shall exchange the parabolic  
	$L^p$-spaces $L^p([0,T];L^p(\Sigma,\rel^n))$ 
	and $X_{T,p}$ from (\ref{X.T}) by the smaller function 
	spaces ${\bf E}_0$ and ${\bf E}_1$, in order to 
	argue correctly throughout this proof. 
	Since the considered flow line $\{u^*_t\}_{t\in [0,T]}$ 
	of equation (\ref{de_Turck_equation}) is smooth and thus 
	contained in the Banach space ${\bf E}_1$, 
	we can choose a small open neighborhood $V_{U_0,T,p}$ about $\{u^*_t\}_{t\in [0,T]}$ in ${\bf E}_1$, 
	and we shall introduce the non-linear operator 
	$\Psi: V_{U_0,T,p} \times B_{r}^3(0) 
	\longrightarrow E_1 \times {\bf E}_0$, sending 
	the triple $(\{f_t\},(\lambda,\mu))$ to the pair 
	\[\left(  \begin{array}{c} \gamma_0(\{f_t\}) - U_0    \\
	\partial_t f_t - (1+\xi'(t)\,\lambda) \, 
	T_{\mu}\circ \Mill_{F_0} \circ (T_{\mu})^{-1}(\{f_{t}\}) 
	- B_{(\lambda,\mu)}(\{f_{t}\})
	\end{array}  \right).                     \]
    We should note here for later use, that statement (\ref{prepaing.Psi}) is equivalent to: 
    \begin{equation}  \label{Psi.lambda,mu} 
    \Psi(u^*_{\lambda,\mu},(\lambda,\mu))=0   \quad 
    \textnormal{in} \quad E_1 \times {\bf E}_0,
    \end{equation}
    for every pair $(\lambda,\mu) \in B_r^3(0)$.    
    Now we shall prove exactly as in Section 4 of  
    \cite{Shao.2013}, that this operator is 
    real analytic, provided $V_{U_0,T,p}$ and $r>0$ 
    are chosen sufficiently small. We can firstly infer from   
    (\ref{B.real.analytic}) and Proposition 2.5 in 
    \cite{Shao.2015} - exactly as on p. 12 in 
    \cite{Shao.2013} - that there holds: 
\begin{equation}  \label{B.mu.real.analytic}
[(\{f_t\},(\lambda,\mu)) \mapsto B_{(\lambda,\mu)}(\{f_t\})] \in C^{\omega}(V_{U_0,T,p} \times B_{r}^3(0),{\bf E}_0). 
\end{equation}
Now we aim to prove, that the technically most 
challenging operator  
\begin{equation}   \label{most.complicated} 
[(\{f_t\},\mu) \mapsto 
(T_{\mu}\circ \Mill_{F_0} \circ (T_{\mu})^{-1})(\{f_{t}\})]
\end{equation}  
within the above definition of $\Psi$ is of class 
$C^{\omega}(V_{U_0,T,p} \times B_{r}^2(0),{\bf E}_0)$. 
To this end, we proceed exactly as on pp. 12--13
in \cite{Shao.2013}, making again use of the fact, that the non-linear differential operator $[f \mapsto \Mill_{F_0}(f)]$ 
has actually a quasilinear structure
$[f \mapsto K^{F_0}(f).(f)]$, which is given by formula (\ref{A.F_0}) respectively (\ref{phi_revisited}), i.e. 
that there holds: 
	\begin{eqnarray}  \label{quasilinear.OP} 
	\Mill_{F_0}(f_t)(x)   
	= (K^{F_0}(f_t).(f_t))(x)      \nonumber \\              
	= -\frac{1}{2}\,\frac{1}{|A^0_{f_t}|^4(x)} \,
	g^{ij}_{f_t} \, g^{kl}_{f_t} 
	\, \nabla_i^{F_0} \nabla_j^{F_0} \nabla_k^{F_0} \nabla_l^{F_0}(f_t)(x)     	                  \\
	+ \tilde \NN^{F_0}(x,D_xf_t(x),D_x^2f_t(x)) 
	\cdot D^3_xf_t(x)                        \nonumber  \\
	+ \tilde \CC^{F_0}(x,D_xf_t(x),D_x^2f_t(x)) 
	\cdot D_x^2f_t(x)                     \nonumber  \\
	+ \tilde \dom^{F_0}(x,D_xf_t(x),D_x^2f_t(x)) 
	\cdot D_xf_t(x)                    \nonumber
\end{eqnarray}  
for $(x,t)\in \Sigma \times [0,T]$ and every 
$\{f_t\} \in V_{U_0,T,p}$ - provided $V_{U_0,T,p}$ 
is chosen sufficiently small in ${\bf E}_1$ - 
where the functions $\tilde \NN^{F_0}$, 
$\tilde \CC^{F_0}$ and 
$\tilde \dom^{F_0}$ in (\ref{quasilinear.OP}) 
have the same algebraic structure as 
the functions $\NN^{F_0}$, 
$\CC^{F_0}$ and $\dom^{F_0}$
in formulae (\ref{A.F_0}) and (\ref{phi_revisited}). 
In order to correctly deal with the 
quasilinearity of the decisive factor 
$[f \mapsto g^{ij}_{f} \, g^{kl}_{f} \, 
\nabla_i^{F_0} \nabla_j^{F_0} \nabla_k^{F_0} \nabla_l^{F_0}(f)]$ of the leading differential 
operator in (\ref{quasilinear.OP}), 
we have to localize this operator on 
the particular coordinate patch $O_{j_p}$ 
containing the point $p$, which we had chosen above. 
As explained on p. 68 in \cite{Shao.2015}, we 
can use a fixed localization system 
$\{(\pi_j,\zeta_j)\}_{j=1,\ldots,N}$ subordinate 
to our original conformal atlas 
$\Area=\{(O_j,\varphi_j)\}_{j=1,\ldots,N}$, 
meeting properties (L1)--(L3) on p. 49 in \cite{Shao.2015}, 
together with the additional localization function 
$\zeta_{\iota}:=\varphi^*_{\iota}(\zeta)$ 
on the auxiliary coordinate patch $O_{\iota}$ 
- similarly to the original bump functions 
$\zeta_j:= \varphi^*_{j}(\zeta)$, 
which are pullbacks w.r.t. the original charts 
$\varphi_{j}$ of a fixed bump function 
$\zeta \in C^{\infty}_c(B^2_1(0),[0,1])$  
appearing in condition (L2) on p. 49 in \cite{Shao.2015} -     in order to define a new, more suitable partition of unity 
subordinate to our original atlas $\Area$ by: 
$$ 
\tilde \pi_{j}:= \pi_{j}^2 - \tilde \pi_{j_p} \, \pi_{j}^2, 
\quad 1 \leq j \leq N, \quad \textnormal{with} \quad   
\tilde \pi_{j_p}:=\zeta_{j_p} + \zeta_{\iota} 
- \zeta_{j_p} \, \zeta_{\iota}. 
$$ 
The particular bump function $\pi:=\sqrt{\tilde \pi_{j_p}}$ 
has support in the original chart $(O_{j_p},\varphi_{j_p}) 
\in \Area$ containing the chosen point $p$, and its construction guarantees that 
$\psi_{\iota}(B_4) \subset [\pi =1]$, provided 
$\varepsilon_0>0$ is chosen appropriately.   
Now first of all, we obtain the operator identity
\begin{equation}   \label{split.identity}  
K^{F_0}(h) = \pi^2 \cdot K^{F_0}(h)  
+ \sum_{j \not = j_p} \tilde \pi_j \cdot K^{F_0}(h)  
\quad \textnormal{in} \,\,\, \Lift(E_1,E_0),
\end{equation}  
for every fixed 
$h \in\textnormal{Imm}_{\textnormal{uf}}(\Sigma,\rel^n)$. 
According to formula (\ref{derivation_formulae}) 
we have for any immersion $h \in E_1$:     
\begin{eqnarray}  \label{new.fourth.order} 
\pi(x) \cdot \Big{(} g^{ij}_{h}(x) \, g^{kl}_{h}(x) \, 
\nabla_i^{F_0} \nabla_j^{F_0} \nabla_k^{F_0} 
\nabla_l^{F_0}(h)(x) \Big{)} 
\nonumber \\
=\pi(x) \cdot \Big{(} g^{ij}_{h}(x) \, g^{kl}_{h}(x) \,\partial_{ijkl}h(x) 
+D^{ijk}_3(x,\partial_xh(x),\partial^2_xh(x)) \cdot \partial_{ijk}h(x)            \nonumber \\
+\tilde D(x,\partial_xh(x),\partial^2_{x}h(x)) \Big{)},\quad 
\end{eqnarray} 
for $x\in O_{j_p}$, where 
$\partial_xh(x):=(\partial_1h,\partial_2h)(x)$ 
and $\partial^2_xh(x):=
(\partial^2_{11}h,\ldots,\partial^2_{22}h)(x)$, 
and where $D^{ijk}_3$ and $\tilde D$ are 
$\textnormal{Mat}_{n,n}(\rel)$- respectively 
$\rel^n$-valued functions on $O_{j_p}$, 
whose components are rational functions of 
the partial dervatives of $h$ up to 
second order - as precisely indicated in (\ref{new.fourth.order}) - and additionally of the partial derivatives of $F_0$ 
up to fourth order, which are here real analytic. 
Using the fixed chart $\psi_{j_p}=(\varphi_{j_p})^{-1}:
B^2_1(0) \stackrel{\cong}\longrightarrow O_{j_p}$, 
the partial derivatives in (\ref{new.fourth.order}) 
correspond to partial derivatives on $B^2_1(0)$, 
and therefore the differential operator on the 
right-hand side of (\ref{new.fourth.order}) 
is a well-defined quasilinear differential operator 
on $O_{j_p}$, which extends smoothly to the 
zero-operator on entire $\Sigma$. Similarly to 
p. 67 in \cite{Shao.2015} and to p. 13 in \cite{Shao.2013}
we infer from formulae (\ref{derivation_formulae}) and 
(\ref{Cramers.rule}), that the right-hand side 
in (\ref{new.fourth.order}) can be written as a quotient,
whose nominator is a sum of finitely many products   
of locally defined, linear differential operators 
up to order $o=4$, i.e. of
\begin{equation}   \label{nominator} 
\Big{[}v \mapsto c^o_{\alpha} \, \partial^{\alpha}(v) \Big{]},
\quad \textnormal{for functions}  \,\,\, 
v \in C^{4}(O_{j_p},\rel^n), 
\end{equation}
and whose denominator is exactly given by: 
\begin{equation}   \label{denominator} 
	\Big{[}v \mapsto 
	\big{(}\det \big{[}(\langle \partial_i(v),\partial_j(v)\rangle)_{ij}\big{]}\big{)}^2 
	\Big{]}, \quad \textnormal{for functions}  \,\,\, 
	v \in C^{1}(O_{j_p},\rel^n), 
\end{equation}
here applied to the restriction of the immersion $h \in E_1$
to the coordinate patch $O_{j_p} \subset \Sigma$,   
where $\alpha=(\alpha_1,\alpha_2) \in (\nat_0)^2$ 
with $1 \leq o:=|\alpha|=\alpha_1+\alpha_2 \leq 4$,
and where especially $c^4_{\alpha}(x) = \pi(x)$ for 
$|\alpha|=4$. In view of the proof of Proposition 3.7 in \cite{Shao.2015}, we shall interpret here the operators in 
(\ref{nominator}) as differential operators mapping $h^{4+\beta}_{\textnormal{cp}}(\Sigma,\rel^n)$ 
into $h^{\beta}_{\textnormal{cp}}(\Sigma,\rel^n)$,  
adopting here the non-standard notation 
$$ 
h^{s}_{\textnormal{cp}}(\Sigma,\rel^n):= 
\{\, u \in h^{s}(\Sigma,\rel^n) \, | \, 
\textnormal{supp}(u) \subseteq \psi_{\iota}(\overline{B_5})\,\}, \,\,\, \textnormal{for} \,\,\, s \in \rel_+\setminus \nat_0,
$$
from Section 3.1 in \cite{Shao.2015},
and we should note here, that the above coefficients 
$c^o_{\alpha}$, with $o=|\alpha| \in \{1,2,3\}$, 
are smooth real-valued functions on the entire patch $O_{j_p}$ 
and even real analytic functions on any open subset 
$U_{j_p} \subset O_{j_p}$ satisfying 
$U_{j_p} \subset [\pi=1]$, because they are composed 
of certain partial derivatives of the components of 
the real analytic immersion $F_0$, 
and they might also contain the smooth bump function 
$\pi$ from (\ref{new.fourth.order}) as a factor. 
As mentioned at the beginning of Section 3.2 in 
\cite{Shao.2015}, we should choose here the open 
set $U_{j_p}$ in such a way, that  
$\psi_{\iota}(\overline{B_3}) \subset U_{j_p} \subset 
\psi_{\iota}(B_4)$ holds, in view of the fact that 
$\psi_{\iota}(B_4)$ is contained in $[\pi=1]$, by construction. 
Hence, exactly as in the proof of Proposition 3.7 in \cite{Shao.2015} we can conclude via Theorem 4.2 in \cite{Escher.Pruess.Simonett} combined with Lemma 3.1 
in \cite{Shao.2015}, that the maps  
\begin{equation}  \label{real.analytic.partial.deriv}
\Big{[}\mu \mapsto 
\Theta_{\mu} \circ \,
\Big{(} c^o_{\alpha} \,\partial^{\alpha} \Big{)} \circ (\Theta_{\mu})^{-1} \Big{]}
\end{equation}  
are of class $C^{\omega}(B_{r}^2(0),
\Lift(h^{4+\beta}_{\textnormal{cp}}(\Sigma,\rel^n),
h^{\beta}_{\textnormal{cp}}(\Sigma,\rel^n)))$, 
for each order $o=|\alpha| \in \{1,2,3,4\}$, 
provided $0<r \ll \varepsilon_0$. 
Now we may apply Lemma 3.8 in \cite{Shao.2015} respectively 
Lemma 5.1 in \cite{Escher.Pruess.Simonett}
to statement (\ref{real.analytic.partial.deriv}) - 
here with the basic Banach space $\Lift(h^{4+\beta}_{\textnormal{cp}}(\Sigma,\rel^n),
h^{\beta}_{\textnormal{cp}}(\Sigma,\rel^n))$ -  
and we infer, that the maps
\begin{equation}  \label{real.analytic.partial.deriv.2}
\Big{[}\mu \mapsto 
T_{\mu} \circ \,\Big{(} 
c^{o}_{\alpha}\, \partial^{\alpha} \Big{)} \circ (T_{\mu})^{-1}
\Big{]}
\end{equation}  
are of class $C^{\omega}(B_{r}^2(0),C^0([0,T],
\Lift(h^{4+\beta}_{\textnormal{cp}}(\Sigma,\rel^n),
h^{\beta}_{\textnormal{cp}}(\Sigma,\rel^n))))$, for 
each order $o=|\alpha| \in \{1,2,3,4\}$, provided $r>0$ 
is sufficiently small. 
See here the definition in (\ref{T.mu}) and also formula (4.8) in \cite{Shao.2015}. Since the ``cut-off'' differential operator in (\ref{new.fourth.order}) extends trivially to 
entire $\Sigma$ and is a rational expression of linear 
differential operators as pointed out in 
(\ref{nominator}) and (\ref{denominator}), we can finally combine the result in (\ref{real.analytic.partial.deriv.2}) 
via Proposition 6.4 in \cite{Shao.Simonett} and
Proposition 2.5 in \cite{Shao.2015}, in order 
to conclude that the map 
\begin{equation}  \label{real.analytic.partial.deriv.3}
\Big{[}(\{f_t\},\mu)  \mapsto
T_{\mu} \circ \, \Big{(}
\pi \cdot \big{(} g^{ij}_{T_{\mu}^{-1}(\{f_t\})}\, g^{kl}_{T_{\mu}^{-1}(\{f_t\})} \, 
\nabla_i^{F_0} \nabla_j^{F_0} \nabla_k^{F_0} 
\nabla_l^{F_0}(T_{\mu}^{-1}(\{f_t\})) \big{)}\Big{)} 
\Big{]} 
\end{equation} 
is of class $C^{\omega}(V_{U_0,T,p} \times 
B_{r}^2(0),{\bf E}_0)$, where we should keep 
in mind the fact, that the derivatives of fourth and third 
order in (\ref{new.fourth.order}) contribute only 
affine linearly to (\ref{real.analytic.partial.deriv.3}).    
On account of formulae 
(\ref{double_nabla})--(\ref{mean_curvat_laplacian}) 
and (\ref{A0.squared}) exactly the same method can be employed, in order to prove that the remaining factor 
$\frac{1}{2|A^0_{f_t}|^4}$ in (\ref{quasilinear.OP}) 
gives rise to a map    
\begin{equation}  \label{real.analytic.leading.coeff}
\Big{[}(\{f_t\},\mu) \mapsto 
T_{\mu} \circ \Big{(} \pi \cdot  
\frac{1}{2\,|A^0_{(T_{\mu})^{-1}(\{f_t\})}|^4}  \,\Big{)} 
\Big{]}
\end{equation}  
of class $C^{\omega}(V_{U_0,T,p} \times 
B_{r}^2(0),C^0([0,T];h^{\beta}(\Sigma,\rel)))$. 
Now, combining the results in (\ref{real.analytic.partial.deriv.3}) 
and (\ref{real.analytic.leading.coeff}) again via 
Proposition 2.5 in \cite{Shao.2015}, we infer that the 
entire leading operator on the right 
hand side of (\ref{quasilinear.OP}) yields a map 
\begin{eqnarray}  \label{real.analytic.leading.operator}
	\Big{[}(\{f_t\},\mu) \mapsto 
	T_{\mu} \circ \Big{(}
	\pi^2 \cdot \frac{1}{2}\,\frac{1}{|A^0_{(T_{\mu})^{-1}(\{f_t\})}|^4} \,g^{ij}_{(T_{\mu})^{-1}(\{f_t\})} 
	\, g^{kl}_{(T_{\mu})^{-1}(\{f_t\})}     \nonumber    \\ 
	\cdot \nabla_i^{F_0} \nabla_j^{F_0} \nabla_k^{F_0} \nabla_l^{F_0}((T_{\mu})^{-1}(\{f_t\})) \Big{)} \,\Big{]} \quad
\end{eqnarray}  
of regularity class $C^{\omega}(V_{U_0,T,p} \times B_{r}^2(0),{\bf E}_0)$.
Since the remaining three terms in (\ref{quasilinear.OP})
of order $<4$ can be handled in exactly the same way, 
the above method shows that the map 
$$
\Big{[}(\{f_t\},\mu) \mapsto 
\Big{(} T_{\mu}\circ \Big{(} \pi^2 \cdot \Mill_{F_0}  \Big{)} \circ (T_{\mu})^{-1}\Big{)}(\{f_{t}\})\Big{]}
$$
is of regularity class $C^{\omega}(V_{U_0,T,p} 
\times B_{r}^2(0),{\bf E}_0)$.  
Since each remaining map  
$[(\{f_t\},\mu) \mapsto 
(T_{\mu}\circ (\tilde \pi_j \cdot \Mill_{F_0}) \circ (T_{\mu})^{-1})(\{f_{t}\})]$, for $j \not = j_p$, 
according to our decomposition in (\ref{split.identity}), 
is trivially of regularity class $C^{\omega}(V_{U_0,T,p} 
\times B_{r}^2(0),{\bf E}_0)$ as well, we can 
sum them up and finally infer from  
$\sum_{j=1}^N \tilde \pi_j = 1$ on $\Sigma$ the    
correctness of our claim in (\ref{most.complicated}).
Together with statement (\ref{B.mu.real.analytic}) 
we conclude again via Proposition 2.5 in 
\cite{Shao.2015}, that the map $\Psi$ is of class 
$C^{\omega}(V_{U_0,T,p} \times B_{r}^3(0),{\bf E}_0)$
indeed, provided $V_{U_0,T,p}$ and $r>0$ are chosen 
sufficiently small.   \\
Furthermore, one can prove exactly as in the second and third part of Theorem \ref{Psi.of.class.C_1}, that the
quasilinear differential operator 
$$
\partial_t-\Mill_{F_0}:V_{U_0,T,p} \subset 
{\bf E}_1 \longrightarrow {\bf E}_0
$$ 
is real analytic - here w.r.t. the norms of the 
Banach spaces ${\bf E}_1$ and ${\bf E}_0$, 
differing from the corresponding $L^p$-spaces - and 
that its Fr\'echet-derivative 
$$
\partial_t - D\Mill_{F_0}(\{f_t\}): 
T_{\{f_t\}}V_{U_0,T,p} = {\bf E}_1 
\longrightarrow {\bf E}_0
$$ 
in any fixed $\{f_t\} \in V_{U_0,T,p}$ 
is a linear, uniformly parabolic operator of fourth 
order, concretely given by formula (\ref{Frechet.of.Psi}) 
again. Moreover, since the considered flow line 
$\{u^*_t\}_{t\in [0,T]} \equiv \{u(t,U_0)\}_{t\in [0,T]}$ 
is $C^{\infty}$-smooth, one can argue as in the proof 
of Proposition \ref{A.priori.estimate} -  
but now via Theorems 3.3 and 3.7 in 
\cite{Shao.Simonett}, here applied simply with interpolation parameter $\gamma = 1$ - that for each fixed time 
$s \in [0,T]$ the linear operator 
$$
\Big{(} \gamma_0,\partial_t - 
D\Mill_{F_0}(u^*_s) \Big{)}: 
{\bf E}_1 \stackrel{\cong} \longrightarrow 
E_1 \times {\bf E}_0
$$ 
is a topological isomorphism. Now, again on account of  
formula (\ref{Frechet.of.Psi}) combined with the 
smoothness of $\{u^*_t\}_{t\in [0,T]}$, we can infer  
from Lemma 2.8 (a) in \cite{Clement.Simonett} - 
here again with interpolation parameter $\mu=1$ - that 
\begin{equation}  \label{max.reg}
\Big{(}\gamma_0, \partial_t 
- D\Mill_{F_0}(\{u^*_t\})\Big{)}: 
{\bf E}_1 \stackrel{\cong} \longrightarrow 
E_1 \times {\bf E}_0
\end{equation} 
is an isomorphism as well, similarly to the statement 
of the fourth part of Theorem \ref{Psi.of.class.C_1}. 
On the other hand, the Fr\'echet-derivative of the 
real analytic map $\Psi:V_{U_0,T,p} \times B^3_r(0) \longrightarrow E_1 \times {\bf E}_0$ w.r.t. its 
first argument in the special point 
$(\{f_t\},(\mu,\lambda))=(\{u^*_t\},(0,0))$   
is exactly the linear, isomorphic operator in 
(\ref{max.reg}), due to $B_{(0,0)} \equiv 0$. 
Since we also know that 
$\Psi(u^*_{\lambda,\mu},(\lambda,\mu)) =0$ 
for all $(\lambda,\mu) \in B_r^3(0)$ by (\ref{Psi.lambda,mu}), 
the Implicit Function Theorem for real analytic
operators [\cite{Zeidler}, Theorem 4.B and 
Corollary 4.23] yields some small open, non-empty ball 
$B^3_{\epsilon_0}(0) \subset 
B^3_{r}(0)$ and a unique real analytic map 
$h:B^3_{\epsilon_0}(0) \longrightarrow V_{U_0,T,p}$, 
such that $h(\lambda,\mu)=u^*_{\lambda,\mu}$ for 
$(\lambda,\mu) \in B^3_{\epsilon_0}(0)$,   
implying especially that 
$(\lambda,\mu) \mapsto u^*_{\lambda,\mu}$ 
is a real analytic function from $B^3_{\epsilon_0}(0)$ 
into ${\bf E}_1$. Hence, Theorem 1.1 
in \cite{Shao.2015} respectively Theorem 4.1 in 
\cite{Shao.Simonett} finally yield the asserted real analyticity of the flow line $\{u^*_t\}_{t\in [0,T]}$ about both any point 
$p \in \Sigma$ and any time $t_0 \in (0,T)$, 
and for every $T \in (0,T_{\textnormal{max}}(U_0))$.   
\qed               
\end{itemize}

\section{Main theorems \ref{Frechensbergo} and \ref{Functional.prop.Frechet.derivative} and their proofs} \label{Linearization} 

In this section we will prove the main results of this article. We shall start with a further investigation of the evolution operator of the ``DeTurck modification'' (\ref{de_Turck_equation}) of the MIWF (\ref{Moebius.flow}) and of its linearization. For any fixed immersion $F_0 \in \textnormal{Imm}_{\textnormal{uf}}(\Sigma,\rel^n)
\cap C^{\infty}(\Sigma,\rel^n)$ and for 
any fixed $p \in (3,\infty)$ we will denote  
by $B_{\rho}^{4-\frac{4}{p},p}(F_0)$ the open ball 
of radius $\rho>0$ about $F_0$ in the space of traces $W^{4-\frac{4}{p},p}(\Sigma,\rel^n)$ from line 
(\ref{Trace.X.T}), and we will abbreviate 
$X_T$ for the basic Banach space $X_{T,p}$ from line 
(\ref{X.T}), and similarly $Y_T:=Y_{T,p}$. 
\begin{theorem}   \label{Frechensbergo}
	Suppose that $\Sigma$ is a smooth compact torus 
	and that $F_0$ is a smooth and umbilic-free immersion, 
	i.e. of class $C^{\infty}(\Sigma,\rel^n) 
	\cap \textnormal{Imm}_{\textnormal{uf}}(\Sigma,\rel^n)$. 
	Furthermore, let $T\in (0,T_{\textnormal{max}}(F_0))$ and 
	$p\in (3,\infty)$ be arbitrarily chosen, where $T_{\textnormal{max}}(F_0)$ has been defined in Definition 
	\ref{Umbilic.free}. Then the following statements hold:  
	\begin{itemize}
	\item[1)] There is some $\rho=\rho(\Sigma,F_0,T,p)>0$, such that for every $F \in W^{4-\frac{4}{p},p}(\Sigma,\rel^n)$ with $\parallel F-F_0 \parallel_{W^{4-\frac{4}{p},p}(\Sigma,\rel^n)}<\rho$ the existence interval $[0,t^+(F))$ of the 
	maximal solution $\{\PP^*(\,\cdot \,,0,F)\}$
	of Cauchy problem (\ref{initial.value.problem}) from the second part of Theorem \ref{short.time.solution} 
	contains the prescribed interval $[0,T]$, and the resulting evolution operator 
	\begin{equation}  \label{solution.operator.1}
	\PP^*(\,\cdot\,,0,\,\cdot\,):
	B_{\rho}^{4-\frac{4}{p},p}(F_0) \subset W^{4-\frac{4}{p},p}(\Sigma,\rel^n) 
	\longrightarrow X_{T}, 
	\end{equation}
	mapping any element $F$ of the open ball $B_{\rho}^{4-\frac{4}{p},p}(F_0)$ about $F_0$ in $W^{4-\frac{4}{p},p}(\Sigma,\rel^n)$ to the 
	unique solution $\{\PP^*(t,0,F)\}_{t \in [0,T]}$ 
	of Cauchy problem (\ref{initial.value.problem}) in 
	$X_{T}$, is of class $C^{\omega}$. 
	\item[2)] For any fixed $\xi \in W^{4,p}(\Sigma,\rel^n)$ the function $[(x,t) \mapsto (D_{F}\PP^*(t,0,F_0).\xi)(x)]$ arising from the Fr\'echet derivative 
	\begin{eqnarray}  \label{Frechet.derivative}
	D_{F}\PP^*(\,\cdot\,,0,F_0):
	W^{4-\frac{4}{p},p}(\Sigma,\rel^n) 
	\longrightarrow X_{T} 
	\end{eqnarray}
		of the evolution operator $[F\mapsto \PP^*(\,\cdot\,,0,F)]$ in line (\ref{solution.operator.1}) in the immersion $F_0$ is the \underline{unique solution} 
		$\{z_t\}$ of the linear uniformly parabolic initial value problem (\ref{Poincare.trick}) in 
		$C^{1}([0,T];L^{p}(\Sigma,\rel^n)) \cap C^{0}([0,T];W^{4,p}(\Sigma,\rel^n))$, 
		which starts moving at time $t=0$ in the given function $\xi \in W^{4,p}(\Sigma,\rel^n)$: 
		\footnote{According to Definition 4.1.1 in \cite{Lunardi} this means exactly, that the function 
		$[(x,t) \mapsto (D_{F}\PP^*(t,0,F_0).\xi)(x)]$ is the \underline{unique strict solution} of initial value problem (\ref{Poincare.trick}) w.r.t. the Banach spaces $X:= L^{p}(\Sigma,\rel^n)$ and 
		$D:= W^{4,p}(\Sigma,\rel^n)$.}
		\begin{eqnarray}  \label{Poincare.trick}
			\partial_t u_t(x)=  (D(\Mill_{F_0})(\PP^*(t,0,F_0)).u_t)(x)           \qquad \quad      \\
			\equiv 	- \frac{1}{2} \mid A^0_{\PP^*(t,0,F_0)}(x) \mid^{-4} \,
			g^{ij}_{\PP^*(t,0,F_0)} \, g^{kl}_{\PP^*(t,0,F_0)} \, \nabla_i^{F_0} \nabla_j^{F_0} \nabla_k^{F_0} \nabla_l^{F_0}(u_t)(x)  
			\nonumber \\ 
			-B_{3}^{ijk}(x,D_x\PP^*(t,0,F_0),
			D_x^2\PP^*(t,0,F_0)) \cdot \nabla^{F_0}_{ijk}(u_t)(x) 
			\nonumber \\
			-B_2^{ij}(x,D_x\PP^*(t,0,F_0),D_x^2\PP^*(t,0,F_0),D_x^3\PP^*(t,0,F_0),
			D_x^4\PP^*(t,0,F_0)) \cdot \nabla^{F_0}_{ij}(u_t)(x)       \nonumber \\
			-B_1^i(x,D_x\PP^*(t,0,F_0),
			D_x^2\PP^*(t,0,F_0),D_x^3\PP^*(t,0,F_0),
			D_x^4\PP^*(t,0,F_0)) \cdot \nabla^{F_0}_{i}(u_t)(x)       \nonumber\\
			\textnormal{for} \quad (x,t) \in \Sigma \times [0,T] \qquad
			\textnormal{and with} \quad  u_{0} = \xi  \quad \textnormal{on} \quad  \Sigma,
			\nonumber
		\end{eqnarray} 
		where 
		\begin{eqnarray*} 
			B_{3}^{ijk}(\,\cdot \,,D_x\PP^*(t,0,F_0),
			D_x^2\PP^*(t,0,F_0)),            \\
			B_2^{ij}(\,\cdot\,,D_x\PP^*(t,0,F_0),
			D_x^2\PP^*(t,0,F_0), 
			D_x^3\PP^*(t,0,F_0),
			D_x^4\PP^*(t,0,F_0)),         \\
			B_1^i(\,\cdot\,,D_x\PP^*(t,0,F_0),
			D_x^2\PP^*(t,0,F_0),	D_x^3\PP^*(t,0,F_0),D_x^4\PP^*(t,0,F_0))
		\end{eqnarray*} 
		are coefficients of $\textnormal{Mat}_{3,3}(\rel)$-valued, contravariant tensor fields of degrees $3,2$ and $1$, whose basic properties have been proved in Theorem \ref{Psi.of.class.C_1} (iii).
		\item[3)] The Fr\'echet derivative in line (\ref{Frechet.derivative}) of the evolution operator 
		$[F\mapsto \PP^*(\,\cdot\,,0,F)]$ in the umbilic-free immersion $F_0$ has the particular property, to yield for every fixed time 
		$t \in [0,T]$ a linear continuous operator 
		\begin{equation}  \label{Frechens.W.4.p}
			D_{F}\PP^*(t,0,F_0): W^{4,p}(\Sigma,\rel^n)  \longrightarrow  W^{4,p}(\Sigma,\rel^n),
		\end{equation}
		and this linear operator has a unique continuous linear extension 
		\begin{equation}  \label{Frechens.L.p}
			D_F\PP^*(t,0,F_0):L^{p}(\Sigma,\rel^n) \longrightarrow L^{p}(\Sigma,\rel^n)
		\end{equation} 
		for every fixed time $t \in [0,T]$, whose range is contained in $W^{4,p}(\Sigma,\rel^n)$, if $t\in (0,T]$.
		\item[4)] The differential equation (\ref{Poincare.trick})
		generates a parabolic fundamental solution 
		$$
		G^{F_0}(t,s):L^p(\Sigma,\rel^n) \longrightarrow  L^p(\Sigma,\rel^n),  \quad \forall \, s \leq t \in [0,T],
		$$
		with range in $W^{4,p}(\Sigma,\rel^n)$ for $t>s$, 
		which extends the family of linear operators 
		$[(t,s) \mapsto D_F\PP^*(t,s,\PP^*(s,0,F_0))]$ continuously and linearly from the space $W^{4,p}(\Sigma,\rel^n)$ to the space $L^p(\Sigma,\rel^n)$, i.e. we have 
		\begin{equation}  \label{Frechet.evol.operator.s}
			G^{F_0}(t,s) = D_F\PP^*(t,s,\PP^*(s,0,F_0)): 
			W^{4,p}(\Sigma,\rel^n)  \longrightarrow  W^{4,p}(\Sigma,\rel^n),
		\end{equation} 
		for every pair $s \leq t \in [0,T]$. In particular, this uniquely extended family of linear operators in 
		$\Lift(L^{p}(\Sigma,\rel^n),
		L^{p}(\Sigma,\rel^n))$ - again denoted by 
		$[(t,s) \mapsto D_F\PP^*(t,s,\PP^*(s,0,F_0))]$ - 
		has the semigroup property:
		\begin{equation}  \label{semigroup.property.b}   
			D_F\PP^*(t,s,\PP^*(s,0,F_0)) \circ D_F\PP^*(s,r,\PP^*(r,0,F_0))
			=D_F\PP^*(t,r,\PP^*(r,0,F_0))
		\end{equation} 
		in $L^{p}(\Sigma,\rel^n)$, for any $0 \leq r \leq s \leq t \leq T$. 
		\item[5)] For every fixed $t \in (0,T]$ and 
		$\xi^* \in (L^{p}(\Sigma,\rel^n))^*$ the function 
		$$
		y_s:= \Big{(}D_{F}\PP^*(t,s,\PP^*(s,0,F_0)) \Big{)}^*(\xi^*), \quad \textnormal{for} \,\,\, s \in [0,t], 
		$$
		is of class  
		\begin{equation} \label{regularity.backwards}
		\{y_s\} \in C^{1}([0,t);(L^{p}(\Sigma,\rel^n))^*) 
		\cap C^{0}([0,t);D(A_0^*)) \cap
		C^{0}([0,t];(L^{p}(\Sigma,\rel^n))^*),
		\end{equation}
		it is the unique classical solution\footnote{See here Definition 4.1.1 in \cite{Lunardi}.} of the adjoint linear \underline{terminal value problem}:
		\begin{eqnarray} \label{adjoint.backward.uniqueness}
			\partial_s y_s = - \Big{(}D(\Mill_{F_0})(\PP^*(s,0,F_0))\Big{)}^*(y_s) 
			\quad \textnormal{in} \,\,\, (L^{p}(\Sigma,\rel^n))^*,
			\,\, \forall\, s \in [0,t), \qquad \\
			y_t = \xi^* \quad  \textnormal{in} \quad (L^{p}(\Sigma,\rel^n))^*, \qquad \nonumber   
		\end{eqnarray} 
		and its time inverse $\{\tilde y_s\} := \{y_{t-s}\}$ is the unique mild solution \footnote{See here Definition 4.1.4 in \cite{Lunardi}.} of the initial value problem: 
		\begin{eqnarray} \label{adjoint.backward.unique}
			\partial_s \tilde y_s =  \Big{(}D(\Mill_{F_0})
			(\PP^*(t-s,0,F_0))\Big{)}^*(\tilde y_s) \quad \textnormal{in} \,\,\, (L^{p}(\Sigma,\rel^n))^*,
			\,\, \forall\, s \in (0,t], \qquad \\
			\tilde y_0 = \xi^*  \quad  \textnormal{in} \quad (L^{p}(\Sigma,\rel^n))^*. \qquad \nonumber   
		\end{eqnarray}  
	\end{itemize} 
\end{theorem}
\proof 
\begin{itemize} 
	\item[1)] On account of the first part of Theorem \ref{real.analytic.flow} there is a unique, maximal flow line $\{\PP^*(t,0,F_0)\}_{t\in [0,T_{\textnormal{max}}(F_0))}$ of evolution equation (\ref{de_Turck_equation_2}), in the sense of
	Definition \ref{Umbilic.free}, (d), because the 
	immersion $F_0$ was supposed to be smooth and umbilic-free. 
	Furthermore we know from Theorem \ref{Psi.of.class.C_1}, that for any fixed $T\in (0,T_{\textnormal{max}}(F_0))$ there is some open neighborhood $W_{F_0,T,p}$ 
	about the smooth solution
	$\{\PP^*(t,0,F_0)\}_{t \in [0,T]}$ of equation (\ref{de_Turck_equation}) within the space $X_{T}$, such that the operator $\Psi^{F_0,T}$ in line (\ref{Psi}) 
	is a $C^{\omega}$-map from $W_{F_0,T,p}$ to $Y_{T}$, 
	whose Fr\'echet derivative in the particular element 
	$\{\PP^*(t,0,F_0)\}_{t \in [0,T]}\in 
	W_{F_0,T,p} \cap C^{\infty}(\Sigma \times [0,T],\rel^n)$ yields a topological isomorphism between 
	$X_{T}$ and $Y_{T}$. Noting also that there holds 
	$$
	\Psi^{F_0,T}(\{\PP^*(t,0,F_0)\}_{t \in [0,T]})=(F_0,0)\in Y_{T},
	$$ 
	by definition of the operator $\Psi^{F_0,T}$ in (\ref{Psi}) 
	and since $\{\PP^*(t,0,F_0)\}_{t\in [0,T]}$ 
	solves equation (\ref{de_Turck_equation}),  
	we infer from the Inverse Mapping Theorem for non-linear $C^{\omega}$-operators [\cite{Zeidler}, Theorem 4.B and 
	Corollary 4.23], that there is some small open ball 
	$B_{\rho}((F_0,0)) \subset Y_{T}$, with $\rho=\rho(\Sigma,F_0,T,p)>0$, 
	and an appropriate further open neighborhood 
	$W_{F_0,T,p}^* \subset W_{F_0,T,p}$ of 
	$\{\PP^*(t,0,F_0)\}_{t \in [0,T]}$ 
	in $X_{T}$, such that 
	$$
	\Psi^{F_0,T}: W_{F_0,T,p}^* \stackrel{\cong}\longrightarrow B_{\rho}((F_0,0)) 
	$$ 
	is a $C^{\omega}$-diffeomorphism. Hence, by definition of the map $\Psi^{F_0,T}$ the restriction of the inverse map $(\Psi^{F_0,T})^{-1}$ to the product $B_{\rho}^{4-\frac{4}{p},p}(F_0) \times \{0\} \subset B_{\rho}((F_0,0))$ yields exactly the unique maximal 
	solutions $\{\PP^*(t,0,F)\}_{t \in [0,T]}$ 
	of the parabolic Cauchy problems in 
	(\ref{initial.value.problem}) restricted to 
	$\Sigma \times [0,T]$ from the second part of 
	Theorem \ref{short.time.solution}, i.e.:  
	\begin{eqnarray}  \label{C1.Evolution.operator}
		X_{T} \ni \{\PP^*(t,0,F)\}_{t \in [0,T]} 
		= (\Psi^{F_0,T})^{-1}((F,0)) \\  
		\forall \, F \in B_{\rho}^{4-\frac{4}{p},p}(F_0) 
		\subset W^{4-\frac{4}{p},p}(\Sigma,\rel^n), \nonumber
	\end{eqnarray}   
	and this map consequently has to be of class $C^{\omega}$ as a non-linear operator from $W^{4-\frac{4}{p},p}(\Sigma,\rel^n)$ to $X_{T}$.  
	\item[2)] Now we show that the Gateaux derivative of the evolution operator $[F \mapsto \PP^*(\,\cdot\,,0,F)]$ of equation (\ref{de_Turck_equation}) in the fixed umbilic-free immersion $F_0$ has to coincide with the 
	\underline{unique strict solution} of the linear 
	parabolic initial value problem (\ref{Poincare.trick}); 
	compare here to the footnote attached to equation (\ref{Poincare.trick}) regarding this terminology.
	To this end, we firstly fix some function $\xi \in
	W^{4,p}(\Sigma,\rel^n) \setminus \{0\}$ arbitrarily,
	in order to define the auxiliary function
	$$
	y(\,\cdot \,,h) := \PP^*(\,\cdot \,,0,F_0 + h \,\xi) \in X_{T},
	$$
	for sufficiently small $|h| \geq 0$. We can immediately see that there holds $y(0,h) = \PP^*(0,0,F_0 + h \,\xi) = F_0 + h \,\xi$. 
	Recalling statement (\ref{solution.operator.1}) of
	the first part of this theorem we set
	$h_*(T,\xi):=\frac{\rho}{\parallel \xi \parallel_{W^{4,p}(\Sigma,\rel^n)}}$
	for the same $\rho>0$ as in statement (\ref{solution.operator.1}),
	in order to guarantee that the initial immersions $y(0,h)$ satisfy
	$\parallel y(0,h) - F_0 \parallel_{W^{4,p}(\Sigma)} <\rho$
	for every $|h|< h_*(T,\xi)$. Furthermore, we have that
	\begin{equation} \label{derivative.P.0}
		\frac{y(0,h)-y(0,0)}{h} = 
		\frac{F_0 + h \,\xi - F_0}{h} = \xi   
		\quad \textnormal{for} \,\,h \not =0.
	\end{equation}
	On account of the second part of Theorem \ref{Psi.of.class.C_1}
	there is an open neighborhood $W_{F_0,T,p}$ 
	of the flow line $\{y(t,0)\}_{t\in [0,T]} = \{\PP^*(t,0,F_0)\}_{t\in [0,T]}$ 
	in the space $X_{T}$, such that the operator 
	\begin{equation}   \label{Dill.real.analytic}
		\Mill_{F_0}:W_{F_0,T,p} \stackrel{C^{\omega}} \longrightarrow  
		L^p([0,T];L^p(\Sigma,\rel^n)),  \,\, \{f_t\}_{t\in [0,T]} \mapsto  \{\Mill_{F_0}(f_t)\}_{t\in [0,T]}
	\end{equation} 
	is real analytic, and on account of the first part of this theorem the evolution operator 
	$\PP^*(\,\cdot \,,0,\,\cdot \,)$ in line (\ref{solution.operator.1}) maps the open $W^{4-\frac{4}{p},p}$-ball $B^{4-\frac{4}{p},p}_{\rho}(F_0)$ about $F_0$ into the open neighborhood $W_{F_0,T,p}$ of $\{y(t,0)\}_{t\in [0,T]}$ real analytically, thus in particular locally Lipschitz continuously. Hence, we conclude that for every $\varepsilon>0$ there is a $\delta(\varepsilon)>0$, 
	such that  
	\begin{equation} \label{y.converges}
		\parallel y(\,\cdot \,,h) - y(\,\cdot \,,0) 
		\parallel_{X_T} 
		\equiv \parallel \PP^*(\,\cdot \,,0, F_0 + h \,\xi)
		-\PP^*(\,\cdot \,,0,F_0) \parallel_{X_T} 
		< \varepsilon,
	\end{equation}
    for any \,$|h| < \delta(\varepsilon)$.
	Therefore, we can choose some positive 
	$h^* < h_*(T,\xi) = \frac{\rho}{\parallel \xi \parallel_{W^{4,p}(\Sigma,\rel^n)}}$
	that small, such that every convexly combined function
	$[(s,x) \mapsto (\sigma \, y(s,h) + (1-\sigma)\, y(s,0))(x)]$ 
	is contained in the open neighborhood
	$W_{F_0,T,p}$ of the flow line $\{y(t,0)\}_{t\in [0,T]}$
	in $X_{T}$, for every fixed $\sigma \in [0,1]$, provided 
	we have $0 \leq |h| < h^*$. 
	Now, since we know from the second part of 
	Theorem \ref{short.time.solution} respectively from 
	the first part of this theorem, that for every fixed 
	$F \in B^{4-\frac{4}{p},p}_{\rho}(F_0)$ the  
	solution $\{\PP^*(t,0,F)\}_{t\in [0,T]}$ of flow equation (\ref{de_Turck_equation}) is an element of the space 
	$X_T = W^{1,p}([0,T];L^p(\Sigma,\rel^n)) \cap 
	L^p([0,T];W^{4,p}(\Sigma,\rel^n))$, the differential equation (\ref{de_Turck_equation}), i.e.
	$$ 
	\partial_t(\PP^*(\,\cdot\,,0,F)) = \Mill_{F_0}(\PP^*(\,\cdot\,,0,F))    \quad 
	\textnormal{in}  \quad L^p([0,T];L^p(\Sigma,\rel^n)),
	$$ 
	and embedding (\ref{trace.time.embedding}) imply 
	the integral equation 
	\begin{eqnarray}  \label{Integrated.flow}
		\PP^*(t,0,F) - F = \int_{0}^{t}  \Mill_{F_0}(\PP^*(s,0,F)) \, ds
		\quad \textnormal{in}  \quad L^p(\Sigma,\rel^n)
	\end{eqnarray}
	for \underline{every} $t\in [0,T]$ and for every fixed 
	$F \in B^{4-\frac{4}{p},p}_{\rho}(F_0)$. 
	Moreover, we recall here that for any $\eta \in  X_T$ 
	the function 
	$(D(\Mill_{F_0})(\sigma \, y(\,\cdot\,,h) + (1-\sigma)\, y(\,\cdot\,,0)))(\eta)$ is not only an abstract element of 
	the Banach space \\
	$L^p([0,T];L^p(\Sigma,\rel^n))$, but 
	it can be ``concretely computed'' in terms of the right hand 
	side of formula (\ref{Frechet.of.Psi}), at least in $L^1$-almost 
	every $s\in [0,T]$. Recalling now statements 
	(\ref{derivative.P.0}) and (\ref{Dill.real.analytic}) and the fact that we have\,
	$\{\sigma \, y(s,h) + (1-\sigma)\, y(s,0))\}_{s\in [0,T]} \in W_{F_0,T,p}$, for every fixed 
	$\sigma \in [0,1]$ and every $0 \leq |h| < h^*$,  
	we can combine equation (\ref{Integrated.flow}) with the generalized mean value theorem, in order to infer for the difference quotients
	$$
	z_h(\,\cdot\,) := \frac{1}{h}\, 
	\big{(}y(\,\cdot\,,h) - y(\,\cdot\,,0)\big{)}
	\in X_{T}, \,\,\, \textnormal{being defined for} \,\, 0<|h|<h^*,
	$$
	the crucial equation
	\begin{eqnarray}  \label{z_h_t}
		z_h(t) = \frac{y(0,h) - y(0,0)}{h}
		+   \int_{0}^{t}  \frac{1}{h}
		\Big{(} \Mill_{F_0}(y(s,h)) - \Mill_{F_0}(y(s,0))  \Big{)} \, ds  \qquad\quad  \\
		= \xi + \int_{0}^t  \int_0^1
		D(\Mill_{F_0})\big{(} \sigma \, y(s,h) + (1-\sigma)\, y(s,0)\big{)}.
		\Big{(} \frac{y(s,h) - y(s,0)}{h} \Big{)} \, d\sigma \,ds      \nonumber \\
		\equiv  \xi  +  \int_{0}^t  \Big{(} \int_0^1
		D(\Mill_{F_0})\big{(} \sigma \, y(s,h) + (1-\sigma)\, y(s,0)\big{)}   \, d\sigma \Big{)} \,\,
		\Big{(} \frac{y(s,h) - y(s,0)}{h} \Big{)}\, \,ds                           \nonumber    \\
		\equiv \xi + \int_{0}^t N_{F_0,h\xi}(s).(z_h(s)) \, ds \quad \textnormal{on}
		\,\,\, \Sigma, \nonumber
	\end{eqnarray}
	for every $t \in [0,T]$ and for $0<|h|<h^*$. Here, we have 
	used the abbreviation
	\begin{eqnarray*}
		N_{F_0, h\xi}(s):= \int_0^1 D(\Mill_{F_0})
		\big{(}\sigma \, y(s,h) + (1-\sigma)\, y(s,0)\big{)} \,d\sigma, 
	\end{eqnarray*}
	for $0<|h|<h^*$, which again has to 
	be interpreted pointwise in $L^1$-almost every $s\in [0,T]$ 
	by means of the right hand side of formula (\ref{Frechet.of.Psi}).
	Now, using again statement (\ref{y.converges})
	and the fact from line (\ref{Dill.real.analytic}), that
	$\Mill_{F_0}: W_{F_0,T,p} \longrightarrow  
	L^p([0,T];L^p(\Sigma,\rel^n))$ is an operator of class $C^{\omega}$, we infer that for every $\varepsilon >0$ 
	there is some $\tilde \delta(\varepsilon)>0$, such that
	\begin{eqnarray} \label{operator.convergence}
		\parallel N_{F_0, h\xi}(\,\cdot\,) -
		D(\Mill_{F_0})(y(\,\cdot\,,0))
		\parallel_{\Lift(X_T,L^p([0,T],
		L^p(\Sigma,\rel^n)))}        \quad    \nonumber \\
		\leq \int_{0}^1  \parallel
		D(\Mill_{F_0})\big{(}\sigma \, y(\,\cdot \,,h) + (1-\sigma)\, y(\,\cdot \,,0)\big{)}  \qquad \quad \qquad \\  
		- D(\Mill_{F_0})(y(\,\cdot \,,0))
		\parallel_{\Lift(X_T,L^p([0,T],L^p(\Sigma,\rel^n)))} \, d\sigma     < \varepsilon                  \nonumber
	\end{eqnarray}
	holds, provided $0<|h|<\tilde \delta(\varepsilon)$. 
	Now, again using the result of the first part of this theorem, i.e. statement (\ref{solution.operator.1}), we can immediately infer from the definition of the family $\{z_h\}$ and from the 
	chain rule, that $\{z_h\}$ converges in $X_T$ to the Fr\'echet derivative $D_F\PP^*(\,\cdot \,,0,F_0).(\xi)$ applied to $\xi$ respectively to the Gateaux derivative $D_{\xi}\PP^*(\,\cdot \,,0,F_0)$ w.r.t. $\xi$ of the evolution operator 
	$[F \mapsto \{\PP^*(t,0,F)\}_{t \in [0,T]}]$ 
	as $h\to 0$:
	\begin{eqnarray}   \label{Convergence.1}
		\lim_{h \to 0} z_h  \equiv \lim_{h \to 0} \frac{y(s,h) - y(s,0)}{h}     \nonumber\\
		=  \lim_{h \to 0}
		\Big{(} \Big{[}(x,t) \mapsto
		\frac{\PP^*(t,0,F_0+h\,\xi)(x) - \PP^*(t,0,F_0)(x)}{h} \Big{]}  \Big{)}   \nonumber \\
		\equiv [(x,t) \mapsto  D_{\xi}\PP^*(t,0,F_0)(x)]       
		=: z^* \qquad \textnormal{in} \quad X_T.   \qquad                  
	\end{eqnarray}
	We note here that it is an immediate consequence of the above construction of the function $z^*$ and of embedding (\ref{trace.time.embedding}), that there holds
	\begin{equation}  \label{Gateaux.home}
	[(x,t) \mapsto D_{\xi}\PP^*(t,0,F_0)(x)] \equiv z^*
	\in C^0([0,T],W^{4-\frac{4}{p},p}(\Sigma,\rel^n)). 
	\end{equation}
	Inserting the convergences (\ref{operator.convergence}) and (\ref{Convergence.1}) into formula (\ref{z_h_t}), we infer the following integral equation from the triangle inequality for ``Bochner's integral'' and from H\"older's inequality, when letting tend $h \searrow 0$:
	\begin{eqnarray}   \label{Integral.equation.limit}
		z^*_t =  \xi  +  
		\int_{0}^t D(\Mill_{F_0})(y(s,0)).(z^*_s) \, ds   \qquad  \\
		\equiv \xi + \int_{0}^t D(\Mill_{F_0})(\PP^*(s,0,F_0)).(z^*_s) \, ds \quad  \textnormal{in}  \quad L^p(\Sigma,\rel^n)       \nonumber
	\end{eqnarray}
	for \underline{every} $t \in [0,T]$, taking here also (\ref{Gateaux.home}) into account. Since we know already from convergence (\ref{Convergence.1}), that the limit function $z^*$ is an element of $X_T$, 
	and since we know from the second part of Theorem \ref{Psi.of.class.C_1}, that the Fr\'echet derivative 
	$D(\Mill_{F_0})(\PP^*(\,\cdot\,,0,F_0))$ maps the Banach space $X_T$ continuously into $L^{p}([0,T];L^p(\Sigma,\rel^n))$, we may differentiate equation (\ref{Integral.equation.limit}) 
	w.r.t. $t$ and obtain together with computation (\ref{D_Psi_f}) 
	an equivalent reformulation of equation (\ref{Integral.equation.limit}), 
	stating that the function $[(x,t) \mapsto z^*_t(x)] \equiv
	[(x,t) \mapsto D_{\xi}\PP^*(t,0,F_0)(x)] \in X_T$ from line (\ref{Convergence.1}) solves the following linear parabolic differential system of equations:
	\begin{eqnarray}   \label{linear.initial.value.problem}
		\partial_t (z^*_t) 
		=D(\Mill_{F_0})(\PP^*(t,0,F_0)).(z^*_t)     \quad \qquad  \\                
		= - \frac{1}{2} \mid A^0_{\PP^*(t,0,F_0)} \mid^{-4} \,
		g^{ij}_{\PP^*(t,0,F_0)} \, g^{kl}_{\PP^*(t,0,F_0)} \, 
		\nabla_i^{F_0} \nabla_j^{F_0} \nabla_k^{F_0} \nabla_l^{F_0}(z^*_t)  
		\nonumber \\ 
		-B_{3}^{ijk}(\,\cdot\,,D_x\PP^*(t,0,F_0),D_x^2\PP^*(t,0,F_0)) \cdot \nabla^{F_0}_{ijk}(z^*_t) \nonumber \\ 
		-B_2^{ij}(\,\cdot\,,D_x\PP^*(t,0,F_0),D_x^2\PP^*(t,0,F_0),D_x^3\PP^*(t,0,F_0),
		D_x^4\PP^*(t,0,F_0)) \cdot \nabla^{F_0}_{ij}(z^*_t)       \nonumber \\
		-B_1^i(\,\cdot\,,D_x\PP^*(t,0,F_0),D_x^2\PP^*(t,0,F_0),D_x^3\PP^*(t,0,F_0),
		D_x^4\PP^*(t,0,F_0)) \cdot \nabla^{F_0}_{i}(z^*_t)                \nonumber 
	\end{eqnarray} 
	in $L^{p}([0,T];L^p(\Sigma,\rel^n))$, with $z^*_0 = \xi$ on $\Sigma$, if $\xi \in W^{4,p}(\Sigma,\rel^n)$, just as asserted in line (\ref{Poincare.trick}). Since we know that 
	$\{\PP^*(t,0,F_0)\}_{t\in [0,T]}\in 
	C^{\infty}(\Sigma \times [0,T],\rel^n)$,
	we infer from the fourth part of Theorem \ref{Psi.of.class.C_1}, that the Fr\'echet derivative $D(\Psi^{F_0,T})(\PP^*(\,\cdot\,,0,F_0))$ is a topological isomorphism between $X_T$ and $Y_T$, which implies that the function $[(x,t) \mapsto z^*_t(x)] \equiv [(x,t) \mapsto D_{\xi}\PP^*(t,0,F_0)(x)]$ in line (\ref{Convergence.1}) has to 
	be the \underline{unique solution} of initial value problem (\ref{linear.initial.value.problem}) within the space $X_T$. Furthermore, we know already from the proof 
	of the fourth part of Theorem \ref{Psi.of.class.C_1}, 
	that the linear operator
	``$\partial_t - D(\Mill_{F_0})(\PP^*(\,\cdot\,,0,F_0))$'' 
	satisfies all requirements of Proposition \ref{A.priori.estimate}.
	Hence, Proposition \ref{A.priori.estimate} guarantees us, 
	that the linear differential operators 
	\begin{eqnarray} \label{A.t}
	\omega \,\textnormal{Id}_{L^p(\Sigma,\rel^n)} - A_t 
	:= \omega \,\textnormal{Id}_{L^p(\Sigma,\rel^n)} 
	- D(\Mill_{F_0})(\PP^*(t,0,F_0)): \nonumber \\
	W^{4,p}(\Sigma,\rel^n)  \longrightarrow  L^p(\Sigma,\rel^n) 
	\end{eqnarray}
	are of type $(M,\vartheta)$ for every fixed $t \in [0,T]$, where the angle $\vartheta \in (\pi/2,\pi)$ is arbitrarily fixed, and where $M=M(\Sigma,F_0,T,p,\vartheta)\geq 1$ and $\omega=\omega(\Sigma,F_0,T,p,\vartheta)$ $>0$ are sufficiently large real numbers. In particular, the differential operators 
	$A_t \equiv D(\Mill_{F_0})(\PP^*(t,0,F_0))$ in line (\ref{A.t}) 
	are ``sectorial'' in $L^p(\Sigma,\rel^n)$ - in the sense of Definition 2.0.1 in \cite{Lunardi} - on account of estimate (\ref{sectorial.estimate}), uniformly for every $t \in [0,T]$.
	Since the operators $A_t$ also depend smoothly on $t$, the family of linear operators $\{A_{t}\}$ satisfies all principal hypotheses of Section 6.1 in \cite{Lunardi} for the couple of 
	Banach spaces $X := L^p(\Sigma,\rel^n)$ and $D := W^{4,p}(\Sigma,\rel^n)$, adopting here the notation 
	of Chapter 6 in \cite{Lunardi}. Hence, we may apply the 
	results of Sections 6.1 and 6.2 of \cite{Lunardi},
	and we thus obtain a ``parabolic fundamental solution'' respectively an ``evolution operator'' $\{G(t,s)\} \equiv \{G^{F_0}(t,s)\}\subset \Lift(X,X)$, for $0 \leq s \leq t \leq T$, which is characterized by the following $5$ properties:
	\begin{itemize}
		\item[1)] There exists some constant $C(\Sigma,F_0,T,p)>0$, such that
		\begin{equation}  \label{long.time.esimate.G}
			\parallel G(t,s).(\xi) \parallel_{L^p(\Sigma,\rel^n)} \leq
			C(\Sigma,F_0,T,p) \, \parallel \xi \parallel_{L^p(\Sigma,\rel^n)},
		\end{equation}
		for every $\xi \in L^p(\Sigma,\rel^n)$ and for every 
		$0  \leq  s  \leq  t \leq T$. In particular, the operator 
		$G(t,s)$ yields a linear continuous map
		\begin{equation}  \label{fundament.G.map}
			G(t,s): L^p(\Sigma,\rel^n)	\longrightarrow L^p(\Sigma,\rel^n),
		\end{equation}
		for $0 \leq  s \leq t \leq T$; and it has the additional property that
		$\textnormal{range}(G(t,s)) \subset W^{4,p}(\Sigma,\rel^n)$, if $s<t$.
		\item[2)] For every $0  \leq  r  \leq  s  \leq  t \leq T$ there holds the semigroup property:
		\begin{equation}  \label{semigroup.property}
			G(t,s) \circ G(s,r) = G(t,r)
		\end{equation}
		and also $G(s,s) = \textnormal{Id}$ on $L^p(\Sigma,\rel^n)$,
		$\forall \, s\in [0,T]$.
		\item[3)] For every fixed $s \in [0,T)$ and for every $\xi \in L^p(\Sigma,\rel^n)$
		the function $[t \mapsto G(t,s).(\xi)]$ is continuous on $[s,T]$ and continuously differentiable in $(s,T]$, and
		the derivative $\partial_t G(t,s).(\xi)$ satisfies the linear differential equation
		\begin{eqnarray}  \label{generate}
			\partial_t G(t,s).(\xi) = (A_t \circ G(t,s)).(\xi),  
			\qquad \forall \, t \in (s,T].
		\end{eqnarray}
		\item[4)] For every fixed $t \in (0,T]$ and every $\xi \in W^{4,p}(\Sigma,\rel^n)$ the function $[s \mapsto G(t,s).(\xi)]$ 
		is continuous on $[0,t]$ and 
		continuously differentiable in $[0,t)$, and the derivative 
		$\partial_s G(t,s).(\xi)$ satisfies the linear differential equation
		\begin{equation}  \label{generate.2}
			\partial_s G(t,s).(\xi)  = - (G(t,s) \circ A_s).(\xi)  \qquad \forall \, s\in [0,t).
		\end{equation}
		\item[5)] There is some constant $C = C(\Sigma,F_0,T,p)>0$ such that there holds the estimate
		\begin{equation}  \label{long.time.esimate.G.2}
			\parallel G(t,s).(\xi) \parallel_{W^{4,p}(\Sigma,\rel^n)}
			\leq \frac{C}{t-s} \,  \parallel \xi \parallel_{L^p(\Sigma,\rel^n)},
		\end{equation}
		for every $\xi \in L^p(\Sigma,\rel^n)$, 
		for $0 \leq s < t \leq T$.
	\end{itemize}
	It is important to mention here, that the above ``parabolic fundamental solution'' $G^{F_0}(t,s)$ is actually \underline{uniquely determined} by the above $5$ properties according to Section 3 in \cite{Amann.1} respectively Theorem 4.4.1 of Chapter II in \cite{Amann.1995}. 
	Moreover, we can infer from the above properties of the 
	differential operators $A_t:= D(\Mill_{F_0})(\PP^*(t,0,F_0))$ 
	and from Corollary 6.1.9 and Proposition 6.2.2 in \cite{Lunardi}, that the above ``parabolic 
	fundamental solution'' $[(x,t) \mapsto (G^{F_0}(t,0).\xi)(x)]$ is the \underline{unique strict solution} of the initial value problem 
	\begin{eqnarray}    \label{initial.value.A}
	\partial_t(z^*_t) =  D(\Mill_{F_0})(\PP^*(t,0,F_0)).(z^*_t) \equiv
	(A_t).(z^*_t),  \quad \textnormal{for}  \,\, 
	t \in [0,T],             \nonumber   \\
	z_{0} = \xi  \in W^{4,p}(\Sigma,\rel^n) 
	\end{eqnarray}
	w.r.t. the Banach spaces $X:=L^{p}(\Sigma,\rel^n)$ and $D:= W^{4,p}(\Sigma,\rel^n)$, i.e. 
	that the mild solution $[(x,t) \mapsto (G^{F_0}(t,0).\xi)(x)]$ of problem (\ref{initial.value.A}) satisfies: 
	\begin{equation}  \label{regularity.alpha}
		[(x,t) \mapsto (G^{F_0}(t,0).\xi)(x)] \in C^1([0,T];L^p(\Sigma,\rel^n)) \cap C^{0}([0,T];W^{4,p}(\Sigma,\rel^n))
	\end{equation}
	and that there is no further such solution of initial value problem (\ref{initial.value.A}).  
	Now, since the Banach space in (\ref{regularity.alpha}) is contained in the space $X_T$, and since the function $[(x,t) \mapsto z^*_t(x)] \equiv [(x,t) \mapsto D_{\xi}\PP^*(t,0,F_0)(x))]$ solves the same initial value problem \underline{uniquely} in $X_T$, we can conclude the identity
	\begin{equation}   \label{Gateaux.equat.G}
		D_{\xi}\PP^*(t,0,F_0) = G^{F_0}(t,0).\xi   \quad \textnormal{on}
		\quad  \Sigma, \quad \forall \, t \in [0,T],
	\end{equation}
	for every $\xi\in W^{4,p}(\Sigma,\rel^n)$,  
	i.e. that the Gateaux derivative $[(x,t) \mapsto D_{\xi}\PP^*(t,0,F_0)(x)]$
	in direction of any fixed function $\xi \in W^{4,p}(\Sigma,\rel^n)$
	is given by the application of the unique parabolic fundamental solution $G^{F_0}(t,0)$ to any initial function $\xi \in W^{4,p}(\Sigma,\rel^n)$, 
	which is generated by the smooth family of uniformly elliptic differential operators $\{D(\Mill_{F_0})(\PP^*(t,0,F_0))\}_{t\in [0,T]}$ in (\ref{linear.initial.value.problem}) and (\ref{A.t}). Combining statements 
	(\ref{regularity.alpha}) and (\ref{Gateaux.equat.G}) we 
	verify the asserted uniqueness of the solution \\ 
	$[(x,t) \mapsto D_{\xi}\PP^*(t,0,F_0)(x)]$ of equation (\ref{Poincare.trick}) and its asserted regularity: $[(x,t) \mapsto D_{\xi}\PP^*(t,0,F_0)(x)] \in 
	C^1([0,T];L^p(\Sigma,\rel^n)) \cap C^{0}([0,T];W^{4,p}(\Sigma,\rel^n))$, 
	improving the previous regularity statement 
	(\ref{Gateaux.home}) significantly. 
	\item[3)] Combining the second property of the parabolic 
	fundamental solution \\
	$\{G^{F_0}(t,s)\}$ with statements (\ref{long.time.esimate.G.2}) and (\ref{Gateaux.equat.G}) 
	we see that 
	\begin{equation}  \label{mapping.G}
		D_F\PP^*(t,0,F_0) = G^{F_0}(t,0): W^{4,p}(\Sigma,\rel^n) 
		\longrightarrow W^{4,p}(\Sigma,\rel^n)
	\end{equation}
	is a continuous operator for every fixed $t \in [0,T]$, 
	and that this linear operator can be uniquely extended by 
	means of statement (\ref{fundament.G.map}) to a continuous linear operator 
	$$
	G^{F_0}(t,0):L^p(\Sigma,\rel^n) \longrightarrow L^{p}(\Sigma,\rel^n),
	$$
	for every $t \in [0,T]$, whose range is contained in $W^{4,p}(\Sigma,\rel^n)$ for every $t\in (0,T]$. We shall 
	term this unique extension 
	$G^{F_0}(t,0)$ again ``$D_F\PP^*(t,0,F_0)$'' in the sequel. 
	\item[4)] We have already shown in the second part of this 
	theorem, that the family of linear operators 
	$A_t := D(\Mill_{F_0})(\PP^*(t,0,F_0))$, $t\in [0,T]$, generates a unique parabolic fundamental solution $G^{F_0}(t,s)$, and that the function in identity  (\ref{Gateaux.equat.G}) is the unique strict solution to 
	initial value problem (\ref{initial.value.A}), for any initial function $\xi \in W^{4,p}(\Sigma,\rel^n)$. Since the immersion $\PP^*(s,0,F_0)$ 
	is of class $C^{\infty}(\Sigma,\rel^n)$ 
	for every $s \in [0,T]$ on account of the first part 
	of Theorem \ref{real.analytic.flow}, we may interchange the immersion $F_0$ with the immersion $\PP^*(s,0,F_0)$ and the initial time $0$ with another initial time $s \in [0,T]$ in the first part of this theorem, and we infer that the evolution operator $[F \mapsto \PP^*(\,\cdot\,,s,F)]$ is of class $C^{\omega}$ in a small ball 
	$B_{\rho}^{4-\frac{4}{p},p}(\PP^*(s,0,F_0))$ about $\PP^*(s,0,F_0)$
	in $W^{4-\frac{4}{p},p}(\Sigma,\rel^n)$.   
	Moreover, we infer from the second part of this theorem 
	that the resulting Fr\'echet derivative $D_{F}\PP^*(t,s,\PP^*(s,0,F_0))$
	of the evolution operator $[F \mapsto \PP^*(\,\cdot\,,s,F)]$ in the 
	immersion $\PP^*(s,0,F_0)$ yields functions 
	$[(x,t) \mapsto D_{F}\PP^*(t,s,\PP^*(s,0,F_0)).\xi(x)]$, 
	for $\xi \in W^{4,p}(\Sigma,\rel^n)$, which are the unique strict solutions of the initial value problem  
	\begin{eqnarray}    \label{initial.value.A.2}
	\partial_t(z^*_t) =  D(\Mill_{F_0})(\PP^*(t,s,\PP^*(s,0,F_0))).
	(z^*_t), \quad \textnormal{for}  \,\, t \in [s,T],             \nonumber   \\
	z_{s} = \xi  \in W^{4,p}(\Sigma,\rel^n)
	\end{eqnarray}
	w.r.t. the B-spaces $X:= L^{p}(\Sigma,\rel^n)$ and $D:= W^{4,p}(\Sigma,\rel^n)$. Because of the identity 
	$$ 
	\PP^*(t,0,F_0) = \PP^*(t-s,0,\PP^*(s,0,F_0))=
	\PP^*(t,s,\PP^*(s,0,F_0)),  \,\,\,
	\forall \,0\leq s \leq t \leq T,
	$$
	we have here:
	\begin{equation}  \label{mugging}
	D(\Mill_{F_0})(\PP^*(t,s,\PP^*(s,0,F_0)))
	= D(\Mill_{F_0})(\PP^*(t,0,F_0))  \equiv A(t), \qquad 
	\end{equation} 
	$\forall \, 0\leq s\leq t\leq T$.          
	Hence, fixing $s \in [0,T]$ initial value problem (\ref{initial.value.A.2}) becomes the problem
	\begin{eqnarray}    \label{initial.value.A.3}
		\partial_t(z^*_t) = A(t).(z^*_t), \quad \textnormal{for}  
		\,\, t \in [s,T],       \quad
		z_{s} = \xi  \in W^{4,p}(\Sigma,\rel^n)
	\end{eqnarray}
	whose unique strict solution with values in $W^{4,p}(\Sigma,\rel^n)$ 
	is given by the functions $[(x,t) \mapsto G^{F_0}(t,s).\xi(x)]$, 
	since we know from the proof of the second part of this theorem, that our particular parabolic fundamental solution $\{G^{F_0}(t,s)\}$ 
	is uniquely generated by the family of operators 
	$\{A_t\}_{t\in [0,T]}$, so that we can argue here again by means of property (3) of $\{G^{F_0}(t,s)\}$ combined with Corollary 6.1.9 and Proposition 6.2.2 in \cite{Lunardi}. Hence, combining this with the fact, that $[(x,t) \mapsto D_{F}\PP^*(t,s,\PP^*(s,0,F_0)).(\xi)(x)]$ is the unique strict solution of initial value problem 
	(\ref{initial.value.A.3}) with values in $W^{4,p}(\Sigma,\rel^n)$ as well, we obtain the identity
	\begin{equation}  \label{Frechet.G}
		D_{F}\PP^*(t,s,\PP^*(s,0,F_0)).(\xi) = G^{F_0}(t,s).(\xi),  \quad 
		\forall \,0\leq s \leq t \leq T,
	\end{equation} 
	for any $\xi \in W^{4,p}(\Sigma,\rel^n)$, which proves  
	- again on account of estimate (\ref{long.time.esimate.G.2}) - that the linear operators $G^{F_0}(t,s):L^p(\Sigma,\rel^n) \longrightarrow L^p(\Sigma,\rel^n)$ in (\ref{fundament.G.map}) are the unique continuous linear extensions of the linear operators
	$$ 
	D_F\PP^*(t,s,\PP^*(s,0,F_0)): 
	W^{4,p}(\Sigma,\rel^n) \longrightarrow W^{4,p}(\Sigma,\rel^n),
	$$
	induced by the Fr\'echet derivative $D_F\PP^*(\,\cdot\,,s,\PP^*(s,0,F_0))$, 
	from $W^{4,p}(\Sigma,\rel^n)$ to $L^p(\Sigma,\rel^n)$, 
	for every $0\leq s \leq t \leq T$. 
	Finally, combining this insight with formula (\ref{semigroup.property}), we verify immediately, that the family of linear operators 
	$[(t,s) \mapsto D_F\PP^*(t,s,\PP^*(s,0,F_0))]$ 
	has the asserted semigroup property (\ref{semigroup.property.b})   
	both on $W^{4,p}(\Sigma,\rel^n)$ and on $L^{p}(\Sigma,\rel^n)$. 
	\item[5)] For ease of notation we introduce the abbreviation   
	\begin{equation}  \label{new.starting.time}
		D_{F}\PP^*(t_2,t_1,F_0) := D_{F}\PP^*(t_2,t_1,\PP^*(t_1,0,F_0))
	\end{equation}
	for every pair of times $t_2 \geq t_1$ in $[0,T]$, 
	and we know from the fourth part of this theorem, that 
	these are linear operators mapping
	$W^{4,p}(\Sigma,\rel^n)$ into itself - induced by the unique strict solutions of Cauchy problem 
	(\ref{initial.value.A.2}) - which can be uniquely and continuously extended to $L^{p}(\Sigma,\rel^n)$
	by the parabolic fundamental solution $\{G^{F_0}(t,s)\}_{s\leq t}$ to Cauchy problem (\ref{initial.value.A.3}), so that we shall simply write in the sequel: 
	\begin{equation}  \label{Frechet.L.p}
		D_{F}\PP^*(t_2,t_1,F_0): 
		L^{p}(\Sigma,\rel^n) \longrightarrow L^{p}(\Sigma,\rel^n)
	\end{equation}
	for $t_2 \geq t_1$ in $[0,T]$. 
	We shall note here, that the notation introduced in line
	(\ref{new.starting.time}) is a slight abuse of our original notation, but this abbreviation will not cause any confusion in the sequel, because from the very beginning the flow line $\PP^*(\,\cdot\,,0,F_0)$ of equation (\ref{de_Turck_equation}) was supposed to start moving in the immersion $F_0$ at time $t=0$, and not at any later time $t_1>0$. Now, since the linear operators in line (\ref{Frechet.L.p}) are well-defined in entire
	$L^{p}(\Sigma,\rel^n)$ and continuous w.r.t.
	the $L^{p}$-norm on both sides of (\ref{Frechet.L.p}), 
	they are especially closed linear operators in 
	$L^{p}(\Sigma,\rel^n)$, for every pair $t_2 \geq t_1$ in $[0,T]$.
	Theorem 5.29 in Chapter III of \cite{Kato} therefore
	guarantees, that every operator in line (\ref{Frechet.L.p}) 
	has a uniquely determined maximal, densely defined 
	and closed linear adjoint operator in $(L^{p}(\Sigma,\rel^n))^*$:
	\begin{equation} \label{adjoint}
		(D_{F}\PP^*(t_2,t_1,F_0))^*:
		(L^{p}(\Sigma,\rel^n))^* \longrightarrow (L^{p}(\Sigma,\rel^n))^* 
		\cong L^{q}(\Sigma,\rel^n),
	\end{equation}
	for $q = \frac{p}{p-1}$.   
	Moreover we know, that the linear operators 
	$A_{t}$ from statements (\ref{A.t}) and (\ref{mugging}) depend smoothly on the time $t\in [0,T]$, that they are uniformly sectorial in the Banach space $L^p(\Sigma,\rel^n)$ on account of statement (\ref{A.t}) and that the domain of $A_{t}$, which is $W^{4,p}(\Sigma,\rel^n)$, is densely contained in $L^p(\Sigma,\rel^n)$. In particular, the operators $A_{t}$ are closed linear operators in $L^p(\Sigma,\rel^n)$. 
	Therefore, again Theorem 5.29 in Chapter III of \cite{Kato} guarantees, that $A_{t}$ has a uniquely determined maximal, densely defined and closed linear adjoint operator 
	\begin{eqnarray}  \label{domain.A.star}  
		(A_t)^*:D((A_t)^*) \subset (L^p(\Sigma,\rel^n))^* 
		\longrightarrow  (L^{p}(\Sigma,\rel^n))^*,
	\end{eqnarray}  
	for every $t\in [0,T]$, and furthermore the adjoint operators in (\ref{domain.A.star}) constitute a smooth family of uniformly sectorial linear operators in 
	$(L^p(\Sigma,\rel^n))^* \cong L^q(\Sigma,\rel^n)$  
	on account of Proposition 1.2.3 in Chapter I of \cite{Amann.1995}.
	Now we use the adjoint operators in (\ref{adjoint}), 
	in order to define the family of functions $\{y_s\} \in C^{0}([0,t],(L^{p}(\Sigma,\rel^n))^*)$ by
	\begin{equation}   \label{y.s}
		y_s := \Big{(} D_{F}\PP^*(t,s,F_0) \Big{)}^*(\xi^*),
		\quad \textnormal{for} \,\,\,s \in [0,t].
	\end{equation} 
	It is important to note here, that the fourth property of 
	the parabolic fundamental solution 
	$[(t,s) \mapsto G^{F_0}(t,s)]$ 
	together with equation (\ref{Frechet.G}) actually guarantees, that the above defined family of functions $\{y_s\}$ is of class $C^{0}([0,t];(L^{p}(\Sigma,\rel^n))^*)$. 
	Since the family of operators $\{D_{F}\PP^*(t_2,t_1,F_0)\}_{t_1\leq t_2}$ in line (\ref{Frechet.L.p}) has the semigroup property
	(\ref{semigroup.property.b}), the family of functions
	$\{y_s\}$ evolves according to the evolution
	of the ``adjoint Fr\'echet derivatives'' of the
	evolution operator $\PP^*(\,\cdot \,,\,\cdot \,,F_0)$ of the modification (\ref{de_Turck_equation}) of the MIWF:
	\begin{eqnarray}   \label{move.it}
		\Big{(} D_{F}\PP^*(\tilde s,s,F_0) \Big{)}^*(y_{\tilde s})
		= \Big{(}D_{F}\PP^*(\tilde s,s,F_0)\Big{)}^* \circ
		\Big{(} D_{F}\PP^*(t,\tilde s,F_0) \Big{)}^*(\xi^*)         \qquad\,\,   \\
		= \Big{(} D_{F}\PP^*(t,\tilde s,F_0) \circ 
		D_{F}\PP^*(\tilde s,s,F_0)\Big{)}^*(\xi^*)           =  \qquad    \nonumber      \\
		= \Big{(} D_{F}\PP^*(t,s,F_0) \Big{)}^*(\xi^*)
		= y_s	\quad  \textnormal{for} \,\,\, 0 \leq s \leq \tilde s < t.      \nonumber
	\end{eqnarray}
	Since we know from the fourth part of this theorem 
	- formula (\ref{Frechet.evol.operator.s}) - that the family of Fr\'echet derivatives 
	$\{D_{F}\PP^*(\tilde s,s,F_0)\}_{s \leq \tilde s}$ is the unique ``parabolic fundamental solution'' $\{G(\tilde s,s)\}_{s \leq \tilde s}$ for the linear system (\ref{Poincare.trick}), we can apply here formula (\ref{generate.2}), in order to conclude that for every fixed $y \in (L^{p}(\Sigma,\rel^n))^*$ the function
	$s \mapsto \big{(} D_{F}\PP^*(\tilde s,s,F_0) \big{)}^*(y)$ is continuously differentiable on every interval $[0,\tilde s)$ with $\tilde s <t$ and that it satisfies the linear differential equation
	\begin{equation}  \label{Deriv.Evolution.2nd.comp}
		\partial_s \Big{(} 
		\Big{(} D_{F}\PP^*(\tilde s,s,F_0)\Big{)}^*(y) \Big{)}
		=- \Big{(} D_{F}\PP^*(\tilde s,s,F_0)\circ A_{s} \Big{)}^*(y), \,\,\,  \forall \, s \in [0,\tilde s).
	\end{equation}
	We can therefore take the derivative of
	the family of functions $\{y_s\}$ w.r.t. $s$ in equation (\ref{move.it}), and we infer from a combination of equation (\ref{Deriv.Evolution.2nd.comp}),
	here applied to the function $y:=y_{\tilde s} \in (L^{p}(\Sigma,\rel^n))^*$, again with equation (\ref{move.it}) and with statement (\ref{domain.A.star}), 
	that the family $\{y_s\}$ is of class
	$C^{1}([0,t);(L^{p}(\Sigma,\rel^n))^*)\\ 
	\cap C^{0}([0,t];(L^{p}(\Sigma,\rel^n))^*)$
	and that $\{y_s\}$ is a solution of the following linear 
	differential equation:
	\begin{equation} \label{adjoint.evolution}
		\partial_s y_s
		= \partial_s \Big{(}\Big{(} 
		D_{F}\PP^*(\tilde s,s,F_0)\Big{)}^*(y_{\tilde s}) \Big{)}
		= - \big{(} A_{s} \big{)}^*(y_s),   \quad \forall\, 
		s \in [0, \tilde s),
	\end{equation}
	for every $\tilde s<t$,
	proving that $\{y_s\}$ is actually a mild solution of 
	equation (\ref{adjoint.backward.uniqueness}) on $[0,t)$.
	We note, that by equation (\ref{adjoint.evolution})  
	$y_s$ is not only contained in $(L^{p}(\Sigma,\rel^n))^*$, but even in the smaller space $D(A_s^*) \cong D(A_0^*)$ for  
	every $s \in [0,t)$. Now we consider the functions 
	\begin{equation}   \label{mild.solution}
		\tilde y_s:= y_{t-s}= 
		\Big{(} D_{F}\PP^*(t,t-s,F_0) \Big{)}^*(\xi^*), 
		\quad \textnormal{for}\,\,\, s\in [0,t]. 
	\end{equation}
	Obviously, we know 
	already that $\{\tilde y_s\}$ is of class  
	$C^{1}((0,t];(L^{p}(\Sigma,\rel^n))^*)\cap
	C^{0}([0,t];(L^{p}(\Sigma,\rel^n))^*)$, and we infer from (\ref{adjoint.evolution}) that it solves 
	the initial value problem:  
	\begin{equation}  \label{sectorial.evolution} 
		\partial_s \tilde y_s = 
		\big{(} A_{t-s} \big{)}^*(\tilde y_s)   
		\quad \forall\, s \in (0,t],  \quad \textnormal{and} 
		\quad \tilde y_0 = \xi^*,
	\end{equation}  
	i.e. that $\{\tilde y_s\}$ is the \underline{unique mild solution} of this initial value problem, as asserted in (\ref{adjoint.backward.unique}), taking also its definition in (\ref{mild.solution}) into account, and moreover that 
	$\tilde y_s \in D((A_{t-s})^*) \cong D(A_0^*)$ for every 
	$s\in (0,t]$. We know already that the adjoint linear operators
	in (\ref{domain.A.star}) constitute a $C^{\infty}$-smooth 
	family of densely defined and uniformly sectorial linear operators in $(L^p(\Sigma,\rel^n))^*\cong L^q(\Sigma,\rel^n)$. 
	Hence, the family of linear operators 
	$\{(A_{t-s})^*\}_{s\in [0,t]}$ in (\ref{sectorial.evolution})
	satisfies all requirements of Section 6.1 in \cite{Lunardi}, and we may thus apply Lemma 6.2.1 in \cite{Lunardi} to 
	the mild solution $\{\tilde y_s\}$ 
	of problem (\ref{sectorial.evolution}), in order to obtain 
	even better regularity of $\{\tilde y_s\}$, namely: 
	\begin{equation}  \label{regularity.backwards.2}
		\{\tilde y_s\} \in C^{1}((0,t];(L^{p}(\Sigma,\rel^n))^*) \cap 
		C^{0}((0,t];D(A_0^*)) \cap
		C^{0}([0,t];(L^{p}(\Sigma,\rel^n))^*).
	\end{equation}      
	This means exactly that $\{\tilde y_s\}$ is a 
	classical solution of Cauchy problem (\ref{sectorial.evolution}), 
	in the sense of Definition 4.1.1 in \cite{Lunardi}. 
	We can therefore improve the above regularity statement about the original family of functions $\{y_s\}$, namely: 
	$$
	\{y_s\} \in C^{1}([0,t);(L^{p}(\Sigma,\rel^n))^*) \cap 
	C^{0}([0,t);D(A_0^*)) \cap
	C^{0}([0,t];(L^{p}(\Sigma,\rel^n))^*),
	$$      
	which means that $\{y_s\}$ is a classical solution of the adjoint terminal value problem:   
	$$
	\partial_s y_s = -\big{(} A_{s} \big{)}^*(y_s)   
	\quad \forall\, s \in [0,t), \quad \textnormal{with} \quad 
	y_t = \xi^*,   
	$$
	as asserted in (\ref{adjoint.backward.uniqueness}). 
	Finally, if $\{x_s\}$ is any given classical solution of 
	this terminal value problem, then its time inverse 
	$\{x_{t-s}\}$ is a classical solution of initial value  
	problem (\ref{sectorial.evolution}).  
	Since we also know that $\overline{D((A_{t-s})^*)}=(L^p(\Sigma,\rel^n))^*$ 
	for every $s\in [0,t]$, that $\xi^* \in (L^p(\Sigma,\rel^n))^*$ and that the family of linear operators $\{(A_{t-s})^*\}_{s\in [0,t]}$ 
	in (\ref{sectorial.evolution}) satisfies all requirements of 
	Section 6.1 in \cite{Lunardi}, we may infer here from 
	Corollary 6.2.4 in \cite{Lunardi},
	that the family $\{x_{t-s}\}$ actually has to be the \underline{unique mild solution} of initial value problem (\ref{sectorial.evolution}), 
	and that therefore $x_{t-s}=\tilde y_s$ for every
	$s\in [0,t]$. 
	Hence, together with (\ref{mild.solution}) we arrive at: 
	$x_s = \tilde y_{t-s} = y_s$, for $s\in [0,t]$, proving  
	that such $\{x_s\}$ has to be the concretely 
	given family of functions in (\ref{y.s}), which has shown 
	the ``uniqueness part'' of the assertion in (\ref{adjoint.backward.uniqueness}).
\qed   \\
\end{itemize}
\noindent
Combining Theorems \ref{real.analytic.flow} and 
\ref{Frechensbergo}, we can now prove the final result 
of this article. The application of the strong regularity
result of Theorem \ref{real.analytic.flow} (ii) 
in the proof of the second part of Theorem \ref{Functional.prop.Frechet.derivative} below cannot be 
substituted by Theorem 1.1 of \cite{Malgrange.Lions}, 
whose requirements are not satisfied in our 
mathematical situation. Exactly this technical 
difficulty motivated the author to prove 
Theorem \ref{real.analytic.flow} (ii).        	
\begin{theorem}  \label{Functional.prop.Frechet.derivative}    
	Suppose that $\Sigma$ is a smooth compact torus, being 
	endowed with a complex structure, and suppose that 
	$F_0 \in \textnormal{Imm}_{\textnormal{uf}}(\Sigma,\rel^n)
	\cap C^{\omega}(\Sigma,\rel^n)$
	is a real analytic and umbilic-free immersion.  
	Furthermore, let $T\in (0,T_{\textnormal{max}}(F_0))$ and 
	$p\in (3,\infty)$ be arbitrarily chosen, where $T_{\textnormal{max}}(F_0)$ has been defined in Definition 
	\ref{Umbilic.free}. Then the following statements hold: 
	\begin{itemize}
		\item[1)] The Gateaux derivative $[(x,t) \mapsto (D_F\PP^*(t,0,F_0).\xi)(x)]$ 
		of the evolution operator in line (\ref{solution.operator.1}) 
		is real analytic in time $t\in (0,T)$, for every fixed 
		$\xi \in W^{4,p}(\Sigma,\rel^n)$,
		and the unique classical solution 
		$[(x,s) \mapsto (D_{F}\PP^*(T,s,F_0))^*(\xi^*)(x)]$
		of the adjoint linear terminal value problem (\ref{adjoint.backward.uniqueness}) is real analytic 
		in time $t\in (0,T)$, for every fixed 
		$\xi^* \in (L^{p}(\Sigma,\rel^n))^*$.
		\item[2)] For any fixed $t_0 \in (0,T)$ we denote 
		$F_1:=\PP^*(t_0,0,F_0)$. Then the family of 
		linear differential operators
		$\{A_t\}_{t\in [t_0,T]} \equiv \{D(\Mill_{F_0})(\PP^*(t,t_0,F_1))\}_{t\in [t_0,T]}$ generating the linear evolution equation (\ref{Poincare.trick}) on $\Sigma \times [t_0,T]$ in the second part of Theorem \ref{Frechensbergo}, and also its ``$L^p$-adjoint family'' of linear operators
		$\{-(A_t)^*\}_{t\in [t_0,T]}$ - along the flow line 
		$[t \mapsto \PP^*(t,t_0,F_1)]$ of the modified MIWF (\ref{de_Turck_equation}) - have the ``backwards-uniqueness-property'' in the sense of Section 1 in \cite{Malgrange.Lions}. This means here precisely the following two statements: \\
		a) For any function $\xi \in W^{4,p}(\Sigma,\rel^n)$ the function
		$[(x,t) \mapsto D_{F}\PP^*(t,t_0,F_1).\xi(x)]$
		from the second part of Theorem \ref{Frechensbergo} vanishes identically on $\Sigma \times [t_0,T]$, if there should hold: $D_{F}\PP^*(t,t_0,F_1).\xi \equiv 0$ on $\Sigma$ at time $t=T$, i.e. if $[(x,t) \mapsto D_{F}\PP^*(t,t_0,F_1).\xi(x)]$
		solves the linear terminal value problem:
		\begin{equation}  \label{backward.uniqueness.L}
			\partial_t z_t = D(\Mill_{F_0})(\PP^*(t,t_0,F_1))(z_t),  \,\,\, \textnormal{for}
			\,\,\,  t \in [t_0,T], \,\,\, \textnormal{with} \,\,\, z_T \equiv 0
			\,\,\, \textnormal{on} \,\,\, \Sigma.
		\end{equation}
		b) For every $t \in (t_0,T]$ and every $\xi^* \in (L^{p}(\Sigma,\rel^n))^*$ the unique classical solution  
		$$
		y_s:= \Big{(}D_{F}\PP^*(t,s,\PP^*(s,t_0,F_1)) 
		\Big{)}^*(\xi^*), \,\,\, \textnormal{for} \,\, s \in [t_0,t], 
		$$
		of the adjoint linear terminal value problem:
		\begin{eqnarray} \label{adjoint.backward.uniqueness.2}
			\partial_s y_s = - \Big{(}D(\Mill_{F_0})(\PP^*(s,t_0,F_1))\Big{)}^*(y_s)  
			\,\,\, \textnormal{in} \,\,\, (L^{p}(\Sigma,\rel^n))^*,
			\,\,\, \textnormal{for} \,\, s \in [t_0,t),  \nonumber \\ \textnormal{with} \quad  y_{t} = \xi^*  \,\,\, \textnormal{in}
			\,\,\, (L^{p}(\Sigma,\rel^n))^*  \quad
		\end{eqnarray}
		from the fifth part of Theorem \ref{Frechensbergo}
		vanishes identically on $\Sigma \times [t_0,t]$, if it also satisfies $y_{t_0}=0$ in $(L^{p}(\Sigma,\rel^n))^*$ at time $s=t_0$.
		\item[3)] For any fixed $t_0 \in (0,T)$ and every 
		$t \in [t_0,T]$ the unique linear continuous extension 
		$$
		G^{F_0}(t,t_0): L^{p}(\Sigma,\rel^n) \longrightarrow  L^{p}(\Sigma,\rel^n)
		$$
		of the Fr\'echet derivative $D_{F}\PP^*(t,t_0,F_1)$ 
		from line (\ref{Frechet.evol.operator.s}) has dense range in $L^{p}(\Sigma,\rel^n)$, where we have set again 
		$F_1:=\PP^*(t_0,0,F_0)$.
	\end{itemize}
\end{theorem}
\proof
\begin{itemize}
	\item[1)] 
	From Theorem \ref{real.analytic.flow} (ii) and from 
	the definition of the differential operator $\Mill_{F_0}$ in formula (\ref{de_Turck_equation}) 
	we can immediately deduce the fact, that all coefficients on
	the right hand side of equation (\ref{Poincare.trick}) respectively (\ref{linear.initial.value.problem}) have to be real analytic functions on $\Sigma \times (0,T)$.
	Since we also know from the proof of the second part of Theorem \ref{Frechensbergo}, that the linear operators 
	$A_t = D(\Mill_{F_0})(\PP^*(t,0,F_0))$
	from line (\ref{A.t}) are sectorial in $L^p(\Sigma,\rel^n)$, 
	uniformly for every $t \in [0,T]$, 
	we may apply Theorems 8.1.1 and 8.3.9 in \cite{Lunardi} to the linear parabolic system (\ref{Poincare.trick})
	and infer that its unique strict solution 
	$[(x,t) \mapsto (D_{F}\PP^*(t,0,F_0).\xi)(x)]$,
	starting in some $\xi \in W^{4,p}(\Sigma,\rel^n)$ at time $t=0$ - which is provided by the second part of Theorem \ref{Frechensbergo} - is real analytic in
	$t\in (0,T)$, with values in $W^{4,p}(\Sigma,\rel^n)$.  
	Similarly, we can argue that the unique classical solution 
	\begin{eqnarray*}
	y_{s} := [(x,s) \mapsto (D_{F}\PP^*(T,s,\PP^*(s,0,F_0)))^*(\xi^*)(x)] \\
	\equiv [(x,s) \mapsto (G^{F_0}(T,s))^*(\xi^*)(x)]
	\end{eqnarray*}
	of the adjoint linear terminal value problem (\ref{adjoint.backward.uniqueness}) 
	is real analytic in $s \in (0,T)$, for 
	every fixed $\xi^* \in (L^{p}(\Sigma,\rel^n))^*$.
	For, we know from the proof of the fifth part of Theorem 
	\ref{Frechensbergo}, that the linear operators $A(t)^*$ from (\ref{domain.A.star}) satisfy all requirements of Section 6.1 in \cite{Lunardi}, and additionally 
	that they depend real analytically on $t\in (0,T)$, just as the operators $A(t)$ do here. We can therefore argue again
	by means of Theorems 8.1.1 and 8.3.9 in \cite{Lunardi}, that 
	the unique classical solution 
	$\{\tilde y_s\}=\{y_{T-s}\}$ of equation (\ref{adjoint.backward.unique}), 
	starting in $\tilde y_0=\xi^* \in (L^p(\Sigma,\rel^n))^*$ and of the regularity stated in line (\ref{regularity.backwards.2}), 
	is real analytic in $s\in (0,T)$, with values in $D(A_0^*)$. 
	This implies immediately that also the unique classical solution 
	$\{y_s\} \equiv (D_{F}\PP^*(T,s,\PP^*(s,0,F_0)))^*(\xi^*)=
	\{\tilde y_{T-s}\}$ of problem (\ref{adjoint.backward.uniqueness}) 
	is real analytic in $s\in (0,T)$, with 
	values in $D(A_0^*)$, for any 
	$\xi^* \in (L^{p}(\Sigma,\rel^n))^*$.  
	\item[2)] a) Suppose that the unique strict solution 
	$[(x,t) \mapsto D_{F}\PP^*(t,t_0,F_1).\xi(x)]$
	of the linear system (\ref{linear.initial.value.problem}) 
	would also satisfy $z^*_{T} \equiv 0$ on $\Sigma$, at some 
	fixed time $T\in (0,T_{\textnormal{max}}(F_0))$ and 
	for some $\xi \in W^{4,p}(\Sigma,\rel^n)$. 
	Since the flow line $\{\PP^*(\,\cdot\,,0,F_0)\}$ of the modified MIWF (\ref{de_Turck_equation}) exists smoothly in the time interval $(0,T_{\textnormal{max}}(F_0))$, 
	and since the statements of Theorem \ref{real.analytic.flow} 
	(ii) and of the first part of this theorem hold in every fixed time interval $(0,T)$, for $T \in (0,T_{\textnormal{max}}(F_0))$, we may vary here our choice of starting and final time, namely setting $t_0':=T$ and $T':=T+\varepsilon$, for any sufficiently small $\varepsilon>0$. Using the uniqueness of the strict
	solution to the initial value problem (\ref{Poincare.trick}) w.r.t. the spaces $X:=L^p(\Sigma,\rel^n)$ and  
	$D:=W^{4,p}(\Sigma,\rel^n)$ by the second part of Theorem \ref{Frechensbergo}, the unique short-time extension, again denoted by $\{z_t^*\}$, of the solution $\{z^*_t\}:= [(x,t) \mapsto D_{F}\PP^*(t,t_0,F_1).\xi(x)]$ of equation (\ref{Poincare.trick}) from the interval $[t_0,T]$ to the larger interval $[t_0,T+\varepsilon]$ is given by the unique short-time solution of the initial value problem (\ref{Poincare.trick}),
	which starts at time $t=T$ in the function $z^*_{T} \in W^{4,p}(\Sigma,\rel^n)$.
	Since this short-time solution is explicitly given by the function
	$$
	[(t,x) \mapsto D_{F}\PP^*(t,T,\PP^*(T,t_0,F_1)).(z^*_T)(x)]  \quad \textnormal{for} \quad (x,t) \in \Sigma \times [T,T+\varepsilon]
	$$
	on account of equation (\ref{Frechet.G}),
	we conclude from $z^*_{T} \equiv 0$ on $\Sigma$
	that $z^*_{t} \equiv 0$ on $\Sigma \times [T,T+\epsilon]$.
	Hence, on account of the real analytic dependence  
	of the function $[t \mapsto z^*_t]$ on the time 
	$t \in [t_0,T+\epsilon)$ by the first part of 
	this theorem, we can therefore conclude that 
	$z^*_t \equiv 0$ on $\Sigma \times [t_0,T+\epsilon]$.
	This proves already the first assertion of the second part of the theorem. \\
	b) Now we consider the unique classical solution 
	of the adjoint linear terminal value problem 
	(\ref{adjoint.backward.uniqueness.2}), provided by the fifth part of Theorem \ref{Frechensbergo}. The unique classical extension $\{z^*_s\}$ of the mild and classical solution 
	\begin{eqnarray*}
	[(x,s) \mapsto (D_{F}\PP^*(t,s,\PP^*(s,t_0,F_1)))^*(\xi^*)(x)] \\
	\equiv                    
	[(x,s) \mapsto (G^{F_0}(t,s))^*(\xi^*)(x)] \quad
	\textnormal{for} \quad (x,s) \in \Sigma \times [t_0,t]
	\end{eqnarray*}
	of the terminal value problem (\ref{adjoint.backward.uniqueness}) 
	from the interval $[t_0,t]$ to the larger interval $(0,t]$ is given by the unique classical solution of the adjoint terminal value problem (\ref{adjoint.backward.uniqueness}), which stops at time $t=t_0$ in the function $z^*_{t_0} \in (L^{p}(\Sigma,\rel^n))^*$.
	Since this classical solution is explicitly given by the function
	\begin{eqnarray*} 
		[(x,s) \mapsto  (D_{F}\PP^*(t_0,s,\PP^*(s,0,F_0)))^*
		(z^*_{t_0})(x)]           \\
		\equiv                                         
		[(x,s) \mapsto (G^{F_0}(t_0,s))^*(z^*_{t_0})(x)] \quad 
		\textnormal{for} \,\,\, (x,s) \in \Sigma \times (0,t_0]
	\end{eqnarray*} 
	on account of the fifth part of Theorem \ref{Frechensbergo},
	we conclude from the assumption ``$z^*_{t_0} \equiv 0$ on $\Sigma$'', that $z^*_{s} \equiv 0$ on $\Sigma \times (0,t_0]$.
	Hence, on account of the real analytic dependence  
	of the function $[s \mapsto z^*_s]$ on the time 
	$s \in (0,t)$ by the first part of this theorem, we can 
	therefore conclude that the original 
	solution $\{z^*_{s}\}_{s\in [t_0,t]}$ of 
	the adjoint linear terminal value problem (\ref{adjoint.backward.uniqueness.2}) has to vanish identically on $\Sigma \times [t_0,t]$, which proves also the second assertion of the second part of this theorem. 
	\item[3)] We fix some $t \in (t_0,T]$ arbitrarily, 
	and we consider some arbitrary function 
	$\eta \in (L^p(\Sigma,\rel^n))^* \cong L^q(\Sigma,\rel^n)$, $q = \frac{p}{p-1}$, 
	which satisfies:   
	\begin{eqnarray} \label{contradiction}
		\int_{\Sigma}  \langle \eta, 
		D_{F}\PP^*(t,t_0,F_1).(\xi) \rangle_{\rel^n} \,d\mu_{\PP^*(t,t_0,F_1)} = 0  \quad \forall \,\xi \in W^{4,p}(\Sigma,\rel^n). \qquad 
	\end{eqnarray}
	We try to prove in the sequel, that any such function $\eta$ has to vanish on $\Sigma$. To this end, we recall from the fourth part of Theorem \ref{Frechensbergo}, that the family of solutions of the homogeneous, linear system (\ref{Poincare.trick}) respectively (\ref{linear.initial.value.problem}) constitutes a parabolic fundamental solution $\{G^{F_0}(t,s)\}_{s\leq t}$ in $L^p(\Sigma,\rel^n)$, 
	which coincides with the family of Fr\'echet derivatives 
	$\{D_{F}\PP^*(t,s,\PP^*(s,0,F_0))\}_{s\leq t}$ on $W^{4,p}(\Sigma,\rel^n)$,
	and we recall from the proof of the fifth part of Theorem \ref{Frechensbergo}, that for every pair of times $t_2 \geq t_1$ in $[t_0,T]$ the operators
	\begin{equation}  \label{new.starting.time.2}
	D_{F}\PP^*(t_2,t_1,F_1) := D_{F}\PP^*(t_2,t_1,\PP^*(t_1,t_0,F_1))
	\end{equation}
	from line (\ref{new.starting.time}) - up to exchanging $F_1:=\PP^*(t_0,0,F_0)$ with $F_0$ - can be considered as continuous linear operators 
	\begin{equation}  \label{Frechet.L.p.2}
	D_{F}\PP^*(t_2,t_1,F_1): 
	L^{p}(\Sigma,\rel^n) \longrightarrow  L^{p}(\Sigma,\rel^n),
	\end{equation}
	which have uniquely determined maximal, densely defined 
	and closed linear adjoint operators in $(L^{p}(\Sigma,\rel^n))^*$:
	\begin{equation} \label{adjoint.2}
		(D_{F}\PP^*(t_2,t_1,F_1))^*:
		(L^{p}(\Sigma,\rel^n))^* \longrightarrow (L^{p}(\Sigma,\rel^n))^* 
		\cong L^{q}(\Sigma,\rel^n),
	\end{equation}
	see statement (\ref{adjoint}).
	We use the adjoint operators in line (\ref{adjoint.2})
	and the test function $\eta$ in line (\ref{contradiction}),
	in order to define the family of functions
	$\{y_s\} \in C^{0}([t_0,t];(L^{p}(\Sigma,\rel^n))^*)$ by
	$$
	y_s := \Big{(} D_{F}\PP^*(t,s,F_1) \Big{)}^*(\eta),
	\quad \textnormal{for} \,\,\,s \in [t_0,t].
	$$
	We can immediately infer from the fifth part of Theorem 
	\ref{Frechensbergo}, that we have here
	\begin{equation}  \label{regularity.backwards.3}
		\{y_s\} \in C^{1}([t_0,t);(L^{p}(\Sigma,\rel^n))^*) \cap C^{0}([t_0,t);D(A_0^*)) \cap C^{0}([t_0,t];(L^{p}(\Sigma,\rel^n))^*)
	\end{equation} 
	and that $\{y_s\}$ is the unique classical solution of the linear equation
	\begin{equation} \label{adjoint.evolution.2}
	\partial_s y_s = - \big{(} A_{s} \big{)}^*(y_s),   \quad \forall\, s \in [t_0,t),
	\end{equation}
	satisfying the terminal condition $y_t=\eta$ in $(L^{p}(\Sigma,\rel^n))^*$, 
	which is just the adjoint linear terminal value problem (\ref{adjoint.backward.uniqueness.2}), with $\xi^*$ replaced by $\eta$. Moreover, on account of the definition of the family $\{y_s\}$ and by assumption (\ref{contradiction}) we have:
	\begin{eqnarray} \label{first.strike}
		\int_{\Sigma} \langle y_{t_0}, \xi \rangle \, d\mu_{F_1}
		= \int_{\Sigma}  \langle \eta,
		D_{F}\PP^*(t,t_0,F_1).(\xi)  \rangle \,
		d\mu_{\PP^*(t,t_0,F_1)} = 0  \quad \\ 
		\forall \, \xi \in W^{4,p}(\Sigma,\rel^n).  \quad \nonumber
	\end{eqnarray}
	Since we have $y_{t_0} \in (L^{p}(\Sigma,\rel^n))^* \cong L^{q}(\Sigma,\rel^n)$ by construction, equation (\ref{first.strike}) implies that 
	$y_{t_0} = 0$ in $(L^{p}(\Sigma,\rel^n))^*$.
	Combining this insight with the fact that 
	$\{y_s\}$ is a classical solution of the adjoint linear terminal value problem (\ref{adjoint.backward.uniqueness.2}) of regularity class in line (\ref{regularity.backwards.3}), 
	we can apply the second uniqueness result 
	of the second part of this theorem, guaranteeing that there holds $y_s \equiv 0$ \,on $\Sigma$ for every $s \in [t_0,t]$. Combining this again with the definition of the family of functions $\{y_s\}$, and with the fact that $\eta \in (L^p(\Sigma,\rel^n))^*\cong 
	L^q(\Sigma,\rel^n)$, we can now conclude that
	\begin{eqnarray} \label{weak.convergence.2}
		0 \equiv \int_{\Sigma} \langle y_t, \xi \rangle_{\rel^n} \, d\mu_{\PP^*(t,t_0,F_1)}
		= \int_{\Sigma} \langle \eta, D_{F}\PP^*(t,t,F_1).(\xi) 
		\rangle_{\rel^n} \, d\mu_{\PP^*(t,t_0,F_1)} \nonumber \\
		=\int_{\Sigma} \langle \eta, \xi \rangle_{\rel^n} \,
		d\mu_{\PP^*(t,t_0,F_1)}    \qquad \qquad
	\end{eqnarray}
	for an arbitrary function $\xi \in W^{4,p}(\Sigma,\rel^n)$.
	Hence, there has to hold $\eta = 0$ in $L^q(\Sigma,\rel^n)$ by (\ref{weak.convergence.2}). 
	Since this conclusion holds for every $\eta \in L^q(\Sigma,\rel^n)$ satisfying equation (\ref{contradiction}), we have proved that the linear operator 
	\begin{eqnarray*}
		G^{F_0}(t,t_0) \equiv D_{F}\PP^*(t,t_0,F_1):
		L^{p}(\Sigma,\rel^n) \longrightarrow W^{4,p}(\Sigma,\rel^n),
	\end{eqnarray*}
	appearing in (\ref{contradiction}), must have dense range in $L^{p}(\Sigma,\rel^n)$, which proves the assertion of the last part of the theorem in the considered case of a fixed time $t\in (t_0,T]$. In the special case ``$t=t_0$'' this assertion is trivially true. 
\end{itemize}
\qed 
\noindent



\end{document}